\documentclass{article}
\usepackage[left=1in,right=1in,top=1in,bottom=1in]{geometry}
\usepackage{amsmath}
\usepackage{amsthm}
\usepackage{amsfonts}
\usepackage{amssymb}

\title{Phase Space Analysis on some Black Hole Manifolds}
\author{P. Blue, A. Soffer}
\date{\today}

\numberwithin{equation}{section}

\newcommand{\hide}[1]{}

\newtheorem{theorem}{Theorem}[section]
\newtheorem{corollary}[theorem]{Corollary}
\newtheorem{lemma}[theorem]{Lemma}
\newtheorem{proposition}[theorem]{Proposition}
\newtheorem{definition}[theorem]{Definition}

\newtheorem{notation}[theorem]{Notation}
\newtheorem{faketheoremaux}{Theorem}
\newenvironment{faketheorem}[1]{\renewcommand{\thefaketheoremaux}{\ref{#1}}\begin{faketheoremaux}}{\end{faketheoremaux}}
\newtheorem{fakecorollaryaux}{Corollary}

\newcommand{\Reals}{\mathbb{R}}

\newcommand{\Integers}{\mathbb{Z}}
\newcommand{\Naturals}{\mathbb{N}}
\newcommand{\Fourier}[1]{\mathcal{F}[#1]}

\newcommand{\Dt}{\frac{d}{dt}}
\newcommand{\dt}{\frac{\partial}{\partial t}}
\newcommand{\dr}{\frac{\partial}{\partial \rho_*}}
\newcommand{\drr}{\frac{\partial^2}{\partial r_*^2}}

\newcommand{\dalembertian}{\Box}
\newcommand{\defin}{\equiv}

\newcommand{\starman}{\mathfrak{M}}
\newcommand{\solset}{\mathbb{S}}

\begin{document}


\newcommand{\RN}{Reissner-Nordstr{\o}m }
\newcommand{\dAlembertianRN}{\square_{\text{RN}}}
\newcommand{\SLap}{\Delta_{S^2}}

\newcommand{\linH}{H}

\newcommand{\uasoluinS}{ $u\in\solset$ is a solution to the wave equation, $\ddot{u}+\linH u=0$}
\newcommand{\foruasoluinS}{ for $u\in\solset$ a solution to the wave equation, $\ddot{u}+\linH u=0$}

\newcommand{\BigOOne}{\text{O($1$)}}

\newcommand{\Lu}{\nabla_{S^2}u}
\newcommand{\Lv}{\nabla_{S^2}v}
\newcommand{\dmu}{d^3\mu}

\newcommand{\drs}{\partial_{\rho_*}}

\newcommand{\EDens}{e}
\newcommand{\PsDens}{p_{\rho_*}}
\newcommand{\PlDens}{\vec{p}_{\omega}}

\newcommand{\confmult}{\mathcal{C}}
\newcommand{\confdens}{e_\mathcal{C}}
\newcommand{\confchrg}{E_\mathcal{C}}
\newcommand{\confdenstr}{e_{\mathcal{C},(t,\rho_*)}}

\newcommand{\retext}{{\bf Rewrite this text!}}



\newcommand{\Hplus}{\mathcal{H}_+}
\newcommand{\Hminus}{\mathcal{H}_-}

\newcommand{\dFromE}{Fr^2d\rho_* d^2\omega}
\newcommand{\dTrueExpanded}{F^\frac12 r^2 d\rho_* d^2\omega}
\newcommand{\dTrue}{d\mu_{\text{normalised}}}
\newcommand{\FourGrad}{\tilde{\nabla}_4}

\newcommand{\xu}{s_-}
\newcommand{\xv}{s_+}
\newcommand{\du}{\partial_{s_-}}
\newcommand{\dv}{\partial_{s_+}}

\newcommand{\Conf}{\mathcal{C}_{needs fixin'}}

\newcommand{\id}{\text{Id}}

\newcommand{\verbose}[1]{\footnote{ #1} }

\newcommand{\PQfakeroot}{\beta_{\text{approx}}}

\newcommand{\Pl}{ \mathbf{P}_l}
\newcommand{\weakloc}{ \left(1+\rho_*^2\right)}
\newcommand{\chia}{\chi_\alpha}
\newcommand{\supp}{\text{supp}}
\newcommand{\linHexpanded}{\sum_{i=1}^3 H_i}


\newcommand{\qm}{{\bf x}_m}
\newcommand{\pn}{{ {\boldsymbol \xi}_n }} 
\newcommand{\phalf}{{\boldsymbol \xi}_\frac12}
\newcommand{\qn}{{\bf x}_n}
\newcommand{\qhalf}{{\bf x}_\frac12}

\newcommand{\qmnfn}{X}
\newcommand{\qmn}[1]{\qmnfn_\downarrow(#1)}
\newcommand{\qmnWITHDERIV}[1]{\qmnfn_\downarrow'(#1)}
\newcommand{\qmnEXPANDED}[1]{\left( 1 + \left(\frac{#1}{2M}\right)^2\right)^{-1}}
\newcommand{\qmf}[1]{\qmnfn_\uparrow(#1)}

\newcommand{\Oint}{O(L^1)}
\newcommand{\ErrorTermsn}{O(\|L^{\frac{1+2n}{2}}\chia u\|^2,L^1)}
\newcommand{\bndd}{B}

\newcommand{\pnnfn}{\Phi}
\newcommand{\pnnaNOARG}{\pnnfn_{a,\epsilon}}
\newcommand{\pnna}[1]{\pnnaNOARG(#1)}
\newcommand{\pnnaWITHADERIV}[1]{\pnnaNOARG'(#1)}
\newcommand{\pnnaWITHTWODERIV}[1]{\pnnaNOARG''(#1)}
\newcommand{\pnnbNOARG}{\pnnfn_{b,\epsilon}}
\newcommand{\pnnb}[1]{\pnnbNOARG(#1)}
\newcommand{\pnncNOARG}{\pnnfn_{c,\epsilon}}
\newcommand{\pnnc}[1]{\pnncNOARG(#1)}
\newcommand{\pnn}[1]{{\pnnfn}_{|#1|\leq1}}
\newcommand{\pnnNOARG}{\pnn{\bullet}}

\newcommand{\pnifn}{\Psi}
\newcommand{\pniaNOARG}{\pnifn}
\newcommand{\pnia}[1]{\pniaNOARG(L^{-\delta},#1)}
\newcommand{\pniaWITHl}[1]{\pniaNOARG(l^{-\delta},#1)}
\newcommand{\pnibNOARG}{\pnifn_{1}}
\newcommand{\pnib}[1]{\pnibNOARG(L^\delta #1)}
\newcommand{\pnibWITHl}[1]{\pnibNOARG(l^\delta #1)}
\newcommand{\pnibDERIVWITHl}[1]{\pnibNOARG'(l^\delta #1)}
\newcommand{\pnibDERIVNOARGNOSCALE}{\pnibNOARG'}
\newcommand{\pnicNOARG}{\pnifn_{2}}
\newcommand{\pnic}[1]{\pnicNOARG(#1)}

\newcommand{\pni}[1]{{\pnifn}_{L^{-\delta}\leq|#1|\leq1}}
\newcommand{\pniNOARG}{\pni{\bullet}}
\newcommand{\pniWITHl}[1]{{\pnifn}_{l^{-\delta}\leq|#1|\leq1}}
\newcommand{\pniGENERICINDEX}[1]{{\pnifn}_{i}(#1)}

\newcommand{\pniaOverF}[1]{(\pnia{\bullet} F_2^{-1})(#1)}
\newcommand{\pnicOverF}[1]{(\pnicNOARG F_2^{-1})(#1)}
\newcommand{\pnicOverFNOARG}{(\pnicNOARG F_2^{-1})}

\newcommand{\pniaOverpnna}[1]{(\pnia{\bullet} \pnnaNOARG^{-1})(#1)}
\newcommand{\pniaOverpnnc}[1]{(\pnia{\bullet} \pnncNOARG^{-1})(#1)}
\newcommand{\pniaOverpnncWITHl}[1]{(\pniaWITHl{\bullet} \pnncNOARG^{-1})(#1)}

\newcommand{\Gammanm}{\Gamma_{n,m}}
\newcommand{\Gammanhalf}{\Gamma_{n,\frac12}}
\newcommand{\Gammann}{\Gamma_{n,n}}

\newcommand{\ad}{\text{Ad}}
\newcommand{\opnorm}{{}}




\newcommand{\tildeutruesolu}{ $\tilde{u}$ is a solution to the wave equation \eqref{tildeLW} }
\newcommand{\tildestarman}{\tilde{\starman}}
\newcommand{\tildedmu}{d^3\tilde{\mu}}



\newcommand{\sgn}{\text{sgn}}
\newcommand{\qmnt}[1]{\tilde\qmnfn_\downarrow(#1)}
\newcommand{\qmntWITHDERIV}[1]{\tilde\qmnfn_\downarrow'(#1)}

\newcommand{\tl}{\tilde{l}}
\newcommand{\Vtl}{V_l}
\newcommand{\Ptl}{{\bf P}_{l}}
\newcommand{\al}{(\alpha_l)_*}

\newcommand{\PolyZero}{P}
\newcommand{\PolyL}{P_L}
\newcommand{\Polyl}{P_{l}}

\newcommand{\I}{I}
\newcommand{\IL}{I_L}
\newcommand{\Il}{I_{\tl}}

\newcommand{\threeprterms}{L^{5m}}

\maketitle

\begin{abstract}

The Schwarzschild and \RN solutions to Einstein's equations describe space- times which contain spherically symmetric black holes. We consider solutions to the linear wave equation in the exterior of a fixed black hole space- time of this type. We show that for solutions with initial data which decay at infinity and at the bifurcation sphere, a weighted $L^6$ norm in space decays like $t^{-\frac{1}{3}}$. This weight vanishes at the event horizon, but not at infinite. 
\end{abstract}

\tableofcontents

\section{Introduction}
\label{sIntroduction}

Black holes are very important objects in Relativity, but very little is known about their dynamics and interaction. Even the question of stability under small perturbations is a challenging, open problem. Any solution to such questions will require an understanding of the interaction between black holes and gravitational radiation. The structure of the vacuum Einstein equations, which govern these dynamics in the absence of other matter, make it impossible to consider only radial perturbations, so that more general perturbations must be considered immediately. We hope that the study of linear, uncoupled waves outside the black hole will help provide an understanding of these problems. 

Even the stability of the empty space was a very difficult problem. It was solved originally by Christodoulou and Klainerman \cite{ChristodoulouKlainerman}, and Lindblad and Rodnianski have since developed a simpler proof \cite{LindbladRodnianski}. Decay estimates for solutions of non linear wave equations played an important role in both proofs.

In Relativity, the structure of space- time is determined by a Lorentzian pseudo- metric, $g$, which, in the absence of matter, satisfies the Einstein equations, 
\begin{align}
R_{\mu\nu} - \frac12 R g_{\mu\nu} = 0. 
\end{align}

The Schwarzschild and \RN solutions are singular solutions to the Einstein equations which describe space- times containing a spherically symmetric black hole. The Schwarzschild solution is a special case of the \RN solutions. We will restrict our attention to the exterior region, outside the black hole. We discuss the geometry of the exterior region further in subsection \ref{ssRNsolu}. Because the structure of the Einstein equations makes it impossible for these to have spherically symmetric perturbations, it is expected that the study of black- hole stability will require analysis of the more general, non spherically symmetric, rotating, Kerr-Newman black holes. The linearisation of the Einstein equations around the Kerr solution has been shown to have no unstable modes\cite{Whiting}. 

In 1957, Regge and Wheeler investigated the linear stability of the Schwarzschild black- hole, and were able to reduce the problem to the study of a second order, scalar equation\cite{ReggeWheeler}. After appropriate transformations, this equation and the geometrically defined, scalar wave equation differ only by a multiple of the potential term appearing in each. Using these transformations, the exterior region of the black- hole can be decomposed into the product of time and a three dimensional space. 

Because of the intrinsic interest in the geometric wave equation and because of its possible applications to the study of black hole stability, we consider the wave equation
\begin{align}
\dalembertian \tilde{u} =&0, \\
\tilde{u}(1)=&\tilde{u}_0,\\
 \dot{\tilde{u}}(1)=&\tilde{u}_1  . 
\end{align}
For this equation, we prove the following result. 

\begin{faketheorem}{cL6Control}
\label{ftL6ControlIntro}
If \tildeutruesolu, and $u=r\tilde{u}$, 
\begin{align}
\| F^\frac12 \tilde{u} \|_{L^6(\tildestarman)} 
\leq& t^{-\frac{1}{3}} C(\|u_0\|_2^2 + E[u_0,u_1] + \confchrg[u_0,u_1] + \|L^\epsilon u\|^2 + E[L^\epsilon u_0, L^\epsilon u_1])^\frac12  ,\\
\| (\rho_*^2+1)^\frac{-1}{2} \tilde{u} \|_{L^2(\tildestarman)}
\leq& t^{-\frac{1}{2}} C(\|u_0\|_2^2 + E[u_0,u_1] + \confchrg[u_0,u_1] + \|L^\epsilon u\|^2 + E[L^\epsilon u_0, L^\epsilon u_1])^\frac12 ,
\end{align}
where $E[v,w]$ and $\confchrg[v,w]$ are the energy and conformal charge which are defined in section \ref{sMethodsPD}, $\rho_*$ is the Regge-Wheeler radial co-ordinate defined in section \ref{sWaveEqn}, and $L$ is an angular derivative operator defined in section \ref{sAngularModulation}. The measures are defined in section \ref{sWaveEqn}. 
\end{faketheorem}

Previously, we were able to apply our techniques to prove results for both the wave equation on the Schwarzschild manifold and the Regge-Wheeler equation, and expect the same to be true here\cite{BlueSoffer,BlueSofferII}. 

This decay rate is slower than the rate found in Euclidean space, which is $t^\frac{-2}{3}$\cite{GinibreVelo}. This is consistent with current results for other equations on the Schwarzschild solution, which only prove slower rates of decay than those found in $\Reals^{3+1}$. The $L^\infty$ norm has been shown to decay like $t^{\frac{-5}{6}}$ for both the massive Dirac equation \cite{FKSY} and the massive Klein- Gordon equation \cite{KoyamaTomimatsu}. For the Schr\"odinger equation, the $L^\infty$ has been shown to decay like $t^{\frac{-1}{4}+\epsilon}$\cite{LabaSoffer}. All these results hold only for radial initial data, or data which is the sum of finitely many spherical harmonics. 

The rigorous study of scattering for linear waves on black hole backgrounds was begun by Dimock and Kay who proved existence and completeness of wave operators for the linear wave equation \cite{Dimock} and the Klein-Gordon equation\cite{DimockKay}. Asymptotic completeness and global existence for the cubic semi linear wave equation have been proven for a wide class of black hole manifolds, including the Kerr- Newman solutions\cite{Bachelot,BachelotNicolas,Nicolas}. DeBi\`evre, Hislop, and Sigal have proven asymptotic completeness on a large class of non compact manifolds, including the \RN manifolds \cite{DeBievreHislopSigal}. However, none of these methods give decay estimates even in the linear case. For sufficiently super critical \RN black holes, Strichartz estimates have also been proven\cite{BPST}. Little is known about fields coupled to the Einstein equations, but for a spherically symmetric, scalar field, Dafermos and Rodnianski have shown the field decays along the event horizon according to an inverse cubic, Price law \cite{DafermosRodnianski}. They also proved a weaker rate of decay for decoupled, radial, semi linear wave equations in the exterior of \RN solutions \cite{DafermosRodnianskiII}. 

We start by introducing an analogue of the conformal charge used by Ginibre and Velo \cite{GinibreVelo}. However, we are prevented from completing the argument which holds for the wave equation in Euclidean space by the presence of a closed geodesic surface in the three dimensional space which describes the exterior region of the black hole. The absence of such geodesics is a non- trapping condition which is commonly imposed in scattering theory. The presence of this geodesic surface and the gravitational lensing they cause is already known\cite{EllisPhotonSphere,EllisLensing}. This surface is called the photon sphere. 

To over come this obstacle, it is sufficient to prove estimates on the angular derivative of the solution in the region near the closed geodesics. We prove:
\begin{faketheorem}{tPhaseSpaceInduction}
\label{ftPhaseSpaceInductionIntro}
If $\varepsilon>0$, then if \tildeutruesolu and $u=r\tilde{u}$
\begin{align}
\int_1^\infty \|L^{1-\frac\varepsilon2} \chia u\|^2 dt < C(\|u_0\|_{L^2}^2 + E[u_0,u_1]) ,
\end{align}
where $L$ acts like one angular derivative and $\chia$ is a function with compact support near the closed geodesics. 
\end{faketheorem}
This is sufficient to prove theorem \ref{ftL6ControlIntro}. The loss of $L^\epsilon$ is responsible for the additional factors of $L^\epsilon$ appearing in theorem \ref{ftL6ControlIntro}. 

Estimates on the angular derivative have already been used to prove Strichartz and point wise in time $L^p$ estimates in Euclidean space\cite{Sterbenz,MachiharaNakamuraNakanishiOzawa} and on non-trapping manifolds \cite{HassellTaoWunsch}. 

Our method is to introduce an analogue for the wave equation of the Heisenberg equation for the Schr\"odinger equation. There is a self-adjoint operator $\linH$ which determines the time evolution of the solution. We refer to this operator as the Hamiltonian. Other self-adjoint operators are referred to as observables, and the analogue of the Heisenberg identity relates the time derivative of a particular inner product involving the observable to the expectation value of the commutator between the Hamiltonian and the observable. For the Schr\"odinger equation, the expectation value of the observable is differentiated in time, but for the wave equation, the analogue is a more complicated inner product involving both a solution and its time derivative. For an operator to be a propagation observable, we require a few other conditions. 

Our first propagation observable is a radial derivative operator directed away from the geodesic surface. This is used to prove a smoothed Morawetz estimate, which is like theorem \ref{ftPhaseSpaceInductionIntro}, but with zero powers of the angular derivative and a non compactly supported, but decaying, localisation function. A function, $g$, of the radial variable, denoted $\rho_*$, is used to direct this operator away from the geodesics. 

Following this, we define the positive operator $L$ by 
\begin{align}
L^2 = 1-\Delta_{S^2} ,
\end{align}
where $\Delta_{S^2}$ is the Laplace- Beltrami operator on the sphere, which acts as the derivative in the angular directions. By rescaling the argument of $g$ by powers of $L$, we generate a new propagation observable, which we use to prove theorem \ref{ftPhaseSpaceInductionIntro} with $\varepsilon=\frac12$. We call this rescaling by $L$ angular modulation. 

To prove the result for all positive $\varepsilon$, we use phase space analysis. We introduce the phase space variables, which are defined in terms of the original radial co-ordinate, $\rho_*$, and radial derivative, $\dr$. 
\begin{align}
\qm=& L^m \rho_* , \\
\pn=& L^{n-1} \dr .
\end{align}
We localise a propagation observable in phase space by multiplying by a compactly supported function of the phase space variable, $\qmnfn(\qm)$ and $\pnnfn(\pn)$. We also localise with functions which are not compactly supported but which decay away from a certain region. By both rescaling the argument of $g$ and localising in $\pn$, we can control powers of $L$ in various regions of phase space. However, these estimates require control of different powers of $L$ in other regions of phase space. In a process of phase space induction, we are able to combine the estimates across all regions, and use the estimate in one region to control the remainder terms in another. This is sufficient to prove theorems \ref{ftPhaseSpaceInductionIntro} and \ref{ftL6ControlIntro}.

%
%
%

\subsection{Structure of the paper}
\label{ssStructure}

Section \ref{sWaveEqn} introduces the wave equation on the \RN solution, some common notation and transformations to simplify the equation, and discusses important regions of the geometry, including the photon sphere. Section \ref{sMethodsPD} introduces the energy, conformal charge, and a Sobolev estimate. In that section, we prove pointwise in time $L^p$ estimates can be reduced to bounding weighted space-time integrals of solutions and their angular derivatives. Section \ref{sRelConsid} translates some of our conditions into terminology commonly used by Physicists. In particular, we show that the finite $L^2$, energy, and conformal charge conditions, require that initial date, but not its derivatives, vanish at the bifurcation sphere. In fact, functions which grow like an arbitrarily small power of the distance from the bifurcation sphere are permitted. 

The remainder of the paper deals with proving weighted space- time integral bounds. Section \ref{sHeisenbergType} provides the basic propagation observable frame work, using a commutator formalism analogous to the Heisenberg equation from quantum mechanics. Section \ref{sMorawetz} introduces a propagation observable $\gamma$ which majorates a decaying weight. This is analogous to a smoothed Morawetz estimate. This is sufficient to control the weighted space- time integral needed in the proof of the $L^p$ estimates. The section concludes with an $L^p$ estimate for radial data and with a similar result for non-radial data, which reduces the problem to one of controlling space- time integrals of one angular derivative in $L^2$. This localisation is in an arbitrarily small neighbourhood of the photon sphere. Section \ref{sAngularModulation} uses a new technique of angular modulation to control $\frac34$ angular derivatives in $L^2$ in the desired region. The method is to rescale the previous propagation observable by fractional powers of the angular derivative operator. 

In the final section, section \ref{sPhaseSpaceAnalysis}, we introduce a family of propagation observables to control $1-\epsilon$ angular derivatives. The propagation observables are like the one from the angular modulation argument, but are also localised in both the radial variable and radial derivative. We refer to this as phase space analysis. Subsection \ref{ssPhaseSpaceVariablesEtc} introduces these observables and rescaled versions of the radial variable and radial derivative, which we call the phase space variables. Subsection \ref{ssCommutatorExpansions} provides a commutator theorem which allows us to rearrange localisations in the non commuting phase space variables. It also provides lemmas which cover particularly common cases. Subsection \ref{ssPhaseSpaceLocalisedEstimates} begins the argument by using the commutators to bound fractional powers of the angular derivative. These estimates have decaying weights in the rescaled radial variable and derivative. Subsection \ref{ssDerivativeLocalisedEstimates} removes the weights in the rescaled radial variable, leaving estimates localised in the radial derivatives only. All these estimates involve remainder terms which include lower powers of the angular derivative, but not localised in the same region. Subsection \ref{ssPhaseSpaceInduction} combines the results inductively to control $1-\epsilon$ angular derivatives with out phase space localisation. It concludes with an $L^p$ estimate with $\epsilon$ loss of angular derivatives in the energy. 

Through out this paper, the notation $C$ is used to denote constants which may vary from equation to equation. Indices are used to separate different constants within a line of an equation and to refer to constants later in a proof.

%
%
\section{The wave equation on the \RN space}
\label{sWaveEqn}

%
%
\subsection{The \RN solution}
\label{ssRNsolu}

The \RN solution to the Einstein equations represents the space-time outside a spherically symmetric, charged, massive body. It is the unique spherically symmetric, static\footnote{The static condition is redundant since spherically symmetric, vacuum solutions are necessarily static. }, asymptotically flat solution to the vacuum Einstein- Maxwell equations, which govern the structure of space-time in the presence of gravity and electro- magnetic fields. This solution can be represented in terms of the co-ordinates $(t, r, \theta, \phi)$ by the Lorentzian pseudo-metric:
\begin{align}
ds^2 =& F dt^2 - F^{-1} dr^2 -r^2 (d\theta^2 +\sin^2(\theta) d\phi^2) , \\
F =& 1-\frac{2M}{r} + \frac{Q^2}{r^2} ,
\end{align}
where $M\geq0$ is the central mass and $Q\in\Reals$ is the central charge, both as measured by observers at infinity. 

The polynomial $r^2F$ has two roots at 
\begin{align}
r_\pm = M \pm \sqrt{M^2-Q^2} .
\end{align}
In the sub critical case, $0<|Q|<M$, these roots determine the location of two event horizons which permit geodesics and energy to cross only in the inward direction, towards decreasing $r$, and not in the outward direction. In the critical case, $|Q|=M$, there is a single event horizon, which has the same property. In the super critical case, $|Q|>M$, the roots are complex, there are no event horizons to prevent geodesics which start at the singularity $r=0$ from being extended to $r=\infty$, and the central singularity is called a naked singularity. We will not study the super critical case. In the non super critical cases, the effects of the event horizon and the structure of the interior region are complicated and have been studied extensively \cite{EllisHawking, Wald}. The interior region is called a black hole since anything, including a light ray, may fall in but can not escape. 

We will restrict our attention to the exterior region of the non super critical solution, 
\begin{align}
t\in&\Reals,& r>r_+=&M+\sqrt{M^2-Q^2}, & (\theta, \phi)\in& S^2. 
\end{align}
The famous Schwarzschild solution corresponds to the chargeless case, $Q=0$,
\begin{align}
ds^2=& F dt^2 - F^{-1} dr^2 -r^2 (d\theta^2 +\sin^2(\theta) d\phi^2) , \\
F =& 1-\frac{2M}{r} , \\
t\in&\Reals,& r>r_+=&2M, & (\theta, \phi)\in& S^2
\end{align}
The Schwarzschild solution, $Q=0$, has only one event horizon at $r=2M$, because the second root, $r_-$ coincides with the central singularity $r=0$. In the exterior region, the Schwarzschild solution is representative of all the sub critical solutions. The critical solution is similar; however, $F$ vanishes quadratically instead of linearly towards the horizon, and this affects the rate of decay of other quantities. 

In studying the exterior region of the \RN or Schwarzschild solutions it is common to introduce a new radial co-ordinate, the Regge-Wheeler tortoise co-ordinate, $r_*$, defined by
\begin{align}
\frac{dr}{dr_*} = F .
\end{align} 
This co-ordinate extends from $-\infty$ to $\infty$ and has the effect of ``pushing the horizon to negative infinity''. The original radial co-ordinate, $r$, is now treated as a function of $r_*$. The exterior region of the \RN solution is now represented by
\begin{align}
ds^2=& F dt^2 - F dr_*^2 -r^2 (d\theta^2 +\sin^2(\theta) d\phi^2) , \\
F =& 1-\frac{2M}{r}+\frac{Q^2}{r^2} , \\
t\in&\Reals,& r_*\in&\Reals, & (\theta, \phi)\in& S^2
\end{align}
In the Schwarzschild case, $r_*$ can be expressed simply in terms of $r$ and $M$ and has simple asymptotic behaviour. 
\begin{align}
\label{ersExplicitSchwarzschild}
r_* =& r + 2M \log\left(\frac{r-2M}{M}\right) + C_* ,\\
\lim_{r_*\rightarrow\infty} r_*^{-1} r =& 1 ,\\
\lim_{r_*\rightarrow-\infty} r_*^{-1} \log\left(\frac{r-2M}{M}\right) =& \frac{1}{2M} .
\end{align}
In the literature, $C_*$ is commonly taken, for simplicity, to be $-2M\log(2)$ \cite{Wald, MTW} or $-2M\log(M/2)$ \cite{EllisHawking}. In equation \eqref{eC*foralpa*},we will use a particular choice of $C_*$, coming from the geometry of the \RN solution. 

In the general sub critical case $,|Q|<M$, the expression for $r_*$ is slightly more complicated, but the asymptotic are the same. 
\begin{align}
\label{ersExplicitSubcrit}
r_*=& r +\frac{r_+^2}{r_+-r_-}\ln(\frac{r-r_+}{M})-\frac{r_-^2}{r_+-r_-}\ln(\frac{r-r_-}{M})+C_* , \\
\lim_{r_*\rightarrow\infty} r_*^{-1} r =& 1 , \label{3.2}\\
\lim_{r_*\rightarrow-\infty} r_*^{-1} \ln(\frac{r-r_{+}}{M})= &\frac{r_{+}-r_{-}}{r_{+}^2}\label{3.1} .
\end{align}

In the critical case, $|Q|=M$, the expression for $r_*$ and the asymptotics are inverse linear, instead of logarithmic, towards the event horizon. 
\begin{align}
\label{ersExplicitCrit}
r_*=&r+2M\ln(\frac{r-M}{M})-\frac{M^2}{r-M}+C_* , \\
\lim_{r_*\rightarrow\infty} r_*^{-1} r =& 1 , \label{3.2dup}\\
\lim_{r_*\rightarrow-\infty} r_* (r-M) =& -M^2 \label{eCritRhoNegInf} .
\end{align}

%
%
\subsection{The wave equation}
\label{ssWaveEqn}

We wish to study the linear wave equation, 
\begin{align}
\label{tildeLW}
\dAlembertianRN \tilde{u} =& 0, & \tilde{u}(1) =& \tilde{u}_0, & \frac{\partial}{\partial t}{\tilde{u}}(1)=&\tilde{u}_1 ,
\end{align}
in the exterior region of the \RN solution. In terms of the tortoise co-ordinate the d'Alembertian is
\begin{align}
\dAlembertianRN =& F^{-1}(\frac{\partial^2}{\partial t^2}-r^{-2}\frac{\partial}{\partial r_*}r^2\frac{\partial}{\partial r_*}) -r^{-2}\SLap ,
\end{align}
where $\SLap$ is the Laplace-Beltrami operator on the sphere. The substitution
\begin{align}
u =& r\tilde{u}
\end{align}
simplifies the wave equation to
\begin{align}
\label{LW}
\frac{\partial^2}{\partial t^2} u + \linH u =& 0, & u(1)=&u_0, & \frac{\partial}{\partial t} u(1) = u_1 , 
\end{align}
where the operator $H$ is composed of the following terms
\begin{align}
\linH =& \sum_{i=1}^3 H_i , \\
H_1=& -\frac{\partial^2}{\partial r_*^2} ,\\
H_2=& V(r) = \frac{1}{r} F \frac{dF}{dr} , \label{eDefnV}\\
H_3=& V_L(r)(-\SLap) ,\\
V_L(r)=& \frac{1}{r^2} F .
\end{align}
We refer to $V$ as the potential and $V_L$ as the angular potential. In the exterior region, $r>r_+$, the angular potential has a single critical point at (in terms of the $r$ variable) \footnote{For the critical \RN solution, the angular potential has a second critical point at the event horizon $r=M$. }
\begin{align}
\alpha = \frac{3M+\sqrt{9M^2-8Q^2}}{2} 
\end{align}
The critical point, $\alpha$, is extremely important. Geodesics at $r=\alpha$ and tangent to this surface will remain on the surface forever, and geodesics which approach this surface almost tangentially can approach the surface as $t\rightarrow\infty$ or orbit the black hole arbitrarily many times before escaping to $r_*\rightarrow\pm\infty$ \cite{EllisPhotonSphere,EllisLensing}. This geodesic surface at $r=\alpha$ is sometimes called the photon sphere. The region near $r=\alpha$ will also be the most difficult in which to prove decay of solution to the wave equation. This region only presents a problem for non- radial data. 

The value of $r_*$ corresponding to $r=\alpha$ will be denoted $\alpha_*$, and we will introduce the new radial co-ordinate
\begin{align}
\rho_* = r_*-\alpha_*
\end{align}
This corresponds to taking $r_*=\rho_*$ with the integration constant $C_*$ in equation \eqref{ersExplicitSchwarzschild}, \eqref{ersExplicitSubcrit}, or \eqref{ersExplicitCrit} chosen so that 
\begin{align}
\label{eC*foralpa*}
\alpha_*=0 . 
\end{align}
For this reason, $\frac{\partial}{\partial r_*} = \frac{\partial}{\partial \rho_*}$, etc. 

For simplicity, we typically use the measure
\begin{align}
\dmu = dr_* d^2\mu_{S^2} .
\end{align}
The space $\starman$ refers to $\Reals\times S^2$ with the measure $\dmu$. This defines $L^2(\starman)$ and, more generally, $L^p(\starman)$ for any $p\geq1$. Unless otherwise specified all norms and inner products are with respect to $L^2(\starman)$. 

 The function space $H^1(\starman)$ is the collection of functions for which $u'$ and $V_L^\frac12 \nabla_S^2 u$ are in $L^2(\starman)$. This space is not particularly useful, and we typically consider functions with finite energy or conformal charge, as defined in section \ref{sMethodsPD}. 

When dealing with the original function, $\tilde{u}$, we also use the measure
\begin{align}
\tildedmu = r^2 dr_* d^2\mu_{S^2} ,
\end{align}
and use $L^2(\tildestarman)$ to refer to $\Reals\times S^2$ with this measure. 

The two $L^2$ spaces coincide, in the sense that
\begin{align}
\|u\|_{L^2(\starman)} = \|\tilde{u} \|_{L^2(\tildestarman)}
\end{align}

The Schwartz space, $\solset$, refers to functions which are infinitely differentiable and for which, for a given $v\in\solset$, 
\begin{align}
\forall i,j,k \in \Integers^+, \exists C_{i,j,k} | \rho_*^i \dr^j \nabla_{S^2}^k v| < C_{i,j,k}
\end{align}

The Fourier transform on $\Reals\times S^2$ is defined by first making a spherical harmonic decomposition in the angular variable, and then applying the one dimensional Fourier transform on each spherical harmonic.

%
%
\section{Methods for pointwise decay}
\label{sMethodsPD}

In Minkowski space, $\Reals^{3+1}$, there is a conformal charge which is conserved and which dominates $t^2\|r^{-1}\nabla_{S^2} u\|^2$. Using a Sobolev estimate, the decay of the angular component of the $H^1(\Reals^3)$ norm implies decay of the the $L^6(\Reals^3)$ norm \cite{GinibreVelo}. We introduce an analogous conformal charge for the \RN solution. The growth of this charge is dominated by the sum of localised space-time integrals of $u$ and its angular derivative. Using a Sobolev estimate and the conformal charge, we reduce the proof of a pointwise in time $L^p(\starman)$ estimate to bounding localised space-time integrals of a solution and its angular derivative.

%
%

\subsection{Densities, energy conservation, and the conformal charge}
\label{ssDensitiesToConformal}

There are many formalisms for studying wave equations. In this subsection, we will introduce various densities and show that they satisfy differential relations. Integrating these relations will show that the energy is conserved and give an identity for the time derivative of the conformal charge. 

\begin{definition}
Given a pair of functions, $(v,w)\in\solset\times\solset$, the energy, radial momentum, and angular momentum densities are defined respectively by
\begin{align}
\EDens[v,w]=& \frac12 (w^2+ v'^2+ Vv^2+V_L \Lv\cdot\Lv) ,\\
\PsDens[v,w]=&wv' ,\\
\PlDens[v,w]=&w\Lv .
\end{align}
The energy is defined to be
\begin{align}
E[v,w]=\int \EDens[v,w] \dmu
\end{align}
\end{definition}

These densities satisfy the following differential relations. 
\begin{lemma}
\label{lDensityRelations}
If \uasoluinS, then the following relations hold 
\begin{align}
0 =& \dt \EDens[u(t),\dot{u}(t)]-\drs \PsDens[u(t),\dot{u}(t)] - \nabla_{S^2}\cdot (V_L \PlDens[u(t),\dot{u}(t)]) , \\
0 =& \dt \PsDens[u(t),\dot{u}(t)] -\drs\frac12(\dot{u}(t)\dot{u}(t)+u'(t)u'(t)-u(t)Vu(t)-\Lu(t)\cdot V_L\Lu(t))-\nabla_{S^2}\cdot(u'(t)V_L\Lu(t)) \\
 & - \frac12 u(t)V'u(t) -\frac12 \Lu\cdot V_L'\Lu , 
\end{align}
and the energy is conserved
\begin{align}
\Dt E[u(t),\dot{u}(t)]=& 0 .
\end{align}
\end{lemma}
\begin{proof}
The relations are proven using the method of multipliers, in which both sides of  the wave equation are multiplied by a quantity, typically a differential operator acting on $u$, and then the right hand side is rearranged. 

The multiplier $\dt u$ gives the relation for the derivative of the energy momentum.  
\begin{align}
0=&(\dt u)(\ddot{u}-u''+Vu-V_L(\Delta_{S^2})u)\\
=& \dot{u}\ddot{u}+\dot{u}'u' + \drs(\dot{u}u') +\dot{u}V u +\Lu\cdot(V_L \Lu)+\nabla_{S^2}\cdot(\dot{u}V_L\Lu)\\
=&\frac12\dt(\dot{u}\dot{u}+u'u'+uVu+\Lu\cdot V_L\Lu) -\drs(\dot{u}u') - \nabla_{S^2}\cdot(\dot{u}V_L\Lu)\\
=& \dt \EDens[u,\dot{u}] - \drs \PsDens[u,\dot{u}] - \nabla_{S^2}\cdot (V_L \PlDens[u,\dot{u}])
\end{align}

Since the integral of a pure spatial derivative is identically zero, integrating this result gives that $\Dt E[u(t),\dot(t)]$ is zero, and that the energy is conserved. 

The relation for the time derivative of the radial momentum comes from multiplier $\drs u$. 
\begin{align}
0=&
(\drs u)(\ddot{u}-u''+Vu+V_L(-\Delta_{S^2})u)\\
=&u'\ddot{u}-u'u''+u'Vu +u' V_L(-\Delta_{S^2})u\\
=& \dt(u'\dot{u})-\dot{u}'\dot{u}-u'u''+u'Vu+\Lu'\cdot V_L\Lu -\nabla_{S^2}\cdot(u'V_L \Lu)\\
=& \dt p_{\rho_*} -\drs\frac12(\dot{u}\dot{u}+u'u'-uVu-\Lu\cdot V_L\Lu)\\
&- \frac12 uV'u -\frac12 \Lu\cdot V_L'\Lu \\
&-\nabla_{S^2}\cdot(u'V_L\Lu)
\end{align}
\end{proof}

We now define the conformal charge in terms of the energy and momentum densities. In Minkowski space, the conformal multiplier is found by conjugating the time derivative with the discrete inversion symmetry of Minkowski space. The conformal multiplier is then used to define the conformal charge and its density by the same process which defines the energy and energy density from the time derivative. The \RN solution does not have a discrete inversion symmetry, so we define our conformal multiplier by formally taking the Minkowski conformal multiplier and replacing the Minkowski radial variable by $\rho_*$, the \RN radial variable. 

\begin{definition}
The conformal multiplier, the conformal charge density for a pair of functions $(v,w)\in\solset\times\solset$, and the conformal charge for the same pair are defined respectively by
\begin{align}
\confmult =& (t^2+\rho_*^2)\dt + 2t\rho_*\drs , \\
\confdens[v,w] =& (t^2+\rho_*^2) e[v,w] + 2t\rho_* \PsDens[v,w] , \\
\confchrg[v,w] =& \int \confdens[v,w] \dmu .
\end{align}
\end{definition}

The conformal charge density can be rewritten as a manifestly positive quantity. This form is more useful for making estimates. 

\begin{lemma}
\label{ConfPartII}
For any pair $(v,w)$, 
\begin{align}
\confdens[v,w] =&
 \frac14(t-\rho_*)^2(w-v')^2
+\frac14(t+\rho_*)^2(w+v')^2\\
&+\frac12(t^2+\rho_*^2)Vv^2
+\frac12(t^2+\rho_*^2)V_L(\Lv\cdot\Lv)
\end{align}
\end{lemma}
\begin{proof}
Only the time and radial derivative terms need to be rearranged. 
\begin{align}
(t-\rho_*)^2(w-v')^2
=&t^2w^2-2t^2wv'+t^2v'^2\\
&-2t\rho_*w^2+4t\rho_*wv'-2t\rho_*v'^2\\
&+\rho_*^2w^2-2\rho_*^2wv'+\rho^2v'^2
\end{align}
From this, 
\begin{align}
(t-\rho_*)^2(w-v')^2+(t+\rho_*)^2(w+v')^2
=&2(t^2+\rho_*^2)(w^2+v'^2)+8t\rho_*wv'
\end{align}
\end{proof}

In Minkowski space, since the time derivative is the generator of the time translation symmetry, the conformal multiplier is a composition of symmetries, and hence a symmetry itself. From Noether's theorem, the conformal charge it generates is conserved. Our conformal multiplier is not constructed from symmetries and does not generate a conserved quantity. Heuristically, the change of the \RN conformal charge should only involve the potentials, since, formally, the conformal multiplier is the same as in $\Reals^{3+1}$, and the wave equation differs only by the presence of potentials. The potentials appear in expressions of the form $2V+\rho_* V'$. We call them the trapping terms. In $\Reals^3$, the analogue of $V_L$ is $r^{-2}$, and the corresponding trapping term, $2V_L+rV_L'$, vanishes. 

\begin{lemma}
\label{lConfPartI}
If \uasoluinS, then
\begin{align}
\Dt \confchrg[u(t),\dot{u}(t)]
=&\int t (2V+\rho_*V')u^2 \dmu 
 +\int t(2V_L+\rho_*V_L') (\Lu\cdot\Lu) \dmu
\end{align}
\end{lemma}
\begin{proof}
We multiply the wave equation by the conformal multiplier, $\confmult u$, and then apply the relations from lemma \ref{lDensityRelations}. 
\begin{align}
0
=& (\confmult u) (\ddot{u}-u''+Vu+V_L(-\Delta_{S^2})u )\\
=& (t^2+\rho_*^2)\dot{u}(\ddot{u}-u''+Vu+V_L(-\Delta_{S^2})u )\\
& + 2t\rho_*u'(\ddot{u}-u''+Vu+V_L(-\Delta_{S^2})u )\\
=& (t^2+\rho_*^2)\left( \dt \EDens -\drs \PsDens - \nabla_{S^2}\cdot (V_L p_\omega) \right) \\
&+2t\rho_*\left( \dt \PsDens -\drs\frac12(\dot{u}\dot{u}+u'u'-uVu-\Lu\cdot V_L\Lu)-\nabla_{S^2}\cdot(u'V_L\Lu) \right)\\
&-2t\rho_*\left( \frac12 uV'u +\frac12 \Lu\cdot V_L'\Lu \right)
\end{align}

Integrating these terms and then integrating by parts in the angular derivatives eliminates the angular gradients. 
\begin{align}
0
=&\int (t^2+\rho_*^2)\left( \dt \EDens-\drs \PsDens - \nabla_{S^2}\cdot (V_L p_\omega) \right) \dmu\\
&+\int 2t\rho_*\left( \dt \PsDens -\drs\frac12(\dot{u}\dot{u}+u'u'-uVu-\Lu\cdot V_L\Lu)-\nabla_{S^2}\cdot(u'V_L\Lu) \right) \dmu\\
&-\int 2t\rho_*\left( \frac12 uV'u +\frac12 \Lu\cdot V_L'\Lu \right) \dmu \\
=&\int (t^2+\rho_*^2)\left( \dt \EDens-\drs \PsDens \right) \dmu\\
&+\int 2t\rho_*\left( \dt \PsDens -\drs\frac12(\dot{u}\dot{u}+u'u'-uVu-\Lu\cdot V_L\Lu) \right) \dmu \\
&-\int 2t\rho_*\left( \frac12 uV'u +\frac12 \Lu\cdot V_L'\Lu \right) \dmu
\end{align}

This is further simplified by integrating by parts in the radial variable and isolating pure time derivatives. 
\begin{align}
0
=& \int \dt \left( (t^2+\rho_*^2) +2t\rho_* \PsDens \right) \dmu \\
&+ \int -2t \EDens -2\rho_* \PsDens \dmu\\
&+ \int +2\rho_*\PsDens +2t\frac12(\dot{u}\dot{u}+u'u'-uVu-\Lu\cdot V_L\Lu) \dmu \\
&-\int 2t\rho_*\left( \frac12 uV'u +\frac12 \Lu\cdot V_L'\Lu \right) \dmu \\
=& \Dt \int (t^2+\rho_*^2)\EDens +2t\rho_* \PsDens \dmu \\
&-\int 2t (Vu^2 + \Lu\cdot V_L\Lu) \dmu\\
&-\int 2t\rho_*\left( \frac12 uV'u +\frac12 \Lu\cdot V_L'\Lu \right) \dmu 
\end{align}
\end{proof}

%
%
\subsection{Sobolev estimates}
\label{ssSobolev}

Our goal in this section is to bound the $L^6$ norm by the energy, the conformal charge, and a negative power of $t$. The main step in this is to prove a Sobolev estimate, which, roughly speaking, controls the $L^6$ norm by the $H^1$ norm. The analogous result in Minkowski space,$\Reals^{3+1}$, is that the spatial $L^6$ norm can be controlled by a third of a factor of the radial component of the $H^1$ norm and two thirds of the angular component of the inhomogeneous $H^1$ norm. The energy controls the radial component of the $H^1$ norm, and the conformal charge controls the product of a positive power of $t$ and the angular component of the $H^1$ norm. It is not obvious from the definition that the conformal charge controls the weighted $L^2$ norm needed in the inhomogeneous part of the angular $H^1$ norm, but it does \cite{GinibreVelo}. A similar result holds for the \RN solution, but, because of the weight appearing in the $L^6$ norm, we require two estimates on weighted $L^2$ norms instead of one. 

The two estimates in the following lemma are the estimates needed for the Sobolev estimate. Note that the terms appearing on the right are independent of $w$, the second argument of $\confchrg$. 

\begin{lemma}
\label{lWeightedL2ByConformal}
There is a constant such that, for all $(v,w)\in \solset\times\solset$,  
\begin{align}
\confchrg[v,w]\geq& C \langle v,\frac{t^2+\rho_*^2}{\rho_*^2+1}v\rangle , \\
\confchrg[v,w]t^{-2} \geq& C \langle v, \frac{1}{\rho_*^2+1}v\rangle .
\end{align}
\end{lemma}
\begin{proof}
This is a smoothed version of the argument used in $\Reals^n$ \cite{GinibreVelo}.  

The ingoing and outgoing wave terms can be isolated and rearranged. We introduce the notation $\confdenstr$ to denote twice these terms. 
\begin{align}
\confdenstr[v,w]
=&\|\rho_* w\|^2+\|t w\|^2+\|\rho_* v\|^2 +\|t v\|^2+\langle w,4\rho_*t v\rangle\\
=&\|tw+\rho_*v'\|^2+\|\rho_*w+tv'\|^2
\end{align}

The weighted term $hu =h(\rho_*)u$ is introduced into these terms. 
\begin{align}
\confdenstr[v,w]
=& \|tw+\rho_* v' -\rho_* h v\|^2+\|\rho_*w + t v' + t h v\|^2\\
&+2\langle \rho_* v',\rho_* h v\rangle - 2\langle t v',t h v\rangle
-\langle \rho_* h v,\rho_* h v\rangle-\langle t h v, t h v\rangle
\end{align}
Dropping the first two terms, which are strictly positive, and integrating by parts in the following pair yields:
\begin{align}
\confdenstr[v,w]
\geq& - \langle v, (2\rho_* h + (\rho_*^2-t^2) h') v\rangle - \langle v, (\rho_*^2 + t^2) h^2 v\rangle
\end{align}

A particular choice of $h$ can now be made in terms of parameters $a>0$ and $\epsilon>0$. 
\begin{align}
h(\rho_*)=& \epsilon \frac{\rho_*}{\rho_*^2+a}\\
h'(\rho_*)=&\epsilon \frac{a-\rho_*^2}{(\rho_*^2+a)^2}
\end{align}
We now substitute this choice of $h$ into the previous calculations. 
\begin{align}
-(2\rho_*h+(\rho_*^2-t^2)h') - (\rho_*^2 + t^2) h^2
=& -\epsilon(\frac{2\rho_*^2}{\rho_*^2+a} + (\rho_*^2-t^2)\frac{a-\rho_*^2}{(\rho_*^2+a)^2})
- \epsilon^2 \frac{(\rho_*^2 + t^2)\rho_*^2}{(\rho_*^2+a)^2}\\
=&-\epsilon \frac{\rho_*^4 + 3a\rho_*^2 -t^2a +t^2\rho_*^2}{(\rho_*^2+a)^2}
- \epsilon^2 \frac{\rho_*^4 + t^2\rho_*^2}{(\rho_*^2+a)^2}\\
=&-(\epsilon+\epsilon^2) \frac{\rho_*^4 }{(\rho_*^2+a)^2}
-\epsilon \frac{3a\rho_*^2}{(\rho_*^2+a)^2}
-\epsilon \frac{ -t^2a }{(\rho_*^2+a)^2}
-(\epsilon+\epsilon^2) \frac{t^2\rho_*^2}{(\rho_*^2+a)^2}
\end{align}
For $\epsilon\in (-1,0)$ the first, second, and forth terms are positive. For $\epsilon\in(-\frac12,0)$ and $\rho_*^2>\frac{a}{4}$, the forth term dominates the third by a factor of $2$. $\epsilon$ can be chosen sufficiently close to zero so that for $\rho_*^2<\frac{a}{4}$,
\begin{align}
0\leq& -\epsilon \frac{ -t^2a }{(\rho_*^2+a)^2} + V(\rho_*^2+t^2) 
\leq -\epsilon \frac{ -t^2a }{(\rho_*^2+a)^2} + H_2(\rho_*^2+t^2)
\end{align}
Thus, with $a=1$, 
\begin{align}
\confchrg[v,w]\geq\confdenstr[v,w]+\langle v, H_2 v\rangle
\geq& (-\epsilon-\epsilon^2) \langle v, \frac{t^2+\rho_*^2}{\rho_*^2+1} v\rangle
\end{align}
This proves the first statement in the theorem. The second part follows trivially by dropping the $\rho_*^2/(\rho_*^2+1)$ term from the first estimate and dividing by $t^2$. 
\end{proof}

We now turn to the Sobolev estimate. Because our main interest is making estimates in terms of the energy and conformal estimate, the estimate is expressed both in terms of $L^2$ norms of the derivatives of the function under consideration and in terms of the energy and conformal charge. Once again, the energy and conformal charge take a second argument $w$ which does not appear in the quantity estimated. 

\begin{lemma}
There is a constant $C$, such that if $(v,w)\in\solset\times\solset$, then
\begin{align}
\| F^\frac12 r^\frac{-2}{3} v \|_{L^6(\starman)} \leq& C (E[v,w]+\confchrg[v,w]t^{-2})^\frac16 \confchrg[v,w]^\frac13 t^{-\frac23}.
\end{align}
\end{lemma}
\begin{proof}
We begin by proving a Sobolev estimate which controls the $L^6$ norm by weighted $H^1$ norms. 

Following the standard argument \cite{stein, evans}, the proof starts with a $W^{1,1}\hookrightarrow L^\frac{3}{2}$ estimate. Take $\psi(\rho_*,\omega)$ infinitely differentiable and of compact support for simplicity. To this function $\psi$ associate a function on the sphere
\begin{equation}
I_1(\omega) \defin \int_\Reals |\dr \psi(\rho_*,\omega)|dr_*
\end{equation}
Integration and the Sobolev estimates for $\Reals$ and $S^2$ can be applied. The Sobolev estimate in 1 dimension follows directly from the fundamental theorem of Calculus. 
\begin{align}
|\psi(\rho_*,\omega)| \leq& I_1(\omega)\\
|\psi(\rho_*,\omega)|^{\frac{3}{2}} \leq& I_1(\omega)^\frac{1}{2}|\psi(\rho_*,\omega)|\\
\int_{S^2}|\psi(\rho_*,\omega)|^{\frac{3}{2}}d\omega \leq& \int_{S^2} I_1(\omega)^\frac{1}{2}|\psi(\rho_*,\omega)| d\omega
\end{align}
The second term in the final product on the right side is estimated by the spherical Sobolev estimate. The Sobolev estimate on $S^2$ follows from using a partition of unity on $S^2$ into co-ordinate charts and then applying the Sobolev estimate on $\Reals^2$ \cite{evans}. The notation $\|\psi\|_{\binom{1}{\omega}}(\rho_*)$ is introduced to denote the $L^1(S^2)$ norm of $\psi(\rho_*,\omega)$ with $\rho_*$ fixed, and the norm of ($TS^2$ valued) vector is equal to the norm of the length $\| \nabla_{S^2}\psi \|_X \defin \| (|\nabla_{S^2}\psi|) \|_X$. \begin{align}
\int_{S^2}|\psi(\rho_*,\omega)|^{\frac{3}{2}}d\omega
\leq&(\int_{S^2}I_1(\omega)d\omega)^\frac{1}{2} (\int_{S^2}|\psi(\rho_*,\omega)|^2d\omega)^\frac{1}{2}\\
\leq& \|\dr \psi\|_1^\frac{1}{2} (\|\nabla_{S^2}\psi\|_{\binom{1}{\omega}}(\rho_*)+\|\psi\|_{\binom{1}{\omega}}(\rho_*))
\end{align}
A radial weight $f^\alpha(\rho_*)$ can be introduced before integrating with respect to $dr_*$. For this proof, the exponents $\alpha$ and $\beta$ are used. They have no relation to the value of $r-\alpha$ governing the location of the photon sphere. 
\begin{align}
\int_\Reals\int_{S^2} f^\alpha(\rho_*) |\psi(\rho_*,\omega)|^\frac32 d\omega dr_* \leq&
\|\dr \psi\|_1^\frac{1}{2} \int_\Reals f^\alpha(\|\nabla_{S^2}\psi\|_{\binom{1}{\omega}}(\rho_*)+\|\psi\|_{\binom{1}{\omega}}(\rho_*))dr_*\\
\|f^{\frac{2\alpha}{3}}\psi\|_{\frac{3}{2}} \leq& \|\dr \psi\|_1^\frac{1}{3} (\|f^\alpha \nabla_{S^2}\psi\|_1+\|f^\alpha\psi\|_1)^\frac{2}{3}
\end{align}
The substitution $\psi=f^\beta |v_1|^4$ transforms this to a $H_1\hookrightarrow L^6$ Sobolev estimate. 
\begin{align}
\|f^{\frac{2\alpha}{3}+\beta}|v_1|^4\|_\frac{3}{2}\leq& (4\int f^\beta |v_1|^3|\dr v_1| + \beta f^{\beta-1}f'|v_1|^4 d^3\mu)^\frac{1}{3}\\
& \times (4\int f^{\alpha+\beta} |v_1|^3|\nabla_{S^2}v_1|d^3\mu+\int f^{\alpha+\beta} |v_1|^4 d^3\mu)^\frac{2}{3}\label{choosingweightsforSobolev}\\
(\int f^{\alpha+\frac{3}{2}\beta} |v_1|^6d^3\mu)^\frac{2}{3} \leq& 4(\int f^{2\beta} |v_1|^6d^3\mu)^\frac{1}{6} (\int|\dr v_1|^2d^3\mu +\int (\beta \frac{f'}{f})^2 |v_1|^2d^3\mu)^\frac{1}{6}\\
& \times (\int f^{2\alpha+2\beta-1} |v_1|^6d^3\mu)^\frac{1}{3}(\int f|(-\Delta_{S^2})^{\frac12}v_1|^2d^3\mu+\int f|v_1|^2d^3\mu)^\frac{2}{6}
\end{align}
We now choose
\begin{align}
\alpha+\frac{3}{2}\beta=&2\beta = 2\alpha+2\beta-1  \implies&
\alpha=&\frac{1}{2} , &
\beta=&1 ,
\end{align}
so that $\|f^\frac{1}{3} v_1\|_6$ can be cancelled. 
\begin{align}
\|f^\frac {1}{3} v_1\|_6\leq& 4(\|\dr v_1\|_2 + \|\frac{f'}{f} v_1\|^2)^\frac{1}{3} (\|f^{\frac{1}{2}} (-\Delta_{S^2})^{\frac12}v_1\|_2+\|f^\frac{1}{2} v_1\|_2)^\frac{2}{3}
\end{align}
The weight $f=r^{-2}$ is now taken to give a result analogous to the Sobolev estimate in $\Reals^3$. 
\begin{align}
\|r^\frac{-2}{3} v_1\|_6\leq&
C(\|\dr v_1\|_2 + \|-2Fr^{-1} v_1\|^2)^\frac{1}{3} (\|r^{-1} (-\Delta_{S^2})^{\frac12}v_1\|_2+\|r^{-1} v_1\|_2)^\frac{2}{3}
\end{align}

To complete the proof of the Sobolev estimate for the \RN solution, the substitution $v_1=F^\frac{1}{2}v$ and the inequality $Fr^{-1}\leq(1+\rho_*^2)^{-\frac12}$ are used. 
\begin{align}
\|F^\frac{1}{2} r^\frac{-2}{3} v\|_6
\leq& C(\|F^\frac{1}{2} \dr u\|_2 + \|\frac12 F^\frac12(2Mr^{-2}-2Q^2) v\| + \|-2F^\frac32 r^{-1} u\|^2)^\frac{1}{3}\\
&\times  (\|F^\frac12 r^{-1} (-\Delta_{S^2})^{\frac12}v\|_2+\|F^\frac12 r^{-1} v\|_2)^\frac{2}{3}\\
\leq& C(\|F^\frac{1}{2} \dr v\|_2 + \|(1+\rho_*^2)^{-\frac12} v\|)^\frac{1}{3}\\
&\times (\|F^\frac12 r^{-1} (-\Delta_{S^2})^{\frac12}v\|_2+\|F^\frac12 r^{-1} v\|_2)^\frac{2}{3}
\label{eSobolevEstimate}
\end{align}

We introduce the dummy function $w$ to act as the second argument of the energy and conformal charge. From the definition of the energy, the conformal charge, and lemma \ref{lWeightedL2ByConformal}, 
\begin{align}
\|F^\frac{1}{2} \dr v\|_2^2 \leq& \| \dr v\|_2^2 \leq E[v,w] , \\
\|F^\frac12 r^{-1} (-\Delta_{S^2})^{\frac12}v\|_2^2 \leq& \confchrg[v,w] t^{-2} , \\
\|F^\frac12 r^{-1} v\|_2^2 \leq
\|(1+\rho_*^2)^{-\frac12} v\|^2 \leq& C \confchrg[v,w] t^{-2} .
\end{align}
Substituting these into equation \eqref{eSobolevEstimate} proves the result. 
\end{proof}

We remark that if we repeat the same argument with $f=Fr^{-2}$ and $v_1=v$, the Sobolev type result, analogous to equation \eqref{eSobolevEstimate}, would be
\begin{align}
\|F^\frac{1}{3} r^\frac{-2}{3} v\|_6
\leq& C(\| \dr v\|_2 + \| (-2r^{-1}+6Mr^{-2}-4Q^2r^{-3}) v\|)^\frac{1}{3}\\
&\times (\|F^\frac12 r^{-1} (-\Delta_{S^2})^{\frac12}v\|_2+\|F^\frac12 r^{-1} v\|_2)^\frac{2}{3} .
\end{align}
The new weighted $L^2$ norm is also controlled by lemma \ref{lWeightedL2ByConformal}, for $t>1$, although, there is no longer the additional factor of $t^{-2}$. 
\begin{align}
\| (-2r^{-1}+6Mr^{-2}-4Q^2r^{-3}) v\|^2 \leq
\| (\frac{t^2+\rho_*^2}{\rho_*^2+1} u)^\frac12 v\| \leq& \confchrg[v,w] . 
\end{align}
Combining these results, we have control of a norm which decays less rapidly towards the event horizon but at the cost of less time decay on the right hand side. 
\begin{align}
\int |v|^6 F^2 r^{-4} \dmu \leq & C (E[v,w]+\confchrg[v,w]) \confchrg[v,w]^2 t^{-2},
\end{align}

%
%
\subsection{Local support of the trapping terms}
\label{ssTrappingLocal}

We refer to terms of the form $2V+\rho_*V_L'$, for both the potential and the angular potential, as trapping terms. In this section, we show that the trapping terms are positive only in a finite interval of $\rho_*$ values in the subcritical case. The functions $W$ and $W_l$ will refer to the positive part of the trapping terms. Through the conformal identity and the Sobolev estimates, this reduces the problem of finding point wise, weighted $L^6$ estimates to proving local decay estimates of the form $\int\int \chi |u|^2 \dmu dt + \int\int \chi |\Lu|^2 \dmu dt \leq C$. 

One of the factors in the derivative of the potential will require careful attention, both here and in subsection \ref{ssMorawetzCommutatorCalculation}. We introduce it with the notation $P_Q(r)$. 

\begin{definition}
For $Q\in[0,M]$, $P_Q(r)$ is defined by
\begin{align}
P_Q(r)=3Mr^3 -4(Q^2+2M^2)r^2 + 15MQ^2r -6Q^4
\end{align}
\end{definition}

We now show that the trapping term for the potential is positive only in a bounded set of $\rho_*$ values by computing the limit at $\pm\infty$. 

\begin{lemma}
\label{lPotentialTrappingCompactlySupported}
The derivative of the potential is given by 
\begin{align}
V'=-2Fr^{-7}P_Q(r)
\end{align}

For $|\rho_*|$ sufficiently large, $2V+\rho_* V'$ is negative. 

There is a compactly supported, positive, bounded function, $W$, such that $W>2V+\rho_* V'$
\end{lemma}
\begin{proof}
The derivative is computed from the definition of $V$ in equation \ref{eDefnV}.
\begin{align}
\drs V =&\frac{dr}{d\rho_*}\frac{\partial}{\partial r}(r^{-1}F\frac{\partial F}{\partial r})\\
=&F\frac{\partial}{\partial r}((r^{-1}-2Mr^{-2}+Q^2r^{-3})(2Mr^{-2}-2Q^2r^{-3}))\\
=&F\frac{\partial}{\partial r}(2Mr^{-3}-2(Q^2+2M^2)r^{-4}+6MQ^2r^{-5}-2Q^4r^{-6})\\
=& -2Fr^{-7}(3Mr^3-4(Q^2+2M^2)r^2+15MQ^2r-6Q^4)
\end{align}

From this, 
\begin{align}
2V+\rho_* V'
=&2r^{-1}F\left( \frac{dF}{dr} -\frac{\rho_*}{r^6} P_Q(r)\right)\\
=&2r^{-1}F\left( \left(\frac{2M}{r^2}-\frac{2Q^2}{r^3}\right) - \frac{\rho_*}{r^6}P_Q(r) \right)\\
=&2r^{-1}F\left( \left(\frac{2M}{r^2}-\frac{2Q^2}{r^3}\right) - \frac{\rho_*}{r}\left( \frac{3M}{r^2} -\frac{4(Q^2-2M^2)}{r^3} +\frac{15MQ^2}{r^4} -\frac{6Q^4}{r^5}\right) \right)
\end{align}
To show that this is negative for sufficiently large values of $|\rho_*|$, we will multiply by a positive factor and then show that the resulting quantity has a negative limit as $\rho_*\rightarrow\pm\infty$. As $\rho_*\rightarrow-\infty$, the subcritical and critical cases must be dealt with separately. 

In both the subcritical and critical cases, for $\rho_*\rightarrow\infty$, $\frac{\rho_*}{r}\rightarrow 1$ and $F\rightarrow 1$, so
\begin{align}
r^3(2V+\rho_* V')\rightarrow& 2M - 3M =-M
\end{align}

For $\rho_*\rightarrow -\infty$, the original radial variable has limit $r\rightarrow r_+=M+\sqrt{M^2-Q^2}$, and the limiting value of $P_Q(r_+)$ is 
\begin{align}
P_Q(M+\sqrt{M^2-Q^2})=& -2(2M^2-Q^2)(M^2-Q^2)-4M\sqrt{M^2-Q^2}(M^2-Q^2)
\end{align}

In the subcritical case, the term with $-\rho_*$ dominates, so it is sufficient to look at $P_Q(r_+)$. 
\begin{align}
\frac{r^7}{2F (-\rho_*)} (2V+\rho_* V')
\rightarrow& P_Q(r)\\
=& -2(2M^2-Q^2)(M^2-Q^2)-4M\sqrt{M^2-Q^2}(M^2-Q^2)\\
<&0
\end{align}

In the critical case, since $P_M(M)=0$ , it is necessary to multiply by a different positive factor before taking the limit. For $\rho_*\rightarrow-\infty$, from the asymptotic behaviour of $\rho_*$ given in equation \eqref{eCritRhoNegInf}
\begin{align}
\frac{r^4}{2F}\frac{1}{r-M} (2V+\rho_* V')
=& 2M\frac{r-M}{r-M} +\frac{\rho_*}{r-M}\frac{1}{r^3}P_M(r) \\
=& 2M +\rho_*(r-M)\frac{3M}{r^3}(r-2M)\\
\rightarrow& 2M-3M\\
=& -M
\end{align}

From these limits, it follows that $2V+\rho_* V'$ is negative for sufficiently large values of $|\rho_*|$. Since $V$ and $\rho_* V'$ are continuous, it follows that $2V+\rho_* V'$ can be bounded above by a compactly supported, positive, bounded function. 
\end{proof}

In the subcritical case, a similar calculation of limits as $\rho_*\rightarrow\pm\infty$ shows that $2V_L+\rho_*V_L'$ is positive only on a bounded set of $\rho_*$ values. In the critical case, the angular trapping term is positive as $\rho_*\rightarrow-\infty$, and we compute the rate at which it vanishes as $r\rightarrow r_+=M$. 

\begin{lemma}
In the subcritical case, for $|\rho_*|$ sufficiently large, $2V_L+\rho_* V_L'$ is negative, and there is a compactly supported, positive, bounded function, $W_L$, such that $W_L>2V_L+\rho_* V_L'$. 

In the critical case, there is a positive, bounded function, $W_L$, such that $W_L>2V_L+\rho_* V_L'$, $W_L$ is identically zero for sufficiently large $\rho_*$, and $W_L$ decays like $(r-M)F$ as $\rho_*\rightarrow-\infty$. 
\end{lemma}
\begin{proof}
As in the previous lemma, $2V_L+\rho_* V_L'$ will be shown to be negative from the fact that when a positive factor is applied, the limit as $\rho_*\rightarrow\pm\infty$ is strictly negative or the quantity diverges to negative infinity. For this lemma, it is simplest to separate the subcritical and critical cases, when evaluating the limits. 

From the definition of $V_L$, 
\begin{align}
\drs V_L
=&\frac{-2F}{r^3}\left(1-\frac{3M}{r}+\frac{2Q^2}{r^2}\right)
\end{align}
so that
\begin{align}
2V-\rho_* V'
=& \frac{2F}{r^3} \left(r-\rho_*(1-\frac{3M}{r}+\frac{2Q^2}{r^2}) \right)
\end{align}

In the subcritical, as $\rho_*\rightarrow\infty$, by the explicit expansion of $\rho_*$ given in equation \ref{ersExplicitSubcrit}, 
\begin{align}
\frac{r^3}{2F}\left(2V-\rho_* V'\right)
=& r - (r+\frac{r_+^2}{r_+-r_-}\log(\frac{r-r_+}{M})-\frac{r_-^2}{r_+-r_-}\log(\frac{r-r_-}{M}) +C)(1-\frac{3M}{r}+\frac{2Q^2}{r^2})\\
=& -\frac{r_+^2}{r_+-r_-}\log(\frac{r-r_+}{M})+\frac{r_-^2}{r_+-r_-}\log(\frac{r-r_-}{M}) + \BigOOne \\
\rightarrow&-\infty
\end{align}

In the critical case, as $\rho_*\rightarrow\infty$, by the expansion in \ref{ersExplicitCrit}, 
\begin{align}
\frac{r^3}{2F}\left(2V-\rho_* V'\right)
=&r -(r+2M \log(\frac{r-M}{M})-\frac{M^2}{r-M}+C)(1-\frac{3M}{r}+\frac{2Q^2}{r^2})\\
=&-2M \log(\frac{r-M}{M}) +\BigOOne\\
\rightarrow&-\infty
\end{align}

For $\rho_*\rightarrow-\infty$, we use both the asymptotics of $\rho_*$ from subsection \ref{ssRNsolu} and the geometry of the angular potential from subsection \ref{ssWaveEqn}. The critical points of $V_L$ are the roots of $r^2-3Mr+2Q^2$, which we denote by 
\begin{align}
\alpha_\pm =& \frac{3M \pm \sqrt{9M^2-8Q^2}}{2} .
\end{align}
Only one of these, $\alpha_+ = \alpha$ is in the exterior region $r>r_+$. We now multiply the angular trapping term by a positive term. 
\begin{align}
\frac{r^5}{2F}\frac{1}{-\rho_*}(2V_L+\rho_* V_L')
=& \frac{r^3}{-\rho_*}+r^2-3Mr+2Q^2\\
=& \frac{r^3}{-\rho_*}+(r-\alpha)(r-\alpha_-)
\label{eComputingAngularTrappingTowardsNegInf}
\end{align}
In the subcritical case, the relevant terms are ordered $\alpha_- < r_+ < \alpha_+$, so that in the limit $\rho_*\rightarrow-\infty$, $r\rightarrow r_+$, and
\begin{align}
\frac{r^5}{2F}\frac{1}{-\rho_*}(2V_L+\rho_* V_L')
\rightarrow (r_+-\alpha_-)(r_+-\alpha_+) 
<0
\end{align}

In the critical case, the relevant terms are not distinct, $\alpha_-=r_+=M < 2M=\alpha_+$,  so the previous argument does not hold. In fact, the angular trapping term is positive as $r \rightarrow r_+=M$. Instead, equation \eqref{eComputingAngularTrappingTowardsNegInf} can be rearranged as 
\begin{align}
(r^5)(2V_L+\rho_* V_L')
=& (r^3-\rho_*(r-\alpha)(r-\alpha_-)) 
\end{align}
From the asymptotics of $\rho_*$ in equation \eqref{ersExplicitCrit},
\begin{align}
r^3 -\rho_*(r-\alpha_-)(r-\alpha_+) 
=& r^3-(r+2M \ln\left(\frac{r-M}{M}\right) -\frac{M^2}{r-M} +C_*)(r-M)(r-2M)\\
=& (r^3+M^2(r-2M)) + \ln\left(\frac{r-M}{M}\right) (r-M)2M(2M-r) \\
&-(r-M)(r+C_*)(r-2m)
\end{align}
On the right, the last term clearly vanishes like $r-M$, and the second term is negative as $\rho_*\rightarrow-\infty$. The first term is a polynomial in $r$ which vanishes at $r=M$, so it must vanish at least linearly. Thus the right hand side is bounded above by a term which vanishes linearly in $r-M$. The potential term vanishes at a rate which is $F$ times faster. 
\end{proof}

Together, the results from this section show that a weighted $L^6$ norm is controlled, pointwise in time, by weighted space- time integrals of a solution and its derivative. 

\begin{proposition}
\label{pReductionToCompactlySupported}
In the subcritical case, $|Q|<M$, there are bounded, compactly supported functions $W$ and $W_L$ such that if \uasoluinS, 
then 
\begin{align}
\int |u|^6(t,\rho_*,\omega) F^3 r^{-4} \dmu
\leq& C (E[u_0,u_1] + \confchrg[u(t),\dot{u}(t)]t^{-2}) \confchrg[u(t),\dot{u}(t)]^\frac23 t^\frac{-4}{3} \\
\confchrg[u(t),\dot{u}(t)] 
\leq &\confchrg[u_0,u_1]  + \int \int 2t(W|u|^2 + W_L|\Lu|^2) \dmu dt
\end{align}

In the critical case, the same result holds with $W_L$ zero for $r$ sufficiently large and vanishing linearly in $r-M$ as $r\rightarrow r_+=M$. 
\end{proposition}

\section{Relativistic considerations on the event horizon}
\label{sRelConsid}

To study the behaviour of waves on the event horizon, it is common to introduce the Eddington-Finkelstein co-ordinates \footnote{The co-ordinates $(s_-,s_+)$ are typically denoted $(u,v)$ and the term Eddington-Finkelstein co-ordinates properly refers to the mixed co-ordinate systems $(u,r,\theta,\phi)$ or $(v,r,\theta,\phi)$. }\cite{DafermosRodnianski, MTW, Wald}
\begin{align}
s_- =& t-\rho_* \in \Reals , \\
s_+ =& t+\rho_* \in \Reals , \\
S_- =& -e^{-\frac{s_-}{4M}} \in (-\infty,0) , \\
S_+ =& e^{\frac{s_+}{4M}} \in (0,\infty) .
\end{align}
In fact, this range of co-ordinates only describes the exterior region of the black hole, , and the subcritical \RN solutions extend smoothly to positive $S_-$ and negative $S_+$ \footnote{In the critical case, $(S_-,S_+)=(0,0)$ is a singular point. Otherwise, the critical solution can be extended to an open set containing $S_->0, S_+>0$ and $S_-<0, S_-<0$. }. The outer event horizon corresponds to both the lines $S_-=0$ and $S_+=0$. The surface $S_-=0$ is the future event horizon, $\Hplus$, and $S_+=0$ is the past event horizon, $\Hminus$. To approach $\Hplus$, $t$ must diverge to infinite. The bifurcation sphere is the sphere given by $(S_-,S_+)=(0,0)$, with the spherical co-ordinates free. 

In this section, we discuss how solutions must decay near the bifurcation sphere if they have finite conformal charge. 

\subsection{Decay on the bifurcation sphere}
\label{ssBifurcationSphere}

To get $u\in L^2(\starman)$, we made several transformations and change of variables in section \ref{sWaveEqn}. To determine the physical constraints imposed by finite energy and finite conformal charge, we must return to the original function $\tilde{u} =r^{-1}u$, and use normalised vectors and their duals,
\begin{align}
\partial_T =& F^\frac{-1}{2} \partial_t , \\
\partial_R =&F^\frac{-1}{2} \partial_{\rho_*}, \\
dT =& F^\frac{1}{2} dt , \\
dR =& F^\frac{1}{2} d\rho_* .
\end{align}
From this, the natural measure on $(\rho_*,\omega)\in\Reals\times S^2$ is 
\begin{align}
\dTrue = F^\frac12 r^2 d\rho_* d\mu_{S^2} .
\end{align}
Since there is no bifurcation sphere in the critical case, we will only consider the subcritical cases in this subsection. 

We begin by writing the $L^2$ norm in terms of $\tilde{u}$ and the natural measure. 
\begin{align}
\| u\|_{L^2(\starman)}^2
= \| \tilde{u}\|_{L^2(\tildestarman)}^2
=& \int_{\tildestarman} |\tilde{u}|^2 r^2 dr_* d\omega \\
=& \int_{\tildestarman} |\tilde{u}|^2 F^{-\frac12} \dTrue .
\end{align}
From the first line, at the bifurcation sphere, the function $\tilde{u}$ must vanish at the boundary, but any power of $F\sim (r-r_+)$ is sufficient, and even inverse logarithmic decay is sufficient. 

We now write the energy in terms of $\tilde{u}$. 
\begin{align}
E[u,\dot{u}]
=& \int |\dot{u}|^2 + |u'|^2 +Fr^{-2}|\nabla_{S^2}u|^2 +2MFr^{-3}|u|^2 \dmu\\
=&\int (|\dot{\tilde{u}}|^2+ |Fr^{-1}\tilde{u}+\tilde{u}'| + Fr^{-2}|\nabla_{S^2}\tilde{u}|^2 + V|\tilde{u}|^2) r^2\dmu\\
=&\int (|\dot{\tilde{u}}|^2+|\tilde{u}'|^2  + Fr^{-2}|\nabla_{S^2}\tilde{u}|^2 )r^2 \dmu\\
&+ \int 2Fr\tilde{u}\tilde{u}' + F^2|\tilde{u}|^2 +Vr^2|\tilde{u}|^2 \dmu .
\end{align}
The co-ordinates $t$ and $\rho_*$ are singular at the bifurcation sphere, so we must rewrite the derivatives in terms of $\tilde{u}_T =\partial_T\tilde{u}$, $\tilde{u}_R =\partial_R\tilde{u}$, and the four gradient $\FourGrad\tilde{u}$. 
\begin{align}
E[u,\dot{u}]
=&\int (F^{-1}|\dot{\tilde{u}}|^2+F^{-1}|\tilde{u}'|^2  + r^{-2}|\nabla_{S^2}\tilde{u}|^2 ) \dFromE\\
&+ \int -(F\frac{dF}{dr}r+F^2)|\tilde{u}|^2 + (F^2+rF\frac{dF}{dr})|\tilde{u}|^2 \dmu\\
=& \int (|\tilde{u}_T|^2 + |\tilde{u}_R|^2 + r^{-2}|\nabla_{S^2}\tilde{u}|^2) \dFromE\\
=& \int |\FourGrad\tilde{u}|^2 F^\frac12 \dTrue .
\end{align}

For the energy to be bounded, it is sufficient that, far from the bifurcation sphere, the four gradient of $\tilde{u}$ is integrable and that near the bifurcation sphere, 
\begin{align}
|\FourGrad\tilde{u}|^2 <& C (r-r_+)^{-1+\epsilon} ,\\
|\FourGrad\tilde{u}| <& C (r-r_+)^{-\frac12+\epsilon} .
\label{eEnergyBifurcationCondition}
\end{align} 

Under this condition, 
\begin{align}
|\partial_r \tilde{u}| 
=& F^{-1} |\partial_{\rho_*} \tilde{u}|\\
=& F^{-\frac12} |\partial_R \tilde{u}|\\
\leq& F^{-\frac12} |\FourGrad\tilde{u}| \\
<& C (r-r_+)^{-1+\epsilon} .
\end{align}
The additional condition that $\tilde{u}$ vanishes at the bifurcation sphere restricts the size of $\tilde{u}$ near the bifurcation sphere, 
\begin{align}
|\tilde{u}| <& C (r-r_+)^\epsilon .
\end{align}

From integrating the length of the radial vector $\partial_R$, the distance from the bifurcation sphere, $s$, is given by
\begin{align}
s \sim (r-r_+)^\frac12 .
\end{align}
Thus, the sufficient condition on the energy in equation \ref{eEnergyBifurcationCondition} permits functions which grow like an arbitrarily small power of the distance from the bifurcation sphere. 

The calculations and conditions for the conformal charge are similar. Once again, we begin be rewriting the conformal charge in terms of $\tilde{u}$. 
\begin{align}
\confchrg[u,\dot{u}]
=& \int (\rho_*^2+1)(|\dot{u}|^2 + |u'|^2 +Fr^{-2}|\nabla_{S^2}u|^2 +2MFr^{-3}|u|^2) \dmu\\
=& \int (\rho_*^2+1)(|\tilde{u}_T|^2 + |\tilde{u}_R|^2 + r^{-2}|\nabla_{S^2}\tilde{u}|^2) \dFromE\\
&+ \int (\rho_*^2+1)(2Fr\tilde{u}\tilde{u}' + F^2|\tilde{u}|^2 +Vr^2|\tilde{u}|^2) \dmu\\
=& \int (\rho_*^2+1)(|\tilde{u}_T|^2 + |\tilde{u}_R|^2 + r^{-2}|\nabla_{S^2}\tilde{u}|^2) \dFromE\\
&+ \int -2Fr\rho_* |\tilde{u}|^2 \dmu\\ 
=& \int (\rho_*^2+1)(|\tilde{u}_T|^2 + |\tilde{u}_R|^2 + r^{-2}|\nabla_{S^2}\tilde{u}|^2) \dFromE\\
&+ \int -\frac{2\rho_*}{r} |\tilde{u}|^2 \dFromE \\
=& \int | \FourGrad\tilde{u}|^2 (\rho_*^2+1)F^\frac12 \dTrue 
+ \int |\tilde{u}|^2 (-\frac{2\rho_*}{r})F^\frac12 \dTrue
\end{align}
The weights in the conformal charge only differ by factors of $\rho_*$ which is logarithmic in $(r-r_+)$ near the bifurcation sphere. Thus, the sufficient condition given for the energy in equation \ref{eEnergyBifurcationCondition} is sufficient to guarantee that the conformal charge is also bounded.

\section{The Heisenberg-type relation}
\label{sHeisenbergType}

For the Schr\"odinger equation from quantum mechanics, 
\begin{align}
-i \psi_t + H\psi =0 ,
\end{align}
there is the well known Heisenberg type relation for the time derivative of the expectation value of a self-adjoint operator, $A$, 
\begin{align}
\Dt \langle \psi, A \psi\rangle =& \langle \psi,i[H,A]\psi\rangle .
\end{align}
This formulation is central to the standard interpretation of quantum mechanics which associates operators to physically observable quantities, and the expectation value to the mean observed value. This formulation was also used in the original proof of scattering for the quantum $n$-body problem \cite{DerezinskiGerard,HunzikerSigal,SigalSoffer,Soffer}. 

We begin by defining the commutator in the form sense. 

\begin{definition}
If $H$ and $A$ are two self-adjoint operators, and $D(A)\cap D(H)$ is dense in $L^2(\starman)$, then the commutator $[H,A]$ is defined to be the form $\Omega$ given by
\begin{align}
\Omega(v,w) = \langle Hv, Aw\rangle - \langle Av, Hw\rangle .
\end{align}
\end{definition}

The Heisenberg-type relation follows from the wave equation and this definition. 

\begin{lemma}
If $A$ is a time-independent, self-adjoint operator, and \uasoluinS, then 
\begin{align}
\Dt( \langle u,A\dot{u}\rangle - \langle\dot{u},A u\rangle) 
=& \langle u, [\linH,A] u\rangle ,
\end{align}
where $\langle u, [H,A] u\rangle$ is understood to mean the quadratic form $[H,A]$ evaluated on the pair $(u,u)$. 
\end{lemma}
\begin{proof}
We begin by computing the left hand side. 
\begin{align}
\Dt( \langle u,A\dot{u}\rangle - \langle\dot{u},A u\rangle) 
=& \langle \dot{u},A\dot{u}\rangle + \langle u,A\ddot{u}\rangle -\langle \ddot{u},Au\rangle -\langle \dot{u},A\dot{u}\rangle \\
=& \langle u, -A\linH u\rangle + \langle \linH u,Au\rangle .
\end{align}
Since $A$ is self-adjoint, 
\begin{align}
\Dt( \langle u,A\dot{u}\rangle - \langle\dot{u},A u\rangle) 
=& -\langle Au, \linH u\rangle + \langle \linH u,Au\rangle .
\end{align}
The right hand side is exactly the commutator. 
\end{proof}

Our method will be to find bounded propagation observables. A simple example of a propagation observable, $A$, which majorates an operator $G$, is one for which
\begin{align}
[H,A] = G^* G 
\end{align}
(where $G^*$ represents the adjoint of $G$). If $A$ is a bounded operator on the energy space, then, from the Heisenberg-type relation, 
\begin{align}
\int \|Gu\|^2 dt \leq \int \Dt ( \langle u,A\dot{u}\rangle - \langle\dot{u},A u\rangle) dt 
\leq 4 \|\dot{u}\| \|Au\| \leq C E[u,\dot{u}] .
\end{align}

The standard definition of a propagation observable is broader than this, and our definition will be even broader. 

\begin{definition} 
Given a pair of operators, $A$ and $G$, with $D(A)\cap D(\linH)\cap D(G)$ dense in $L^2(\starman)$, the operator $A$ is a propagation observable which majorates $G$ if there is a pair of bounded, non-zero operators, $X_1$ and $X_2$, for which 
\begin{align}
X_1 + X_2 =& \id , \\
[\linH,A] =& G^* X_1 G + \text{lower order terms} ,
\end{align}
where ``lower order terms'' refers to the sum of operators $R$ for which either 
\begin{enumerate}
\item $R$ is a bounded operator and for all $u\in\solset$ which solve the wave equation, 
\begin{align}
\int \langle u, R u \rangle dt \leq C ,
\end{align}
\item or $R$ is an operator with domain $D(R) \subset D( (G^*G)^\frac12 )$, and for all $u\in D( (G^*G)^\frac12 )$, 
\begin{align}
\langle u, R u \rangle \leq C \langle u, (G^*G)^{1-\epsilon} u \rangle .
\end{align}
\end{enumerate}
\end{definition}

If $A$ maps $\solset$ to $\solset$, even if it is not self-adjoint, then, since $\linH$ maps $\solset$ to $\solset$, for $v\in\solset$, $\linH(Av)-A(\linH v)$ is well-defined on $\solset$. Since $\linH$ is self-adjoint, by a similar calculation, 
\begin{align}
\Dt( \langle u,A\dot{u}\rangle - \langle\dot{u},A u\rangle) 
=& \langle u, (\linH A- A\linH)u\rangle .
\end{align}
Since products of smooth functions, integer powers of the radial derivative, and integer powers of angular derivatives map $\solset$ to $\solset$, if $A$ is of this form, we can write
\begin{align}
\Dt( \langle u,A\dot{u}\rangle - \langle\dot{u},A u\rangle) 
=& \langle u, [\linH,A]u\rangle , \\
[\linH,A] = \linH A - A \linH .
\end{align}

For anti-self-adjoint operators, since multiplication by $i$ commutes with $\linH$, we can define
\begin{align}
\langle u, [\linH, A] u \rangle 
=& \langle u, [\linH, -i A] i u \rangle .
\end{align}
From this, for an anti-self-adjoint operator, we have the same Heisenberg-type relation.

\section{Morawetz Estimates}
\label{sMorawetz}

The goal of this section is to prove bounds on weighted space-time norms of solutions. We begin by introducing a propagation observable, $\gamma$, show it majorates a weighted quantities, and conclude with a Gronwall's type argument to integrate the Heisenberg-type relation. These estimates are proven using a spherical harmonic decomposition. We show these estimates have a uniform nature in the spherical harmonic parameter, so that we can recover an estimate for general $u$. In later sections, the contribution from $H_2$ is controlled by the local decay estimate, and we can use a uniform multiplier. 

In $\Reals^{3+1}$, the radial derivative can be used as a propagation observable to control the time integral of $|u(t,\vec{0})|^2$. This can be thought of as a weighted space- time integral with the $\delta$- function as a weight. Essentially, our propagation observable is a smooth version of the radial derivative, which leads to a smooth weight. This follows \cite{LabaSoffer, BlueSoffer}. 

The spherical Laplacian is well known to have discrete spectrum. On each spherical harmonic we will use $l$ to refer to the standard spherical harmonic parameter and $\tl$ to refer to the value of the square root of the Laplacian on that harmonic, hence $\tl^2=l(l+1)$. We will also introduce the effective potential on each spherical harmonic, given by
\begin{align}
\Vtl = H_2 + H_3 = V + \tl^2\Vtl .
\end{align}
If necessary, we will use $\Ptl$ to denote the projection on to the spherical harmonic with parameter $l$, but typically, in this section, we will assume that $u$ has support on only one spherical harmonic. In this case, $u$ satisfies
\begin{align} 
\ddot{u}-u''+\tl^2\Vtl u =0 .
\end{align}
We will show in lemma \ref{lVtlUniquePeak} that $\Vtl$ has a single critical point, which is a maximum. We denote this critical point $\al$. 

The Morawetz-type operator is defined on each spherical harmonic in terms of the radial co-ordinate $\rho_*$. The definition involves a weight, $g_{\sigma}$, which is defined for $\sigma>1$ so that $g_\sigma$ remains bounded. 

\begin{definition} 
\label{defngamma}
Given $\sigma>1$, $b>0$, and $l\in\Naturals$, the Morawetz-type multiplier $\gamma_{l,\sigma}$ is defined by
\begin{align}
g_{l,b,\sigma}(\rho_*)\defin&\int_{0}^{b(\rho_*-\al)} \frac{1}{(1+|\tau|)^{\sigma}}d\tau\\
\gamma_{l,b,\sigma} \defin&\frac{1}{2}(g_{l,b,\sigma}(\rho_*)\dr+\dr g_{l,b,\sigma}(\rho_*))\\
=&g_{l,b,\sigma}(\rho_*) +\frac{1}{2}g_{l,b,\sigma}'(\rho_*)
\end{align}
In all cases $\sigma$ will not vary so the notation $\gamma=\gamma_{l,\sigma}$ and $g=g_{l,\sigma}$ will be used. The projection operator onto the $l^{\text{th}}$ spherical harmonic is denoted $\Ptl$. When working on more than one spherical harmonic, we also use $\gamma$ to denote 
\begin{align}
\sum_{l} \gamma_{l,b,\sigma} \Ptl . 
\end{align}
\end{definition}

We note that as $l$ varies, the observable $\gamma_{l,b,\sigma}$ is simply translated. We will also show that $\al\rightarrow0$ from which it follows that $g_{l,b,\sigma}$ has a limit in $L^\infty$. 

\subsection{Preliminary bounds}
\label{ssMorawetzPrelim}

\begin{lemma}
\label{L18.1}
If $\sigma>1$ and $b>0$, then there are constants $C_\sigma$ and $C_b$, independent of $l$, such that for all $u\in\solset$, 
\begin{align}
\langle u,\gamma_{l,b,\sigma}u\rangle =& 0\\
\|\gamma u\|_{L^2} \leq& C_\sigma\|u'\|_{L^2} + C_b \|(1+\rho_*^2)^{-\frac\sigma2}u\|_{L^2}\\
\leq& C_\sigma\sqrt{E[u]} + C_b\|(1+\rho_*^2)^{-\frac\sigma2}u\|_{L^2}
\end{align}
\end{lemma}
\begin{proof}
We work on a fixed spherical harmonic. Since $u\in\solset$, derivatives can be moved about freely. 
\begin{align}
\langle u,\gamma u\rangle =&\langle u,\frac{1}{2}(g\dr+\dr g)u\rangle\\
=&\frac{1}{2}(\langle u,g\dr u\rangle -\langle g\dr u, u\rangle)\\
=&0
\end{align}
The last equality holds since the $L^2$ inner product is symmetric on real valued functions. 

Since $g=\int_0^{b(\rho_*-\al)} (1+|\tau|)^{-\sigma} d\tau$, $g'$ is bounded by $C_b(1+\rho_*^2)^{-\sigma/2}$. Since the integrand is continuous and decays as $\tau^{-\sigma}$ for $\sigma>1$, $g$ is bounded by some constant $C_\sigma$. Using this and the fact that $u\in\solset$, the following holds.
\begin{align}
\|\gamma u\|=&\|gu' +\frac{1}{2}g'u\|\\
\leq& \|gu'\|+\|\frac{1}{2}g'u\|\\
\leq& C_\sigma\|u'\|+C_b\|(1+\rho_*^2)^{-\sigma}u\|
\end{align}
The constant $C_\sigma$ is given by the integral of $(1+|\tau|)^{-\sigma}$ from $0$ to $\infty$. This is independent of $l$. The constant $C_b$ is based on the equivalence of $b(1+b|\rho_*-\al|)^{-\sigma}$ and $(1+\rho_*^2)^{-\sigma/2}$. Since $\al$ converges to zero by lemma \ref{lVtlUniquePeak}, this equivalence is uniform in $l$. 
\end{proof}

\subsection{Computation of Morawetz commutators}
\label{ssMorawetzCommutatorCalculation}

We show that the propagation observable, $\gamma$, majorates $(1-\rho_*^2)^{-\sigma-1}$. We will prove this on each spherical harmonic, and then show there is a uniform lower bound on these estimates, to get an estimate for uniform $l$. 

We first show that on each spherical harmonic, the effective potential has a unique maximum. 

\begin{lemma}
\label{lVtlUniquePeak}
For $l\geq0$, $M>0$ and $0\leq |Q|<M$, the effective potential $\Vtl$ has a unique critical point, which is a maximum. The maxima, $\al$, converge to $0$. 
\end{lemma}
\begin{proof}
\newcommand{\pot}{V}
\newcommand{\potL}{V_L}
\newcommand{\potl}{\Vtl}
The Reissner-Nordstr\o m potentials are
\begin{align}
\pot=& \frac{2M}{r^3} (1-\frac{2M}{r}+\frac{Q^2}{r^2}), \\
\potL=& \frac{1}{r^2} (1-\frac{2M}{r}+\frac{Q^2}{r^2}), \\
\potl=& \pot + \tl^2 \potL .
\end{align}

Our goal is to show that in the outer region, the effective potential $\potl$ has a single critical point. The outer region is given by $r>r_+$ with 
\begin{align}
r_\pm=& M \pm \sqrt{M^2-Q^2} ,
\end{align}
being the zeroes of the factor $1-2M/r+Q^2/r^2$.

We will do this by showing that a positive factor times the derivative of each of the potentials increases from a nonpositive value to infinite from $r=2M$ to $r\rightarrow\infty$ and that the potentials are negative for $r$ from the horizon at $r=r_+$ to $r=2M$. 

The derivatives are given by
\begin{align}
\pot'=& -2Fr^{-7} \PolyZero, & \PolyZero=&3Mr^3-4(Q^2+2M^2)r^2 +15MQ^2r -6Q^4 ,\\
\potL'=& -2Fr^{-7} \PolyL, & \PolyL=& r^4-3Mr^3 +2Q^2r^2 , \\
\potl'=& \pot'+\tl^2\potL'= -2Fr^{-7}\Polyl, & \Polyl =& \PolyZero+ \tl^2 \PolyL . 
\end{align}
We will also use 
\begin{align}
\I=& r^{-2} \PolyZero= 3Mr-4(Q^2+2M^2)+15MQ^2r^{-1}-6Q^4r^{-2} ,\\
\IL=&r^{-2} \PolyL=  r^2 -3Mr +2Q^2 , \\
\Il=&r^{-2} \Polyl=  \I + \tl^2 \IL .
\end{align}
Showing that $\Polyl$ has a single root in the exterior region is equivalent to showing that $\Il$ has a single root in the exterior region. 

We start by showing that $I$ and $\IL$ are increasing in the region $r\geq 2M$. The derivative of $\I$ is
\newcommand{\doriginal}{\partial_r}
\begin{align}
\doriginal \I =& 3M - 15MQ^2r^{-2} + 12Q^4 r^{-3}\\
=& r^{-3}(3Mr^3 -15 MQ^2r+12Q^4) .
\label{eItermderiv}
\end{align}
The derivative of the term in brackets, for $r\geq 2M$ and $Q\leq M$, can be estimated by
\begin{align}
\doriginal (3Mr^3 -15 MQ^2r+12Q^4)
= 9Mr^2-15MQ^2 > 0.
\end{align}
At $r=2M$, the quantity in brackets in \eqref{eItermderiv} is $24M^4-30M^2Q^2+12Q^4>0$, thus $\doriginal\I$ is positive, and $\I$ is increasing from $r=2M$. The derivative of $\IL$, $\doriginal\IL= 2r-3M >0$, is positive for $r\geq 2M$ (and even for $r>3M/2$). 

The value of these quantities at $r=2M$ is 
\begin{align}
\I=& 6M^2-4(Q^2+2M^2) + 7.5Q^2 -\frac{3Q^4}{2M^2} \\
=&\frac{1}{M^2}( -2M^4 + 3.5M^2Q^2 - 1.5 Q^4 ) 
\leq 0 , \\
\IL=& 4M^2 - 6M^2 + 2Q^2 \leq 0 .
\end{align}
Thus each of $\I$ and $\IL$ are increasing for $r\geq2M$, diverge to infinite as $r\rightarrow\infty$, and are nonpositive for $r=2M$. 

We now show that $\I$ and $\IL$ is negative for $r$ between $r_+$ and $2M$. To show $\I$ is negative, it is sufficient to show $\PolyZero= r^2\I$ is negative. This is a cubic, and we observe a number of properties about it before drawing a conclusion. The derivative is quadratic with two roots, 
\begin{align}
\doriginal(\PolyZero)
=& 9Mr^2 - 8(Q^2+2M^2)r +15MQ^2 ,\\
r
=&\frac{4Q^2+8M^2 \pm \sqrt{64M^4 - 71M^2Q^2 +16 Q^4}}{9M} .
\end{align}
The discriminant is bounded above by $64M^2-64M^2Q^2+16Q^4=(8M^2-4Q^2)^2$. Therefore, the upper root is less than $16M/9<2M$. The discriminant is bounded from below by $9Q^4$, so the lower root is bounded above by $(8M^2+Q^2)/9M<M$. The value of $\PolyZero=r^2\I$ at $r=r_+$ is 
\begin{align}
\PolyZero(r_+)
=& -2(2M^2-Q^2+2M\sqrt{M^2-Q^2})(M^2-Q^2) \leq 0 .
\end{align}
Considered as a function on the entire real line, $\PolyZero$ goes from negative infinite to a local maximum value at the lower root of the derivative, and then decreases to a non positive value at $r=r_+$. Since $r^2\I$ is cubic, it can only have one root for $r>r_+$. Since that root is is in $r\geq2M$, the cubic must be negative for $r\in[r_+,2M]$. 

The function $\IL$ is much easier to estimate. It has two roots at
\begin{align}
r=\frac{3M\pm\sqrt{9M^2-8Q^2}}{2} .
\end{align}
For $Q\leq M$, the lower root is less than or equal to $M$ and the upper root is greater than or equal to $2M$. Thus, $\IL$ is negative from $r_+=M+\sqrt{Q^2-M^2}$ to $2M$. 

We know that $\I$ and $\IL$ are negative for $r$ between $r_+=M+\sqrt{Q^2-M^2}$ and $2M$, are increasing for $r\geq2M$, and go to infinite as $r\rightarrow\infty$. Therefore, the weighted sum $\Il=\I+\tl^2\IL$ has the same properties and must have a single root for $r> r_+$. 

Since the derivative goes from negative to positive, these critical points are maxima. Since $\potl \tl^{-2}\rightarrow \potL$ and the maximum of $\potL$ occurs at $0$, if $\al$ denotes the maximum point of $\potl$, then $\al\rightarrow 0$. 
\end{proof}

We now show that the commutators uniformly bound a polynomially decaying term, for an appropriate choice of $b$. This value of $b$ will always be used in $\gamma_{l,\sigma}$ from now on. There are additional positive terms which are also dominated, but we are better able to take advantage of them in the later sections using a slightly modified multiplier. 

\begin{lemma}
\label{lBasicMorawetzCommutatorWasLikeL19andL28}
If $\sigma>1$, $M>0$, and $|Q|\leq M$, then there is a choice of $b$ for which there is a constant $C$ such that for all $u\in\solset$, 
\begin{align}
\langle u,[\linH,\gamma]\rangle \geq& C \langle u,(1+\rho_*^2)^{-\frac\sigma2-1}u\rangle .
\end{align}
\end{lemma}
\begin{proof}
We work on a single spherical harmonic and sum at the end of the argument. The commutator can be computed exactly as 
\newcommand{\drrr}{\frac{\partial^3}{\partial r^3}}
\begin{align}
[\linH,\gamma]
=& -\frac12 ( \drr g\dr + \drrr g-g\drrr -\dr g\drr ) - \tl^2 (g\dr\Vtl-\Vtl g\dr) \\
=& -2\dr g \dr -\frac12 g''' - \tl^2 \Vtl' .
\end{align}

We will use $x$ to denote $\rho_*-\al$ 
\begin{align}
g'=& \frac{b}{(1+b|x|)^\sigma} \\
g''=& \frac{-b^2\sigma \sgn(x)}{(1+b|x|)^{\sigma+1}} \\
g'''=& -b^2\sigma 2\delta_0(x) + \frac{b^3\sigma(\sigma+1)}{(1+b|x|)^{\sigma+2}}
\end{align}

To estimate the contribution from the $\delta$ function, we introduce a smooth, weakly decreasing, compactly support function, $\chi$ which is identically one in a neighbourhood of $x=0$ (this is unrelated to any other cut off function used elsewhere in this paper). Using integration by parts, and temporarily treating $u$ as a function of $x$, 
\begin{align}
0
=& \int \dr (x\chi u^2) dx \\
=&\int \chi u^2 dx + \int x\chi'u^2dx + \int 2x\chi u\dr u dx , \\
\int \chi u^2 dx
\leq& \int (|x\chi'|+x^2\chi) u^2 dx + \int \chi (\dr u)^2 dx ,\\
u(0)^2
=& -\int_0^\infty \dr (\chi u^2) dx \\
=& -\int_0^\infty \chi' u^2 dx - \int_0^\infty 2\chi u\dr u dx \\
\leq& \int_0^\infty |\chi'| u^2 dx + \int_0^\infty \chi u^2 dx + \int_0^\infty \chi (\dr u)^2 dx . 
\end{align}
Applying a symmetric argument for $(-\infty,0]$ and substituting the result for $\int \chi u^2 dx$, we have
\begin{align}
\label{eControlDelta}
2u(0)^2 
\leq \int (|\chi'| + |x\chi'|+x^2\chi)u^2 dx + \int 2\chi (\dr u)^2 dx .
\end{align}

We use a similar integration by parts argument to estimate $(1+b|x|)^{-\sigma-2} u$. 
\begin{align}
0=&\int \dr (g''u^2) dx \\
=&\int g''' u^2 dx + \int 2 g''u \dr u dx \\
b^3 \int \frac{\sigma(\sigma+1)}{(1+b|x|)^{\sigma+2}} dx
\leq& 2\sigma b^2u(0)^2 \\
&+ \left(\int b^3 \frac{\sigma(\sigma+1)}{(1+b|x|)^{\sigma+2}} u^2 dx\right)^\frac12 \left(\int\frac{2\sigma}{\sigma+1} b\frac{2(\dr u)^2}{(1+b|x|)^\sigma} dx \right)^\frac12 
\label{eResultOfIBP}
\end{align}
Substituting equation \eqref{eControlDelta}, 
\begin{align}
b^3 \int \frac{\sigma(\sigma+1)}{(1+b|x|)^{\sigma+2}} dx
\leq& \sigma b^2 \int (|\chi'| + |x\chi'|+x^2\chi)u^2 dx 
+ \sigma b^2 \int 2\chi (\dr u)^2 dx \\
&+ \frac12 \int b^3 \frac{\sigma(\sigma+1)}{(1+b|x|)^{\sigma+2}} u^2 dx 
+ \frac12 \int\frac{2\sigma}{\sigma+1} b\frac{2(\dr u)^2}{(1+b|x|)^\sigma} dx \\
b^3 \int \frac{\sigma(\sigma+1)}{(1+b|x|)^{\sigma+2}} dx
\leq& 2 \sigma b^2 \int (|\chi'| + |x\chi'|+x^2\chi)u^2 dx 
+ 2\sigma b\int (\frac{2}{\sigma+1} \frac{1}{(1+b|x|)^\sigma}+2b\chi) (\dr u)^2 dx
\end{align} 

The commutator term can now be estimated
\begin{align}
\int 2g'&(\dr u)^2 -\frac12 g''' u^2 - g\Vtl' u^2 dx \\
=& \int \frac{2b}{(1+b|x|)^\sigma} (\dr u)^2 dx + \sigma b^2 u(0)^2 
-\frac12 \int \frac{b^3 \sigma(\sigma+1)}{(1+b|x|)^{\sigma+2}} u^2 dx 
- g\Vtl' u^2 dx \\
\geq& \int \frac{2b}{(1+b|x|)^\sigma} (\dr u)^2 dx + \sigma b^2 u(0)^2 \\
&- \sigma b^2 \int (|\chi'| + |x\chi'|+x^2\chi)u^2 dx 
- \sigma b\int (\frac{2}{\sigma+1} \frac{1}{(1+b|x|)^\sigma}+2b\chi) (\dr u)^2 dx
- g\Vtl' u^2 dx \\
\geq& 2b\int (\frac{1}{\sigma+1}\frac{1}{(1+b|x|)^{\sigma+2}} -\sigma b\chi) (\dr u)^2 dx 
+ \int (-\sigma b^2 (|\chi'| + |x\chi'|+x^2\chi) - g\Vtl') u^2 dx 
+ \sigma b^2 u(0)^2 
\end{align}
The term $\sigma b u(0)^2$ is positive and will be ignored. The term $(\sigma+1)^{-1}(1+b|x|)^{-\sigma-2} -\sigma b \chi$ is independent of $l$, so we can choose $b$ sufficiently small so that this is bounded below by $c(1+b|x|)^{-\sigma-2}$. Since $\Vtl' = V' + V_L'\tl^2$, and the $\al$ converge, the value of $\Vtl''$ is uniformly bounded from below at $\al$, and $b$ can be chosen sufficiently small so that the quadratically vanishing quantity $-\sigma b^2 (|\chi'| + |x\chi'|+x^2\chi) - g\Vtl') $ is bounded below by $-c x^2 \chi $, uniformly in $l$. Thus there are uniform constants, $c_1$ and $c_2$, for which, 
\begin{align}
\int 2g'&(\dr u)^2 -\frac12 g''' u^2 - g\Vtl' u^2 dx \\
\geq& c_1 \int \frac{1}{(1+b|x|)^\sigma}(\dr u)^2 - x^2 \chi u^2 dx \\
\geq& c_2 \int \frac{1}{(1+b|x|)^{\sigma+2}} u^2 dx .
\end{align}
Since the $\al$ converge, there is a uniform equivalence between $(1+b|x|)^{-\sigma-2}$ and $(1+\rho_*^2)^{-\sigma/2-1}$. Integrating over the angular variables gives
\begin{align}
\langle u, [\linH,\gamma]u\rangle
\geq& C \langle u, \frac{1}{(1+\rho_*^2)^{\frac\sigma2+1}} u\rangle .
\end{align}
\end{proof}

\subsection{$L^2$ local decay estimate}

We begin with a type of Gronwall's estimate. 

\begin{lemma}
\label{L38.2}
If $\theta:[0,\infty)\rightarrow[0,\infty)$ with uniformly bounded derivative, $\epsilon\in(0,\frac{1}{3})$, and there are constants $C_1$, $C_2$, and $T$, such that for $t>T$
\begin{equation}
\label{TechLemmaIntegralCondition}
\int_0^t\theta(\tau)^2d\tau \leq C_1+C_2t^\epsilon\theta(t)^{1-\epsilon}
\end{equation}
then there is a sequence $\{t_i\}$ such that $\lim_{i\rightarrow\infty}t_i=\infty$ and
\begin{equation}
\lim_{i\rightarrow\infty}t_i^\epsilon \theta(t)^{1-\epsilon}=0
\end{equation}
\end{lemma}
\begin{proof}
Successively stronger bounds on $\theta(t)$ will be proven. 

The bound on the derivative implies $\theta(t)$ and $t^\epsilon \theta(t)^{1-\epsilon}$ are  linearly bounded above. 

Suppose there is no sequence $\{t_i\}$ on which $\theta(t_i)\rightarrow 0$, then $\theta(t)$ is bounded below by a constant. By the integral condition $t^\epsilon\theta(t)^{1-\epsilon}$ is bounded below by a linear function. From this $\theta(t)^2$ is bounded below by a quadratic, and its integral is bounded below by a cubic. This contradicts the linear upper bound of the integral. Therefore, there is a sequence $\{t_i\}$ such that $t_i\rightarrow \infty$ and $\theta(t_i)\rightarrow 0$. 

Since $\epsilon<\frac{1}{3}$, 
\hide{
\begin{eqnarray}
\frac{-\epsilon}{1-\epsilon}&>& -\frac{1}{2} \nonumber
\end{eqnarray}
and
\begin{eqnarray}
\frac{1-\epsilon}{1+\epsilon} &>& \frac{\epsilon}{1-\epsilon}\nonumber
\end{eqnarray}
Thus,} 
there exists $r<1$ such that
\begin{equation*}
-\frac{1-\epsilon}{2-r+r\epsilon} < -\frac{\epsilon}{1-\epsilon} 
\end{equation*}
Now choose $\delta$ negative such that $\delta \geq \frac{-\epsilon}{1-\epsilon} $, from which it follows that
\begin{eqnarray}
r\delta &<& 1+\frac{2\delta}{1-\epsilon}\label{rdeltabound}
\end{eqnarray}

Suppose $\exists K>0, S>0: \forall t>S: \theta(t)>Kt^{\delta}$, then by the integral condition (\ref{TechLemmaIntegralCondition}) there is a positive $C_3$ such that
\begin{equation}
C_3 t^{1+2\delta} \leq C_1 + C_2t^{\epsilon}\theta(t)^{1-\epsilon}
\end{equation}
Since $\delta\geq\frac{-\epsilon}{1-\epsilon}>\frac{-1}{2}$, $1+2\delta$ is strictly positive. Thus for sufficiently large $t$, there are constants $C_4, C_5$ such that
\begin{eqnarray}
C_4t^{1+2\delta} &\leq& t^{\epsilon}\theta(t)^{1-\epsilon} \nonumber\\
C_5t^{1+\frac{2\delta}{1-\epsilon}} &\leq& \theta(t) \nonumber
\end{eqnarray}
If $\delta$ is sufficiently close to zero, then this contradicts the
previous result that $\theta(t)$ goes to zero on a subsequence. If
$\delta$ is larger, then we use the fact that by equation
(\ref{rdeltabound}), $\theta(t)\geq C_5t^{r \delta}$. This implies the
original assumed lower bound is replaced by the larger lower bound
$t^{r\delta}$. Repeated iterations of this process shows that
$\theta(t)$ is bounded below by $t^{r^n\delta}$ for any $n$. Since
$r<1$, this reduces the situation to the $\delta$ close to zero case,
which led to a contradiction. Thus $\theta(t)$ can not be bounded
below by a function of the form $Kt^{\delta}$ for
$\delta\geq\frac{-\epsilon}{1-\epsilon}$. In particular, $\theta(t)$
can not be bounded below by a function of the form
$Kt^{\frac{-\epsilon}{1-\epsilon}}$ and so the desired subsequence
must exist. 
\end{proof}

The propagation observable $\gamma$ majorates the weight $\weakloc^{-\beta}$ for $\beta>\frac32$. This is a regularised version of the Morawetz estimate. As a consequence, the solution is space-time integrable in any bounded region, and must decay in that region. 

\begin{theorem}[Local Decay Estimate]
\label{T39}
\label{LocalDecay}
If $\beta>\frac{3}{2}$, then there is a positive constants $K_{LD}$ and such that \foruasoluinS, 
\begin{align}
\int_1^\infty \|(1+\rho_*^2)^{\frac{-\beta}{2}}u\|^2 dt
\leq& K_{LD} \sqrt{E}(\sqrt{E}+\|u(1)\|)
\end{align}
where $E=E[u_0,u_1]$. 
\end{theorem}
\begin{proof}
First a crude estimate on the norm growth of $\|\weakloc^{\frac{-\beta}{2}}u\|$ is made. In this paragraph, we use the notation $f\defin\weakloc^{\frac{-\beta}{2}}\leq1$. 
\begin{align}
\dt \|fu\|^2=&\dt\langle fu,fu\rangle\\
2\|fu\|\dt\|fu\|=&2\langle fu,fu\rangle\\
\dt\|fu\|\leq& \|f\dot{u}\|\leq\|\dot{u}\|\leq\sqrt{E[u(1)]} \label{eCrudeEstimateOnWeightedNorm}
\end{align}

Initially $\beta$ is restricted to $(\frac{3}{2},2)$, and the variable $\sigma=2(\beta-1)\in(1,2)$ is used. This $\sigma$ is distinct from the $\sigma$ used previously in the commutator argument and definition of $\gamma$. From the commutator estimates lemma \ref{lBasicMorawetzCommutatorWasLikeL19andL28}, there is a constants $C_2$ such that 
\begin{align}
\int_1^T C_{2}\|\weakloc^{\frac{-\beta}{2}}u\|^2 dt\label{39.1}
\leq& -\int_1^T\langle u,[\sum_{i=1}^3 H_i,\gamma]u \rangle dt
\end{align}
From integrating the Heisenberg relation, the integral on the right is bounded by 
\begin{align}
\int_1^T\dt(\langle u,\gamma\dot{u}\rangle-\langle\dot{u},\gamma u\rangle)dt
\leq& [\langle u,\gamma\dot{u}\rangle-\langle\dot{u},\gamma u\rangle]_{t=1}^T \\
\leq& [\dt\langle u,\gamma u\rangle-2\langle\dot{u},\gamma u\rangle]_{t=1}^T\\
\leq& [-2\langle\dot{u},\gamma u\rangle]_{t=1}^T\\
\leq& [\|\dot{u}\|\|\gamma u\|]_{t=1}+[\|\dot{u}\|\|\gamma u\|]_{t=T}
\end{align}
From energy conservation and the bound on $\|\gamma u\|$ in lemma \ref{L18.1}, these norms are bounded by 
\begin{align}
[\|\dot{u}\|\|\gamma u\|]_{t=1}+[\|\dot{u}\|\|\gamma u\|]_{t=T}
\leq&
C_1\sqrt{E}(\sqrt{E}+\|\weakloc^{\frac{-\sigma}{2}}u(1)\|)
\\
&+ C_2\sqrt{E}(\sqrt{E}+\|\weakloc^{\frac{-\sigma}{2}}u(T)\|) 
\label{39.2}
\end{align}

Since $\sigma\in(1,2)$, $q$ can be chosen in $(\frac{1}{\sigma}+\frac{1}{2},\frac{3}{2})\subset(1,\frac{3}{2})$. If $p$ is the conjugate exponent to $q$ and $\kappa\defin\frac{2}{p}$, then 
\begin{align}
\frac{1}{p}>&1-\frac{2}{3}=\frac{1}{3}\\
(\frac{2-\kappa}{2})q=&1\\
\sigma q>&\frac{\sigma +2}{2}
\end{align}

H\"older's inequality can be applied to the last term in equation \eqref{39.2} using conjugate exponent $p$ and $q$. Continuing to use $\|\bullet\|$ as the $L^2$ norm gives, 
\begin{align}
\|\weakloc^{-\frac\sigma2}u\|^2=&\int_\starman \frac{|u|^\kappa |u|^{2-\kappa}}{\weakloc^{\sigma}} d^3\mu \\
\leq& (\int_\starman |u|^{p\kappa}d^3\mu)^{\frac{1}{p}}(\int_\starman\frac{|u|^{(2-\kappa)q}}{\weakloc^{\sigma q}} d^3\mu)^{\frac{1}{q}}\\
\leq& \| |u|^{\frac{p\kappa}{2}} \|^{\frac{2}{p}}
\|\frac{|u|^{\frac{2-\kappa}{2}q}}{\weakloc^{\frac{\sigma q}{2}}} \|^\frac{2}{q}
\end{align}
The fact that $\sigma q>\frac{\sigma+2}{2}$ can now be used.
\begin{align}
\|\weakloc^{-\frac\sigma2}u\|
\leq& \|u\|^{\frac{1}{p}}  \|\weakloc^{\frac{-\sigma-2}{4}}u\|^{1-\frac{1}{p}}\label{40.1}
\end{align}

Since $\|\dot{u}\|^2 \leq E$, the norm of $u$ is controlled, $\| u(t)\|\leq \|u(1)\|+ tE^\frac12$. The computations from the beginning of \ref{39.1} to \ref{40.1} give an integral inequality for the weighted norm of $u$. 
\begin{equation}
\int_1^T \|\weakloc^{\frac{-\sigma-2}{4}}u\|^2dt\leq C_1(E[u]+\|u\|^2)+C_2E[u]^{\frac{1}{2p}}t^\frac{1}{p}\|\weakloc^{\frac{-\sigma-2}{4}}u\|^{1-\frac{1}{p}}
\end{equation}
Since $\dt\|\weakloc^{\frac{-\sigma-2}{4}}u\|$ is uniformly bounded by equation \eqref{eCrudeEstimateOnWeightedNorm}, and since $\frac{1}{p}\in(0,\frac{1}{3})$, the Gronwall's estimate lemma \ref{L38.2} applies, and there is a subsequence on which
\begin{equation}
t^{\frac{1}{p}}\|\weakloc^{\frac{-\sigma-2}{4}}u\|^{1-\frac{1}{p}}\rightarrow 0
\end{equation}
Using the calculations starting from the left hand side of equation \eqref{39.1}, which is monotonically increasing in time, to equation \eqref{40.1}, which is sequentially decreasing, the desired result holds for $\beta\in(\frac{3}{2},2)$. Since $\weakloc^{\frac{-\beta}{2}}|u|$ is monotonically decreasing in $\beta$, the Lebesgue dominated convergence theorem extends the result to all $\beta>\frac{3}{2}$. 
\end{proof}

To conclude this section, we apply this result to the conformal estimate. Since the weights in the local decay theorem will dominate any compactly supported function, if $u$ is a radial function, the trapping terms can be integrated in time to control the growth of the conformal charge and the weighted $L^6$ norm. 

\begin{corollary}
\label{cRadialDecayRate}
In the non super critical cases, $|Q|\leq M$, if \uasoluinS, and $u_0$ and $u_1$ are radial functions, then there is a constant, $C$, depending only on $\|u_0\|^2$, $ E[u_0,u_1]$, and $\confchrg[u_0,u_1]$ such that,  
\begin{align}
\left( \int |u|^6(t,\rho_*,\omega) F^3 r^{-4} \dmu )\right)^\frac{1}{6}
\leq& C t^\frac{-1}{3}
\end{align}
\end{corollary}
\begin{proof}
Since the initial conditions are radial, the solution $u$ will remain radial, and 
\begin{align}
\int W_L |\Lu|^2 d\mu = 0 .
\end{align}

Since the potential trapping term, $W$, is compactly supported, it can be dominated by the weights in the local decay result. This provides a bound on the conformal charge. 
\begin{align}
\confchrg[u(T),\dot{u}(T)]
\leq& \confchrg[u_0,u_1] + \int^t \int 2\tau W |u|^2 \dmu d\tau \\
\leq& \confchrg[u_0,u_1] +  t \int \int 2W |u|^2 \dmu d\tau \\
\leq& \confchrg[u_0,u_1] +  t E^\frac12 (E^\frac12 + \|u_0\|) 
\end{align}

From proposition \ref{pReductionToCompactlySupported}, 
\begin{align}
\int |u|^6(t,\rho_*,\omega) F^3 r^{-4} \dmu
\leq& C (E[u_0,u_1] + \confchrg[u(t),\dot{u}(t)]t^{-2}) \confchrg[u(t),\dot{u}(t)] t^{-4} \\
\leq& C t^{-2}
\end{align}
\end{proof}

\section{Angular Modulation}
\label{sAngularModulation}
To close the conformal estimate it is necessary to estimate one
angular derivative in $L^2$ localised near the photon sphere $r=\alpha$. Although this will not be possible, we will be able to bound fractional powers of the angular derivative. To do this, we introduce a propagation observable, which is an analogue of $\gamma$ but is rescaled, on each spherical harmonic, by a fractional power of $1-\SLap$. In this section, we control $(1-\Delta_{S^2})^\frac32$ in expectation value. 

Inspired by, but deviating from, the standard notation from Quantum Mechanics, 
$L$ is notationally used to refer to the operator square root of
$1-\Delta_{S^2}$; however, since an explicit decomposition of the
function space into eigenfunctions of $1-\Delta_{S^2}$ is known, at
this point it can still be thought of as a purely notational
convenience. A strictly positive version of $-\Delta_{S^2}$ is used so that inverse powers can be taken in section \ref{sPhaseSpaceAnalysis}. This differs slightly from $\tl$, which was used in the previous section but will not be used again. 

\begin{definition}
$L$ is defined to be the operator square root of $1-\Delta_{S^2}$. On each spherical harmonic this acts as multiplication by $\sqrt{1+l(l+1)}$. 
\end{definition}


\subsection{Angular Modulation and Initial Estimates}

A new operator $\gamma_{L^m}$ is introduced to give better estimates
near $\rho_*=0$. Because the contribution from $H_2=V$ involves no
powers of $L$, the local decay result allows us to ignore it, and the new 
multiplier is centred at the peak of $V_L$ uniformly in $L$. Instead, the rescaling parameter will vary with $L^m$. 

\begin{definition} The angularly modulated multiplier $\gamma_{L^m}$ is defined for $m\in[0,1/2]$, $b>0$, and $\sigma>1$, by
\begin{align}
g_{L^m}\defin&\int_0^{ b L^m\rho_*}(1+|\tau|)^{-\sigma} d\tau \\
=&\sum_0^\infty \int_0^{ b (1+l(l+1))^{\frac{m}{2}}\rho_*}(1+|\tau|)^{-\sigma} d\tau
\Pl\\
\gamma_{L^m}\defin&\frac{1}{2}(g_{L^m}\dr+\dr g_{L^m})
\end{align}
\end{definition}

We will take $\sigma\geq2$ since there is no improvement in varying $\sigma$ and it simplifies the estimates. 

\begin{lemma}
\label{lII2.1}
If $m\in[0,1/2]$, $\sigma\geq2$, and $b>0$, then for all $u\in\solset$, 
\begin{equation}
\|\gamma_{L^m} u\|^2\leq C_m(E[u]+\|\weakloc^{-1}u\|^2) .
\end{equation}
\end{lemma}
\begin{proof}
Since $\sigma\geq2$, $g_{L^m}$ is bounded in absolute value by
$C_\sigma$. The derivative is bounded on each spherical harmonic shell. 
\begin{align}
g_{L^m}'= L^m(1+bL^m|\rho_*|)^{-\sigma} 
\leq L^m(1+bL^m|\rho_*|)^{-2}
\leq C(L^mF^{-\frac{1}{2}}r^{-1}+\weakloc^{-1})
\end{align}
(this constant is independent of the spherical harmonic shell). The norm $\|g_{L^m}' u\|^2$ is
bounded by $\|LF^{\frac{1}{2}}r^{-1}u\|^2+\|\weakloc^{-1} u\|^2$ which are
bounded by the energy and local decay norm. 
\begin{align}
\|\gamma_{L^m}u\|^2 \leq&2(
\|g_{L^m}u'\|^2+\frac{1}{2}\|g_{L^m}'u\|^2)\\
\leq& C_m(E[u]+\|\weakloc^{-1}u\|^2) .
\end{align}
\end{proof}


\subsection{Direct Angular Momentum Bounds}
\label{ssDirectAngularModulationEstimates}

The family of multipliers $\gamma_{L^m}$ is not sufficient to
prove a local decay estimate for the angular energy $\langle u,W_L L^2
u\rangle$. Instead only an estimate for $\langle u, W_L
L^{\frac{3}{2}}u\rangle$ can be proven directly without phase space
analysis. The same method as in the proof of local decay is used, and
the best result is found by optimising over $m$. We track some of the additional positive terms more carefully to gain this improvement. 

Whenever a commutator involving a $\rho_*$ localisation or its derivatives, appears in a calculation, that term can be controlled by the previous local decay result. For this reason, we start by ignoring the contribution from $H_2$ and estimate it later. 

\begin{lemma}
\label{lAngModResultLikePreviouslII4.2}
If $m\in[0,\frac12]$, $\sigma\geq2$, $\epsilon>0$, and a continuous, compactly supported function $\tilde{\chi}$, then there is a choice of $b$, constants $C_1$, $C_2$ and $C_3$ and a smooth, compactly supported function $\chi$ such that for all $u\in\solset$, 
\begin{align}
\langle u,[H_1+H_3,\gamma_{L^m}] u\rangle
\geq& C_1\langle u', \frac{L^m}{(1+(\frac{L^m\rho_*}{2M}))^2}u'\rangle\\
&+ (1-\frac\epsilon3) \langle u, -gV_L' (L^2-1) u\rangle 
- C_2 \int b^2\rho_*^2 \chi u^2d\rho_* ,
\label{II3.1a} \\
\geq&  \epsilon \int L^{3m} \tilde\chi u^2 d\rho_* 
+ (1-\epsilon) \int -gV_L' (L^2-1) u^2 d\rho_*  
- C_3 \int b^2\rho_*^2 \chi u^2 d\rho_* 
\label{II3.1b}. 
\end{align}
\end{lemma}
\begin{proof}
We follow the same argument as before, but with the extra parameters from scaling by $L^m$. From the definition of $g$, we have 
\begin{align}
g'=& \frac{L^{m} b }{(1+L^{m}b |\rho_*|)^\sigma} , \\
g''=& \frac{-\sigma L^{2m}b ^2 \sgn \rho_*}{(1+L^{m}b |\rho_*|)^{\sigma+1}} , \\
g'''=& -2L^{2m}b ^2 \sigma \delta(\rho_*) + \sigma(\sigma+1)L^{3m} \frac{b ^3}{(1+L^{m}b |\rho_*|)^{\sigma+2}} .
\end{align}

Once again, we begin by controlling the $\delta$ function appearing in $g'''$. For this argument, $\chi$ denotes a compactly supported function, identically one in an interval of width $L^{-m}$ about the origin and supported on an interval of the same scale. 
\begin{align}
0
=&\int \dr (\rho_*\chi u^2) d\rho_* \\
=&\int \chi u^2 d\rho_* + \int \rho_*\chi'u^2d\rho_* + \int 2\rho_*\chi u\dr u d\rho_* 
\end{align}
We isolate the first term, note that the derivative of a characteristic function is bounded by $L^{2m}\rho_*$ times a (different) characteristic function (also denoted $\chi$), multiply by $L^{3m}$, and apply Cauchy-Schwartz, 
\begin{align}
\label{eGoodTermsControlCompact}
\int L^{3m} \chi u^2 d\rho_*
\leq& \int L^{5m} \rho_*^2\chi u^2 d\rho_* 
+ 2\left(\int L^{5m} \rho_*^2\chi^2u^2d\rho_*\right)^\frac12 \left(\int L^{m} \chi (\dr u)^2 d\rho_*\right)^\frac12 .
\end{align}
We control the value of the function at $\rho_*=0$ using a similar method, integrating only on the positive axis to get the value at the origin. 
\begin{align}
u(0)^2
=& -\int_0^\infty \dr (\chi u^2) d\rho_* \\
=& -\int_0^\infty \chi' u^2 d\rho_* - \int_0^\infty 2\chi u\dr u d\rho_* \\
\leq& C \left( \int \chi u^2 d\rho_* \right)^\frac12
\left(\int L^{4m} \rho_*^2 \chi u^2 d\rho_* + \int \chi (\dr u)^2 d\rho_* \right)^\frac12 .
\end{align}
We now multiply by $L$ to the power $2m$. 
\begin{align}
L^{2m}u(0)^2
\leq& C \left( \int L^{3m}\chi u^2 d\rho_* \right)^\frac12
\left(\int L^{5m} \rho_*^2 \chi u^2 d\rho_* + \int L^{m}\chi (\dr u)^2 d\rho_* \right)^\frac12 .
\end{align}

Applying equation \eqref{eGoodTermsControlCompact} to control the weighted $L^2$ norms, we have, 
\begin{align}
L^{2m} u(0)^2
\leq& C_0 (\int L^{5m}\rho_*^2\chi u^2 d\rho_* 
+ \int L^{m} \chi (\dr u)^2 d\rho_* ) .
\label{eControlDeltaAM}
\end{align}

We now turn to controlling non compactly supported weights coming from $g'''$. 
\begin{align}
0=&\int \dr (g''u^2) d\rho_* \\
=&\int g''' u^2 d\rho_* + \int 2 g''u \dr u d\rho_* \\
L^{3m}b^3 \int \frac{\sigma(\sigma+1)}{(1+b L^{m}|\rho_*|)^{\sigma+2}} d\rho_*
\leq& 2\sigma L^{2m}b ^2u(0)^2 \\
&+ \left(\int L^{3m}b ^3 \frac{\sigma(\sigma+1)}{(1+b L^{m}|\rho_*|)^{\sigma+2}} u^2 d\rho_*\right)^\frac12 \left(\int\frac{2\sigma}{\sigma+1} b L^{m}\frac{2(\dr u)^2}{(1+b L^{m}|\rho_*|)^\sigma} d\rho_* \right)^\frac12 
\label{eResultOfAngModIBP}
\end{align}
Substituting equation \eqref{eControlDeltaAM} and noting $\chi\leq C(1+bL^m|\rho_*|)^\sigma$, 
\begin{align}
L^{3m} b^3 \int &\frac{\sigma(\sigma+1)}{(1+b L^{m}|\rho_*|)^{\sigma+2}} u^2 d\rho_*\\
\leq& \sigma b  C_0 \left(\int b \threeprterms \rho_*^2\chi u^2 d\rho_* 
+ \int bL^{m} \frac{2(\dr u)^2}{(1+b L^{m}|\rho_*|)^\sigma} d\rho_*\right)\\
& +\frac12\left(\int L^{3m}b ^3 \frac{\sigma(\sigma+1)}{(1+b L^{m}|\rho_*|)^{\sigma+2}} u^2 d\rho_* + \int \frac{2\sigma}{\sigma+1} L^{m}b \frac{2(\dr u)^2}{(1+b L^{m}|\rho_*|)^\sigma} d\rho_* \right ) \\
L^{3m}b ^3 \int &\frac{\sigma(\sigma+1)}{(1+b L^{m}|\rho_*|)^{\sigma+2}} u^2 d\rho_*\\
\leq& 2\sigma b  C_0 \int \threeprterms b \rho_*^2\chi u^2 d\rho_* 
+ (\frac{2\sigma}{\sigma+1}+2\sigma b C_0 )\int b L^{m}\frac{2(\dr u)^2}{1+bL^{m} |\rho_*|)^\sigma} d\rho_* 
\label{eInverseCubeControlled}
\end{align}

The commutator term can now be estimated
\begin{align}
\int 2g'&(\dr u)^2 -\frac12 g''' u^2 - gV_L'(L^2-1) u^2 d\rho_* \\
\geq& \int L^{m}b \frac{2(\dr u)^2}{1+b |\rho_*|} d\rho_* 
+L^{2m}b ^2\sigma u(0)^2 \\
&- \int L^{3m}b ^3 \frac{\sigma(\sigma+1)}{2(1+b |\rho_*|)^{\sigma+2}} u^2 d\rho_* 
- \int gV_L' (L^2-1)u^2 d\rho_* \\
\geq& \int (\frac{2}{\sigma+1}-2\sigma b C_0)L^{m} b  \frac{(\dr u)^2}{(1+b |\rho_*|)^\sigma} d\rho_* 
+ L^{2m} b^2\sigma u(0)^2\\
& - \int \sigma L^{5m} b^2 C_0 \rho_*^2 \chi u^2 d\rho_* 
- \int gV_L' (L^2-1)u^2 d\rho_* .
\end{align}
We choose $b$ sufficiently small so that $2/(\sigma+1)-2\sigma b C_0>1/2$ and $-(\epsilon/3) \rho_*V_L'\chi > \sigma b^2 \rho_*^2 \chi$. This choice of $b$ can be computed explicitly from $C_0$ and $V_L'(0)$, but it will be independent of the spherical harmonic. Since $g$ vanishes like $L^{m}\rho_*$ in the support of $\chi$, under these conditions, $-(\epsilon/3)gV_L' L^2 > \sigma L^{5m} b^2 \rho_*^2 \chi$. This gives us sufficient control, except in the case $L^2-1=0$, when we must subtract off an additional localisation term. With these choices, we have 
\begin{align}
\int 2g'&(\dr u)^2 -\frac12 g''' u^2 - gV_L' (L^2-1)u^2 d\rho_* \\
\geq& c \int L^{m} b  \frac{(\dr u)^2}{(1+b L^m|\rho_*|)^\sigma} d\rho_* 
+ cL^{2m} b^2u(0)^2 
- (1-\frac\epsilon3 )\int gV_L' (L^2-1) u^2 d\rho_*  
- C \int b^2\rho_*^2 \chi d\rho_* . 
\end{align}

We ignore the $u(0)^2$ term and divide the remaining positive terms into three pieces. A large piece will be $1-\epsilon$ times $-gV_L' (L^2-1)$. Another piece will be the $(u')^2$ term and $\epsilon/3$ times $-gV_L'(L^2-1)$. In a compact set of width $L^{-m}$, by equation \eqref{eGoodTermsControlCompact}, this controls $L^{3m}$, up to a compactly supported term coming from the $L^2-1=0$ case. The remaining $\epsilon/3$ factor of $-gV_L'$, we use outside $|\rho_*|\leq L^{-m}$. In this outer region, but inside a compact set independent of $L$, $-gV_L' (L^2-1)$ dominates any continuous, compactly supported function, $\tilde\chi$, by a factor of $L^{2-m}>L^{3m}$, again with an extra contribution required for the $L^2-1=0$ case. Here, we use the assumption $m\leq1/2$ again. From this, 
\begin{align}
\int 2g'&(\dr u)^2 -\frac12 g''' u^2 - gV_L' u^2 d\rho_* \\
\geq& c \epsilon \int L^{3m} \tilde\chi u^2 d\rho_* 
+ (1-\epsilon) \int gV_L' (L^2-1) u^2 d\rho_*  
- C \int b^2\rho_*^2 \chi d\rho_* . 
\end{align}
Recall $\tilde\chi$ denotes an arbitrary continuous, compactly supported function, and $\chi$ is supported on $|\rho_*| <C L^{-m}$. 
\end{proof}

The next lemma says the contribution from $H_2$ is a local decay term. 

\begin{lemma}
\label{lII4.1}
If $m\in[0,1/2]$, $\sigma\geq2$, and $b>0$, then there is a positive constant $C$ such that for all $u\in\solset$, 
\begin{equation}
|\langle u,[H_2,\gamma_{L^m}]u\rangle| \leq C\langle u,
 \weakloc^{-2}u\rangle
\label{II4.1}
\end{equation}
\end{lemma}
\begin{proof}
From lemma \ref{lPotentialTrappingCompactlySupported}, the derivative of the potential is $-V'=2Fr^{-7}P_Q(r)$ where $P_Q$ is a cubic polynomial in $r$. 

Since for $\rho_*\rightarrow\infty$,
$r\sim\rho_*$, and for $\rho_*\rightarrow-\infty$, $F$ decays
exponentially, there is a constant $C$ so that
$|V'|<C\weakloc^{-2}$. Since $|g_{L^m}|$ is bounded by a constant and
$[H_2,\gamma_{L^m}]=-g_{L^m}V'$, the result holds. 
\end{proof}

The previous two lemmas shows that $\gamma_{L^m}$ is a propagation observable which majorates powers of $L$. That is, ignoring terms which are space- time integrable, the commutator is bounded below by the product of powers of $L$ and localisation functions. These localisation functions are functions of $L^m\rho_*$. In the following theorem, the Heisenberg-type relation with $\gamma_{L^m}$ is integrated in time. The theorem shows that $\gamma_{L^m}$ is a bounded operator on the energy space. 

The localising function $\tilde\chi$ is included to restrict attention to the region $\rho_*$ bounded but on an $L$ independent scale. In the theorem, taking the optimal value $m=\frac{1}{2}$ gives an estimate for $\langle u,\tilde\chi L^{\frac{3}{2}}u\rangle$. The estimate for other values of $m$ will be used later, in the phase space analysis. 

\begin{theorem}\footnote{This is the equivalent of theorem 7.8 from the previous draft. }
\label{tII5.1}
If $m\in[0,\frac12]$, then there is a constant $C$ such that for all $u\in\solset$ which satisfy the linear wave equation, $\ddot{u}+\linH u=0$, 
\begin{align}
\int_1^\infty \langle u', \frac{L^m}{(1+(L^m\rho_*)^2}u'\rangle
+ \langle u, -gV_L' (L^2-1) u\rangle dt 
\leq& C(E[u]+\|u(1)\|^2) \label{II5.1} \\
\int_1^\infty \langle u,L^\frac{3}{2} \chia u\rangle 
\leq& C(E[u]+\|u(1)\|^2) \label{II5.2}
\end{align}
\end{theorem}
\begin{proof}
Applying the Heisenberg identity to the operator $\gamma_{L^m}$ gives
\begin{equation}
\Dt(\langle\dot{u},\gamma_{L^m}u\rangle-\langle
u,\gamma_{L^m}\dot{u}\rangle) = \langle u, \linHexpanded u\rangle
\end{equation}
Using lemma \ref{lII2.1} and the fact that $\gamma_{L^m}$ is
antisymmetric with respect to the $L^2(\starman)$ norm when acting on
$C^\infty$ functions, this equation can be integrated to get an
estimate on the commutators.
\begin{align}
\int_1^t\langle
u,[\linHexpanded,\gamma_{L^m}]u\rangle d\tau\leq&2\|\dot{u}\|\|\gamma_{L^m}
u\| |_1^t\\
\leq&C(E[u]+\|\weakloc^{-1} u(1)\|^2 +\|\weakloc^{-1} u(t)\|^2)
\end{align}

By the local decay result, theorem \ref{T39}, and lemma \ref{lII4.1}, 
\begin{align}
\int_1^\infty |\langle u,[H_2,\gamma_{L^m}]u\rangle| d\tau
\leq&\int_1^\infty \|\weakloc^{-1}u\|^2d\tau\leq C(E[u]+\| u(1)\|^2) ,\\
\int_1^\infty \int_\starman b^2\rho_*^2 \chi u^2 d\rho_* d\tau
\leq&\int_1^\infty \|\weakloc^{-1}u\|^2d\tau\leq C(E[u]+\| u(1)\|^2) ,\\
\end{align}

From the first estimate for $[H_1+H_3,\gamma_{L^m}]$, equation \eqref{II3.1a}, the following estimate holds on a sequence of times:
\begin{align}
\int_1^t \langle u', \frac{L^m}{(1+(\frac{L^m\rho_*}{2M}))^2}u'\rangle
&+ \langle u, -gV_L' (L^2-1) u\rangle d\tau \\
\leq& C(E[u]+\|u(1)\|^2+\|\weakloc^{-1} u(t)\|^2)\\
&+C(E[u]+\|\weakloc^{-1}u(1)\|^2)
\end{align}

Since the left hand side is monotonically increasing and the right hand side is decreasing on a sequence, by the local decay estimate, theorem \eqref{LocalDecay}, equation \eqref{II5.1} holds for all $t$. 

Equation \eqref{II5.2} is proven by taking $m=\frac{1}{2}$, and
applying the same argument using the estimate in \eqref{eGoodTermsControlCompact}. 
\begin{equation}
\int_1^\infty\langle u,L^\frac{3}{2}\chia u\rangle \leq
C(E[u]+\|u(1)\|^2)
\end{equation}
\end{proof}

This has proven an estimate for $\langle u,L^\frac32\chia u\rangle$. In the later phase space analysis, it will be necessary to estimate mixed derivative norms involving $\dr$ and functions of $L$ and $\rho_*$. The following corollary is a result of theorem \ref{tII5.1}. It controls the expectation value of localised $\dr$ terms. This could have been proven directly from the Morawetz estimate with out angular modulation. 

\begin{corollary}
\label{cII6.1}
There is a positive constant $C$ such that for all $u\in\solset$ which satisfy the linear wave equation, $\ddot{u}+\linH u=0$, 
\begin{align}
\int_1^\infty \langle u', \frac{1}{1+\rho_*^2} u'\rangle d\tau \leq C(E(u)+\|u(1)\|^2) .
\end{align}
In particular, if $u$ is a solution of the linear wave equation, $\ddot{u}+\linH u=0$, and $\chi_1$ is a smooth, compactly supported function, then there is a $C$ such that, 
\begin{align}
\int_1^\infty \langle u',\chi_1 u'\rangle d\tau \leq& C(E(u)+\|u(1)\|^2) , \\
\int_1^\infty \| \dr \chi_1 u\|^2  d\tau \leq& C(E(u)+\|u(1)\|^2) .
\end{align}
\end{corollary}
\begin{proof}
Using the first term in the first estimate of theorem \ref{II5.1} with $m=0$ gives this result. 
\end{proof}

We also have control of the full two angular derivatives in any compactly supported region, with a weight which vanishes quadratically at the photon sphere, $\rho_*=0$. 

\begin{corollary}
\label{cReductionToLonPhotonSphere}
For any continuous, compactly supported function $\tilde\chi$, there is a positive constant $C$ such that for all $u\in\solset$ which satisfy the linear wave equation, $\ddot{u}+\linH u=0$, 
\begin{align}
\int_1^\infty \langle u, \rho_*^2\tilde{\chi} L^2 u\rangle d\tau \leq C(E[u] + \|u(1)\|^2) .
\end{align}
\end{corollary}
\begin{proof}
Again using the first estimate of theorem \ref{II5.1} with $m=0$, this time using the second term, and noting that any continuous, compactly supported function which vanishes quadratically at the photon sphere, in particular $\rho_*^2\chi$, is dominated by $-gV_L'$. On this compactly supported region, the difference between $L^2-1$ and $L^2$ is order one, so it's a local decay term, and also space time integrable. 
\end{proof}

\section{Phase Space Analysis}
\label{sPhaseSpaceAnalysis}

The main conclusion of this section is the phase space induction theorem. It does not control the time integral of the angular energy, but instead loses an arbitrarily small power of $L$. The bounds are still in terms of the initial energy and the $L^2$ norm. 

To prove this result, we introduce the operators $\Gammanm$ to majorate $L^\frac{\frac32+m}{2}$. Each $\Gammanm$ is a phase space localised version of $\gamma_{L^m}$. Our localisation will use the radial variable and the radial derivative rescaled using the powers of $L$. We call these rescaled quantities the phase space variables. Because of the localisation, the commutator of $\Gammanm$ and $\linH$ dominates $L^\frac{\frac32+m}{2}$ only in regions of phase space. This is where our definition of propagation observable differs from the standard one. In the standard definition, the commutator of a propagation observable must dominate the operator it is said to majorate; whereas, we allow this domination to occur in a region of phase space. Under our definition, but not the standard one, $\Gammanm$ majorates $L^\frac{\frac32+m}{2}$. 

The definition of majoration always permits ``lower order terms''. For $\Gammanm$, these will be $L^\frac{1+2n}{2}$. So that the lower order terms truly are lower order, we require $0\leq n\leq m \leq\frac12$ and $n<\frac12$. By combining various $m$ and $n$ choices, we eventually show that $\|L^{1-\varepsilon} \chia u\|^2$ is integrable in time. 

In subsection \ref{ssStructure}, {\it Structure of the Paper}, there is already an outline of this section, which we repeat in brief. Subsection \ref{ssPhaseSpaceVariablesEtc} introduces the phase space variables and $\Gammanm$. Subsection \ref{ssCommutatorExpansions} introduces the commutator expansion theorem, which allows is to expand commutators of functions of the phase space variables, and several lemmas for common cases. Subsection \ref{ssPhaseSpaceLocalisedEstimates} computes the commutators to show $\Gammanm$ dominates powers of $L$ with localisation in both phase space variables. Subsection \ref{ssDerivativeLocalisedEstimates} combines these estimates to prove estimates which are localised in the rescaled radial derivative alone. Finally, subsection \ref{ssPhaseSpaceInduction} combines these in a finite induction to eliminate all phase space localisation. This is sufficient to prove pointwise in time $L^6$ decay.

%
%

\subsection{Phase Space Variables, Localisation, and Multipliers}
\label{ssPhaseSpaceVariablesEtc}

Previously, we introduced a radial variable which was rescaled by the angular momentum. Now, we also introduce a rescaled radial derivative. We refer to these as the phase space variables $\qm$ and $\pn$ respectively. These will be rescaled separately with different parameters $m$ and $n$. By varying these parameters, we will be able to control powers of $L$ in various regions of phase space. 

\begin{definition}
The phase space variables are 
\begin{align}
\qm =& L^m \rho_* \\
\pn =& -i L^{n-1} \dr
\end{align}
\end{definition}

Corollary \ref{cReductionToLonPhotonSphere} provides control on the angular energy away from $\rho_*=0$, so that the angular component of the energy only needs to be controlled in a compact set around $\rho_*=0$. This is the region near $r=\alpha$. We introduce the weight $\chia$ to work in this region; although, the particular choice of $\chia$ will not matter. 

\begin{definition}
The function $\chia$ is defined to be a smooth, compactly supported, radial function, such that $\chia\leq1$, and $\chia=1$ in a neighbourhood of $\rho_*=0$. 
\end{definition} 

For the phase space variables, we have a notion of lower order, involving our choice of $\chia$. 

\begin{notation}
The notation $\Oint$ denotes functions of time which are integrable, and for which the $L^1$ norm in time is bounded by a combination of constants, the energy of $u$, and its initial $L^2(\starman)$ norm. By theorem \ref{LocalDecay}, corollary \ref{cII6.1}, and corollary \ref{cReductionToLonPhotonSphere}, this includes local decay norms, localised radial derivatives, and angular derivatives with vanishing weight at $\rho_*=0$. 

The notation $\bndd$ denotes an arbitrary bounded operator, just as the notation $C$ denotes an arbitrary constant. As with arbitrary constants, the value of $\bndd$ may vary from line to line in an argument. The notation $\bndd_i$ is used to refer to a bounded operator which is referred to later in the same argument. 

For $n\leq\frac12$, the notation $\ErrorTermsn$ denotes inner products and norms
 of $u$ which are either of the following form
\begin{equation}
\|L^\frac{1+2n}{2}\chia u\|^2, \Oint , 
\end{equation}
or bounded by norms of this form. Note that for $q\leq 1+2n$, this includes any 
term of the form
\begin{equation}
\|L^\frac{q}{2} \bndd\chia u\|^2 .
\end{equation} 
Note that in particular if $m\leq\frac12$, then $m+n\leq \frac12+n$
and 
\begin{equation}
\|L^{m+n}\bndd\chia u\|^2=\ErrorTermsn .
\end{equation}
\end{notation}

To begin with, it is not possible to use sharp cut off functions in $\qm$ and $\pn$, so more regular functions with infinite support must be used. The rate of decay for these functions is important, and various secondary functions derived from the also appear. The localisations occur in families, with $\qmnfn$ referring to $\qm$ localisation, $\pnnfn$ referring to $\pn$ localisation in a wide band around zero, and $\pnifn$ referring to $\pn$ localisation in a narrow band away from $0$. 

In our estimate on the angular energy, we lose powers of $L$ from our choice of localisation. We lose a factor of $L^\epsilon$ from our choice of decay for $\pnna{x}$ and a factor of $L^\delta$ from our choice of width for $\pni{x}$. 

\begin{definition}
The weight $g$ is defined by
\begin{align}
g(x)=&\int_0^{bx} \frac{1}{(1+|\tau|)^2} d\tau 
\end{align}
where $b$ is chosen as in lemma \ref{lII2.1}. 
The near and far localisation functions for $\qm$ are 
\begin{align}
\qmn{x} =& \frac{b}{(1+b|x|)^2} , \\
\qmnt{x}=& \frac{1}{1+x^2}. \\
\qmf{x} =& x \sqrt{\frac{g(x)}{x}} .
\end{align}
From this, $\qmn{\qm}$, $\qmnt{\qm}$, and $\qmf{\qm}$ are defined by the spectral theorem. 
\end{definition}

\begin{definition}
The smooth, near localisation for $\pn$ and two functions derived from it are defined by
\begin{align}
\pnna{x} =& (1+x^2)^{-\frac{1-\epsilon}{4}} \\
\pnnb{x} =& x \pnnaWITHADERIV{x} \\
\pnnc{x} =& \sqrt{\pnna{x}(\pnna{x}+2x\pnnaWITHADERIV{x})} .
\end{align}
(The definition of $\pnnc{x}$ requires the argument of the square root to be positive, which is shown in lemma \ref{lYnwelldef}.) 

We use the notation $\chi(A,x)$ to denote the characteristic function of $A$ which is $1$ if $x\in A$ and $0$ otherwise. 

The sharp, near localisation and sharp, interval localisation for $\pn$ are 
\begin{align}
\pnn{x} =& \chi([0,1],x) \\
\pniWITHl{x} =& \chi([l^{-\delta},1],x) .
\end{align}
Smooth, interval localisations will be defined in section \ref{ssDerivativeLocalisedEstimates}. The functions $\pnn{x}$ and $\pniWITHl{x}$ are extended as even functions. 

From these, $\pnna{\pn}, \pnnb{\pn}, \pnnc{\pn}, \pnn{x}$ and $\pni{\pn}$ can be defined by the spectral theorem. 
\end{definition}

To prove a phase space localised version of the Morawetz estimate, we introduce the operator $\Gammanm$. The parameters $n$ and $m$ are restricted to $0\leq n \leq m \leq\frac12$, because this is the range on which an estimate can be proven. 

\begin{definition}
For $0\leq n \leq m \leq \frac12$, the phase space induction multiplier is defined to be 
\begin{align}
\Gammanm 
=& \chia \pnna{\pn} L^{n-\epsilon} \gamma_{L^m} \pnna{\pn} \chia \\
=& \frac12 \chia \pnna{\pn} i \left(g(\qm)\pn + \pn g(\qm)\right) \pnna{\pn} \chia
\end{align}
\end{definition}

We now prove some preliminary results, starting with some estimates on the localisations. 

\begin{lemma}
\label{lYnwelldef}
\label{lIIqLcontrolsGamma}
For all $\epsilon>0$, $n\in[0,\frac12]$, and $v\in\solset$
\begin{align} 
\|\pnnc{\pn} v\|\geq& \epsilon^\frac12  \|\pnna{\pn} v\| . \\
\|L^{1-\epsilon}\pn\pnna{\pn}^2 v\| \leq& \|L^{ \frac{1+n\epsilon-2\epsilon}{1-\epsilon}}v\| + \|\dr v\| 
\label{IIqLcontrolsGamma}
\end{align}
\end{lemma}
\begin{proof} 

The function $\pnnc{x}$ is computed explicitly. 
\begin{align}
\pnnc{x}^2=&\pnna{\pn} (\pnna{\pn}+2x\pnnaWITHADERIV{x})\\
=&(1+x^2)^{-\frac{1-\epsilon}{4}}(1+x^2)^{-\frac{5-\epsilon}{4}}((1+x^2)-2\frac{1-\epsilon}{4}2x^2)\\
=&(1+x^2)^{-\frac{6-2\epsilon}{4}}(1+\epsilon x^2)\\
\geq & \epsilon (1+x^2)^{-\frac{1-\epsilon}{2}}\\
=& \epsilon \pnna{\pn}^2
\end{align}
Since $\pnna{\pn} (\pnna{\pn}+2x\pnnaWITHADERIV{x})$ is strictly positive, the
square root function, $\pnnc{x}$ is well defined. Since
$\pnnc{x}^2\geq\epsilon \pnna{\pn}^2$, the spectral theorem also
determines that $\|\pnnc{\pn} v\|\geq \epsilon^\frac12  \|\pnna{\pn} v\|$. 

Equation \eqref{IIqLcontrolsGamma} is proven first for Schwartz class data. For this data, in a representation where $L$ and $-i\dr$ both act as multiplication operators, $\pnna{\pn}^2\leq |\pn|^{-1+\epsilon}=|-i\dr L^{n-1}|^{-1+\epsilon}$. This can be used to start estimating the norm of $L^{n-\epsilon}\dr\pnna{\pn}^2 v$. This gives equation \eqref{IIqLcontrolsGamma}. 
\begin{align}
\| L^{1-\epsilon}\pn\pnna{\pn}^2 v \|
\leq& \|L^{1-\epsilon} \pn |\pn|^{-1+\epsilon} v\|\\
\leq& \|L^{1-\epsilon} |\pn|^\epsilon v\|\\
\leq& \|L^{1+n\epsilon-2\epsilon}|-i\dr|^\epsilon v\|
\end{align}
Continuing in a representation where both $L$ and $-i\dr$ act as multiplication operators, the identity, for positive $f$ and $g$, $fg \leq f^\frac{1}{1-\epsilon}+g^\frac1\epsilon$ can be applied with $f=L^{1+n\epsilon-2\epsilon}$ and $g=|-i\dr|^{\epsilon}$. 
\begin{align}
\| L^{1-\epsilon}\pn\pnna{\pn}^2 v \| \leq& \|L^{ \frac{1+n\epsilon-2\epsilon}{1-\epsilon}} v\| + \|\dr v\|
\end{align}
\end{proof}

Since $\Gammanm$ is a product of operators, there is a Leibniz (or product) rule for computing the commutator with another operator. It is possible to further simplify this by expressing a sum of complex conjugates as a real part. The main application will be when the other operator is $\linH$ or one of the subterms $H_i$. 

\begin{lemma}
\label{lIIp3.1}
For $u\in\solset$, $0\leq n\leq m\leq\frac12$, and any self-adjoint operator $G$ which commutes with $L$
\begin{align}
\label{IIp3.1}
\langle u,[G,\Gammanm]u\rangle=&\langle u,\chia\pnna{\pn}[G,\gamma_{L^m}]L^{n-\epsilon}\pnna{\pn}\chia u\rangle\\
&+\langle u,\chia[G,\pnna{\pn}]\gamma_{L^m}L^{n-\epsilon}\pnna{\pn}\chia u\rangle
\\ &+\langle u,\chia\pnna{\pn}\gamma_{L^m}L^{n-\epsilon}[G,\pnna{\pn}]\chia u\rangle\\
&+\langle u,\chia\pnna{\pn}\gamma_{L^m}L^{n-\epsilon}\pnna{\pn}[G,\chia] u\rangle
\\ &+\langle u,[G,\chia]\pnna{\pn}\gamma_{L^m}L^{n-\epsilon}\pnna{\pn}\chia u\rangle
\end{align}
\end{lemma}
\begin{proof} This is simply an application of the Leibniz rule for
  commutators. 
\end{proof}

%
%
\subsection{Commutator Expansions}
\label{ssCommutatorExpansions}

To apply the Heisenberg like relation, it is necessary to expand commutators involving localisation in $\pn$. This is done through a version of the commutator expansion lemma previously used in scattering theory\cite{Soffer}. We consider the special case of this expansion for commutators involving localisation in the phase space variables $\pn$ and $\rho_*$ or $\qm$. These expansions are as finite order power series with an error term which involves the Fourier transform of a $k$th order derivative of one of the localising functions. 

First, the adjoint is defined. These are iterated commutators and will appear in the commutator expansion, both as terms in the finite power series expansion and as a term in the remainder. 

\begin{definition}
For two operators $H$ and $A$, the $k$th commutator of $A$ with respect to $H$, $\ad_A^k(H)$ is defined recursively by
\begin{align}
\ad_A^1(H)&=[H,A]\\
\ad_A^k(H)&=[\ad_A^{k-1}(H),A]
\end{align}
Initially the commutator $[H,A]$ is defined only as a form on the domain
of $A$ intersect the domain of $H$ by the formula
\begin{align}
\langle u,[H,A]u\rangle=&\langle Hu,Au\rangle-\langle Au,Hu\rangle
\end{align}
If $\ad_A^{k-1}(H)$ extends to a bounded operator, then $\ad_A^k(H)$
is defined on the domain of $A$. 
\end{definition}

The commutator expansion theorem expands the commutator of an operator, $F_1$ , with a function of a self-adjoint operator, $F_2(A)$. The expansion is as a power series in the adjoint and has remainder involving the $L^1$ norm of Fourier transform of $F_2^{[k]}$. 

\begin{theorem}[Commutator Expansion Theorem]
\label{tIIp2.1}
If $n>0$ is an integer, $A$ is a self-adjoint operator, $F_1$ is a self-adjoint operator satisfying
\begin{enumerate}
\item for $1\leq k\leq n$, $\ad_A^k(F_1)$ extends to a bounded operator,
\end{enumerate}
and $F_2(x)$ is a smooth function satisfying
\begin{enumerate}
\setcounter{enumi}{1}
\item $\|\Fourier{F_2^{[n]}}\|_1\leq\infty$ ,
\label{CommutatorExpansionCondition}
\end{enumerate}
then if $[F_1,F_2(A)]$ is defined as a form on the domain of $A^n$, 
\begin{align}
\label{IIp2.2}
[F_1,F_2(A)]=\sum_{k=1}^{n-1}\frac{1}{k!} F_2^{[k]}(A) \ad_A^k(F_1) +R_n
\end{align}
in the form sense with the remainder $R_n$ satisfying
\begin{align}
\|R_n\|_\opnorm\leq& C\|\Fourier{F_2^{[n]}}\|_1 \|\ad_A^n(F_1)\|_\opnorm
\end{align}
Consequently, $[F_1,F_2(A)]$ defines an operator on the domain of $A^{n-1}$. 
\end{theorem}
\begin{proof}
This is proven in \cite{Soffer}. 
\end{proof}

The commutator expansion theorem is now specialised to the localising functions in the phase space variables. In terms of the previous theorem,  $A=\pn=-i\dr L^{n-1}$ and $F_1$ is a function of $\rho_*$ or of $\qm$. 

There are two applications of this lemma. The first application is the expansion of the commutator of the Hamiltonian, $\linH$, and the phase space induction multiplier, $\Gammanm$, as a power series in the adjoint. The remainder is a higher order term, in the sense that the commutator minus the power series expansion can be multiplied by a positive power of $L$ and still remains a bounded operator. The second application is using the expansion with no terms and using the remainder to directly estimate the entire commutator as a higher order term. This second application is expanded upon in the lemmas which follow this one. 

\begin{lemma}[Commutator order reduction lemma]
\label{lCommutatorOrderReduction}
\label{lIIp2.2}
Let $k$ be a positive integer, $F_1(x)$ be a smooth function satisfying
 \begin{enumerate}
\item for $1\leq j\leq k$, $\|F_1^{[j]}\|_\infty\leq\infty$,
\end{enumerate}
and $F_2(x)$ be a smooth function satisfying
\begin{enumerate}
\setcounter{enumi}{1}
\item for $1\leq j\leq k-1$, $\|F_2^{[j]}\|_\infty\leq\infty$, 
\item $\|\Fourier{F_2^{[k]}}\|_1\leq\infty$
\end{enumerate}
then there is a constant, $C_k$, depending only on $k$, such that
\begin{align}
L^{k(1-n)}[F_1(\rho_*),F_2(\pn)]&=\sum_{j=1}^{k-1}L^{k(1-n)}F_2^{[j]}(\pn)(iL^{n-1})^{j}F_1^{[j]}(\rho_*)+R_k\\
\|R_k\|_{\opnorm}&\leq C_k \|\Fourier{F_2^{[k]}}\|_1 \|F_1^{[k]}\|_\infty
\label{bnddCommutatorOrderReductionNom}
\end{align}
and 
\begin{align}
L^{k(1-m-n)}[F_1(\qm),F_2(\pn)]&=\sum_{j=1}^{k-1}L^{k(1-m-n)}F_2^{[j]}(\pn)(iL^{m+n-1})^{j}F_1^{[j]}(\rho_*)+R_k\\
\|R_k\|_{\opnorm}&\leq C_k \|\Fourier{F_2^{[k]}}\|_1 \|F_1^{[k]}\|_\infty
\label{bnddCommutatorOrderReductionWithm}
\end{align}

\end{lemma}
\begin{proof}

The main idea in this proof is to apply the commutator expansion theorem, theorem \ref{tIIp2.1}. 

The $m=0$ case is a special case of the formula for general $m\geq0$. 

All quantities in equation \eqref{bnddCommutatorOrderReductionWithm} are composed of powers of $L$ and functions of operators which commute with $L$. Therefore, they preserve the spherical harmonic decomposition, and it is sufficient to prove the power series expansion of the commutator on each spherical harmonic. We will use $\tl$ to denote the action of $L$ on each spherical harmonic, rather than the standard spherical harmonic index. 

On a fixed spherical harmonic, take $x=\qm=\tl^m \rho_*$ and $A=\pn=-i\tl^{n-1}\dr$ and apply the commutator expansion theorem. 
\begin{align}
[F_1(\qm),F_2(\pn)] 
=& \sum_{j=1}^{k-1} \frac{1}{j!} F_2^{[k]}(\pn) \ad_\pn^j(F_1(\qm)) + R_k \\
\|R_k\| \leq& C_k \| \Fourier{F_2^{[k]}} \|_1 \| \ad_\pn^k(F_1) \|_\infty
\end{align}

The adjoint can be rewritten in terms of $\tl$, $\rho_*$, and $\dr$. 
\begin{align}
\ad_\pn^j(F_1(\qm))
=& [\ad_\pn^{j-1}(F_1(\qm)),\pn]\\
=& [\ad_\pn^{j-1}(F_1(\tl^m\rho_*),-i\tl^{n-1}\dr]\\
=& i\tl^{n-1} \dr \ad_\pn^{j-1}(F_1(\tl^m\rho_*))
\end{align}
Taking a derivative with respect to $\rho_*$ of a function of $\tl^m\rho_*$ will introduce an additional factor of $\tl^m$. For $j=1$, we have
\begin{align}
\ad_\pn^1(F_1(\qm))
=& i\tl^{n-1} \dr F_1(\tl^m\rho_*)\\
=& i\tl^{m+n-1} F_1'(\tl^m\rho_*)
\end{align}
Each adjoint acts as a derivative, so, by induction, 
\begin{align}
\ad_\pn^j(F_1(\qm))
=&(i\tl^{m+n-1})^j F_1^{[j]}(\tl^m\rho_*)\\
=&(i\tl^{m+n-1})^j F_1^{[j]}(\qm)
\end{align}
Applying this, and multiplying by $\tl^{k(1-n-m)}$, gives that on each spherical harmonic, 
\begin{align}
\label{eCommutatorShellbyShell}
\tl^{k(m+n-1)}[F_1(\qm),F_2(\pn)] 
=& \sum_{j=1}^{k-1} \tl^{k(m+n-1)} \frac{1}{j!} F_2^{[k]}(\pn) (i\tl^{1-n-m})^jF_1^{[j]}(\qm) + R_k \\
\|R_k\| \leq& C_k \| \Fourier{F_2^{[k]}} \|_1 \|F_1^{[j]} \|_\infty
\end{align}
Since the constant in the commutator expansion theorem, theorem \ref{tIIp2.1}, is independent of the choice of $x$, $A$, $F_1$, and $F_2$, the constant $C_k$ is independent of $\tl$. The result on each spherical harmonic can be extended across all harmonics with out any change. In particular, the operator $R_k$ is uniformly bounded across all spherical harmonics. 
\end{proof}

The commutator order reduction lemma, lemma \ref{lIIp2.2}, is used to estimate commutators of localisation in the phase space variables. Usually it is sufficient to know that the commutator is a bounded operator times negative powers of $L$, and that is what the following lemmas assert for certain classes of localising operators. 

The following lemma applies in the common situation when the localisation in $\pn$ is a Schwartz class function. 

\begin{lemma}
\label{lSchwartzClassInDRCommutatorsAreGood}
There is a constant $C$, such that if $0\leq n \frac12$, $0\leq m\leq\frac12$, $F_1\in C^1(\Reals)$, $F_1'\in L^\infty$, and $F_2$ is of Schwartz class, then there is a bounded operator, $\bndd$, satisfying the following
\begin{align}
[F_1(\qm),F_2(\pn)]=& L^{n+m-1}\bndd_1 \\
\|\bndd_1\|_\opnorm \leq& C \|\Fourier{F_2'}\|_1 \|F_1'\|_\infty .
\end{align}
\end{lemma}
\begin{proof}
Since the Schwartz class is preserved by differentiation and the Fourier transform, $\Fourier{F_2'}$ is Schwartz class and, hence, in $L^1$. The expansion to order $0$ from the commutator order reduction lemma, lemma \ref{lCommutatorOrderReduction}, is
\begin{align}
L^{1-m-n} [F_1(\qm),F_2(\pn)] =& R_1 
\label{AnIntermediaryInlSchwartzClassInDRCommutatorsAreGood}\\
\| R_1 \|_{\opnorm} \leq& C \|\Fourier{F_2'}\|_{L^1} \|F_2'\|_{L^\infty}
\end{align}
Since $L$ is bounded below by $1$, and $0\leq n, m\leq\frac12$, the operator $L^{m+n-1}$ is bounded. The operators $L^{m+n-1}$ and $R_1$ are bounded, so they can be composed. Applying $L^{m+n-1}$ to \eqref{AnIntermediaryInlSchwartzClassInDRCommutatorsAreGood} proves the desired result. 
\end{proof}

The previous lemma can not be applied with the smooth, near localisations for $\pn$, which are not in Schwartz class. To prove a similar result for the smooth, near localisations, it is sufficient to show that each satisfies the condition $\Fourier{F_2'}\in L^1$. The following lemma also proves a commutator order reduction result for two unbounded functions of $\pn$, $F_2=\xi$ and $F_2=\xi \pnna{\xi}$, using a first order expansion in lemma \ref{lCommutatorOrderReduction}. 

\begin{lemma}
\label{lSomeOtherFunctionsInDRCommutatorsAreGood}
There is a constant $C$, such that if $0\leq n\leq\frac12$, $0\leq m\leq\frac12$, $F_1$ has its first in $L^\infty$, $\epsilon>0$, and $F_2(x)$ is one of the following functions: $x, \pnna{x}, \pnnb{x}, \pnnc{x}$, then there is a bounded operator, $\bndd_{F_2}$, satisfying the following
\begin{align}
[F_1(\qm),F_2(\pn)]=& L^{n+m-1}\bndd_{F_2} \\
\|\bndd_{F_2}\|_\opnorm \leq& C \|F_1'\|_\infty .
\end{align}

If the previous conditions apply, $F_1$ also has its second derivative in $L^\infty$, and $F_2(x) = x\pnna{x}$, then there is a bounded operator, $\bndd_{F_2}$, satisfying
\begin{align}
[F_1(\qm),F_2(\pn)]=& L^{n+m-1}\bndd_{F_2} \\
\|\bndd_{F_2}\|_\opnorm \leq& C (\|F_1'\|_\infty +\|F_1''\|_\infty).
\end{align}
\end{lemma}
\begin{proof}
The fundamental problem here is that the permitted $F_2$ functions are not in Schwartz class. There are finitely many functions $F_2$ considered in this theorem, so if a $C$ can be found for each of them, then the largest one can be applied to make the result hold uniformly. 

For the function $F_2(x)=x$, $F_2(\pn)=\pn=-i\dr L^{n-1}$. All the operators involved commute with the spherical harmonic decomposition, and on each spherical harmonic, the commutator can be explicitly computed 
\begin{align}
[F_1(\qm),-i\dr L^{n-1}]
=& [F_1(L^m\rho_*),-i\dr L^{n-1}]\\
=&iL^{n+m-1} F_1'(L^m\rho_*) = iL^{n+m-1} F_1'(\qm)
\end{align}
The commutator is a product of a power of $L$ and a multiplication operator. The $L^\infty$ norm of the function is the operator norm of the associated multiplication operator. 

To apply the commutator order reduction lemma, lemma \ref{lCommutatorOrderReduction} to the remaining functions, we prove: If $h\in C^{2}(\Reals)\cap L^2(\Reals)$ and $h''\in L^2(\Reals)$, then $\|\Fourier{h}\|_1 \leq\infty$. 

Since $h\in C^2$, $h''$ is well defined, and since $h, h'' \in L^2(\Reals)$, both $\Fourier{h}$ and $\Fourier{h''}$ are well defined and in $L^2(\Reals)$. Let 
\begin{align}
f(\xi)\defin& |(1+\xi^2)\Fourier{h}(\xi)| ,\\
f_{\text{big}}(\xi) \defin& f(\xi) \chi([1,\infty),f(\xi)) ,\\
f_{\text{small}}(\xi) \defin& f(\xi) \chi([0,1),f(\xi)) ,
\end{align} 
 and observe that 
\begin{align}
f(\xi)=&f_{\text{big}}(\xi) + f_{\text{small}}(\xi) \\
|\Fourier{g}(\xi)|
\leq & |\frac{f_{\text{big}}(\xi)}{1+\xi^2}| + |\frac{f_{\text{small}}(\xi)}{1+\xi^2}|
\leq  |f_{\text{big}}(\xi)| + |\frac{f_{\text{small}}(\xi)}{1+\xi^2}| .
\end{align}
From properties of the Fourier transform it is known that 
\begin{equation}
|\Fourier{h}(\xi)|+|\Fourier{h''}(\xi)|= |\Fourier{h}(\xi)|+|\xi^2\Fourier{h}(\xi)| =|f(\xi)| ,
\end{equation}
and hence that $f$ is in $L^2(\Reals)$. Since $f_{\text{big}}\leq f_{\text{big}}^2\leq f^2$, $f_{\text{big}}$ is in $L^1(\Reals)$, and since $f_{\text{small}} \leq 1$, $\frac{f_{\text{small}}(\xi)}{1+\xi^2}$ is also in $L^1$. The sum of these bounds $|\Fourier{h}|$, so that $\Fourier{h}$ is in $L^1(\Reals)$.

For the three functions, $\pnna{x}$, $\pnnb{x}$, and $\pnnc{x}$, each of these is a power of a polynomial and decays like $x^{\frac{-1-\epsilon}{2}}$. Therefore, $F_2'$ and $(F_2')''$ are each in $L^2$, since they decay like $x^{\frac{-3+\epsilon}{2}}$ and $x^{\frac{-7+\epsilon}{2}}$ respectively, and $\Fourier{F_2'}$ is in $L^1$. By the same argument as in lemma \ref{lSchwartzClassInDRCommutatorsAreGood}, 
\begin{align}
[F_1(\qm),F_2(\pn)]=& L^{n+m-1}\bndd \\
\|\bndd\|_\opnorm \leq& C (\|F_1'\|_\infty).
\end{align}

The last, and most difficult piece, is for $F_2(x)= x\pnna{x}$. In the previous arguments, the commutator was directly estimated. In this case it is expanded to first order, and both the first order term and the remainder are estimated. The derivative of $x\pnna{x}$ decays too slowly to be in $L^2$ and may not have an $L^1$ Fourier transform. 

As usual, all the operators involved commute with the spherical harmonic decomposition. On each spherical harmonic, by lemma \ref{lIIp2.2}
\begin{align}
[F_1(\qm),F_2(\pn)]=& F_2'(\pn) (iL^{n+m-1})(F_1'(\qm)) + L^{2n+2m-2} R_2 .
\end{align}
The first term is a product of operators which commute with powers of $L$ and each can be estimated in norm. 
\begin{align}
\|F_2'(\pn) (-iL^{n+m-1})(-F_1'(\qm))\|_\opnorm = &L^{n+m-1} \|F_2'(\pn)\|_\opnorm \|F_1'(\qm)\|_\opnorm\\
\leq& L^{n+m-1} \|F_2'\|_\infty \|F_1'\|_\infty 
\end{align}
Since $F_2(x)=x\pnna{x}$, $F_2'$ is bounded. 

The second term is estimated by lemma \ref{lIIp2.2}.
\begin{align}
\| R_2\|_\opnorm = C \|\Fourier{F_2''}\|_1 \|F_1''\|_\infty 
\end{align}
Since $F_2=x\pnna{x}$, $F_2''$ is a smooth function, which is in $L^2$ and has its second derivative, $(F_2'')''$, also in $L^2$, we conclude that $\Fourier{F_2''}$ is in $L^1$. 

Since $L^{2n+2m-2}\leq L^{n+m-1}$, the estimate on each of the terms in the expression for the commutator can be combined to give
\begin{align}
\|[F_1(\qm),F_2(\pn)]\|_\opnorm \leq& L^{n+m-1} \|F_2'\|_\infty \|F_1'\|_\infty + L^{2n+2m-2} C \|\Fourier{F_2''}\|_1 \|F_1''\|_\infty \\
\leq& C L^{n+m-1}(\|F_1'\|_\infty+\|F_1''\|_\infty)
\end{align}

On each spherical harmonic there is a bounded operator, $\bndd$, such that 
\begin{align}
[F_1(\qm),F_2(\pn)] =& L^{n+m-1}\bndd\\
\|\bndd\|_\opnorm \leq&  C (\|F_1'\|_\infty+\|F_1''\|_\infty) .
\end{align}
Therefore, $\bndd$ can be extended as an operator on $L^2$ satisfying the same condition. 
\end{proof}

%
%
\subsection{Phase Space Estimates}
\label{ssPhaseSpaceLocalisedEstimates}

In this subsection, we compute the commutators between $\Gammanm$ and the components of the Hamiltonian, $H_i$. This is analogous to subsection \ref{ssMorawetzCommutatorCalculation} with phase space localisation and additional factors of $L$. The $H_2$ commutator is simply lower order. The $H_1+H_3$ commutator dominates the radial derivative localised to small $\qm$ and powers of $L$ localised to large $\qm$.

%
%

The commutator with $H_2$ is computed first. Because $H_2$ is a function of $\rho_*$ and contains no derivatives, all the commutators can be shown to be lower order by lemma \ref{lSomeOtherFunctionsInDRCommutatorsAreGood}. 

\begin{lemma}
\label{lfPSH2}
For \uasoluinS, and $0\leq n\leq m\leq\frac12$, 
\begin{equation}
\langle u,[H_2,\Gammanm]u\rangle = \Oint
\end{equation}
\end{lemma}
\begin{proof}
Lemma \ref{lIIp3.1}, the Leibniz formula for commutators, is used to calculate $[H_2,\Gammanm]$ and then the commutators are shown to be lower order by lemma \ref{lSomeOtherFunctionsInDRCommutatorsAreGood}. Since $\chia$ and $V$ commute, $[H_2,\chia]=0$ and only three terms are left in the expansion of the $[H_2,\Gammanm]$.
\begin{align}
\langle u,[H_2,\Gammanm]u\rangle=\langle u,\chia\ L^{n-\epsilon}(&[V,\pnna{\pn}]\gamma_{L^m}\pnna{\pn}\\
&+\pnna{\pn}[V,\gamma_{L^m}]\pnna{\pn}\\
&+\pnna{\pn}\gamma_{L^m}[V,\pnna{\pn}])\chia u\rangle
\end{align}
The commutator $[V,\gamma_{L^m}]$ only involves the angularly modulated multiplier, and was shown to be a bounded operator in lemma \ref{lII4.1}. Applying additional bounded operators and powers of $L$ yields a bounded operator multiplied by the same power of $L$. 
\begin{equation}
L^{n-\epsilon}\pnna{\pn}[V,\gamma_{L^m}]\pnna{\pn}=L^{n-\epsilon}\bndd
\end{equation}
Since $n\leq\frac32$, by the angular modulation theorem, theorem \ref{tII5.1}, the expectation value of this operator is $\Oint$.
\begin{align}
\langle u, \chia\ L^{n-\epsilon}\pnna{\pn}[V,\gamma_{L^m}]\pnna{\pn} \chia u\rangle
=&\langle u,L^{n-\epsilon}\bndd u\rangle\\
=&\Oint .
\end{align}

The remaining terms are complex conjugates of each other. Only one is considered here, and the other can be dealt with in the same way. The commutator is shown to be a bounded operator times a negative power of $L$ using lemma \ref{lSomeOtherFunctionsInDRCommutatorsAreGood}. $\gamma_{L^m}$ can be expanded and the bounded functions involved in this expansion can be absorbed into the bounded operator from the commutator expansion. 
\begin{align} 
L^{n-\epsilon}\pnna{\pn} \gamma_{L^m}[V,\pnna{\pn} ]=& L^{n-\epsilon}\pnna{\pn} (\dr g_{L^m}-\frac12g_{L^m}')L^{n-1}B\\
=& L^{2n-1-\epsilon}\pnna{\pn} \dr\bndd_1 +L^{2n-1-\epsilon} L^m \frac12 \qmn{\qm} \bndd_2\\
=& L^{2n-1-\epsilon}\pnna{\pn} \dr\bndd_1+L^{m+2n-1-\epsilon}\bndd .
\end{align}
This operator can be substituted into the relevant expectation value. This gives a sum of two terms, both of which can be estimated by the Cauchy-Schwartz inequality. 
\begin{align}
\langle u,\chia L^{n-\epsilon}\pnna{\pn}& \gamma_{L^m}[V,\pnna{\pn} ]\chia u\rangle\\
=&\langle -\dr\pnna{\pn} \chia u,L^{2n-1-\epsilon}\bndd_1\chia u\rangle+\langle \chia u,L^{m+2n-1-\epsilon}\bndd \chia u\rangle\\
\geq& -C(\|\dr\chia u\| \|L^{\frac{2n-1-\epsilon}{2}}\bndd_1\chia u\| +\|L^{\frac{3-\epsilon}{2}n}\chia u\|^2)\\
\geq& -C(\|\dr\chia u\|^2+ \|L^{\frac{2n-1-\epsilon}{2}}\chia u\|^2 +\|L^{\frac{3-\epsilon}{2}n}\chia u\|^2)
\end{align}
Since $2n-1-\epsilon\leq\frac32$ and $(3-\epsilon)n\leq\frac32$, by the angular modulation theorem, theorem \ref{tII5.1}, the two terms involving powers of $L$ are time integrable, ie $\Oint$. By corollary \ref{cII6.1}, $\|\dr\chia u\|^2=\Oint$. Therefore the expectation value being considered is time integrable
\begin{align}
\langle u,\chia L^{n-\epsilon}\pnna{\pn} \gamma_{L^m}[V,\pnna{\pn}]\chia u\rangle =&\Oint
\end{align}
\end{proof}

%
%

Because $H_3=(-\Delta_{S^2})V_L=(L^2-1) V_L$ contains two angular derivatives, which are exactly what is under consideration, this commutator requires the most care. The difference between $L^2-1$ and $L^2$ is a local decay term. The commutator is first expanded using the product rule in lemma \ref{lIIp3.1}. The commutator with $\gamma_{L^m}$ is computed from the angular modulation argument in lemma \ref{lAngModResultLikePreviouslII4.2}. The commutator with the localisation $\pnna{\pn}$ can be expanded to second order using the commutator order reduction lemma, lemma \ref{lIIp2.2}. The first order term is of the same order as the commutator with $\gamma_{L^m}$. It involves a derivative of $\pnna{\pn}$ and is multiplied by the $\dr$ term from $\gamma_{L^m}$. There are two such terms, and when rearranged and combined with the term from the commutator with $\gamma_{L^m}$, they give $\pnnc{\pn}^2=\pnna{\pn}(\pnna{\pn}+2\pnnaWITHADERIV{\pn})$. This term must be positive to give the desired estimate. This dictates the choice of the decay function $\pnna{x}$ and forces $\epsilon>0$. The second order term in the expansion naively appears to be of a high order since $H_3$ introduced two factors of $L$. However, rearrangements yield a cancellation of the highest order part of this term and the remainder is a lower order error term. For both the first and second commutators, there are many rearrangements of phase space localising functions. These are shown to be lower order in the standard way using lemma \ref{lSomeOtherFunctionsInDRCommutatorsAreGood}. The remainder term is shown to be lower order directly from the commutator order reduction lemma, lemma \ref{lIIp2.2}. 

\begin{proposition}
\label{lfPSH1}
For \uasoluinS{} and $0\leq n\leq m\leq\frac12$, 
\begin{align}
\langle u,[H_1+H_3,\Gammanm]u\rangle
\geq& 
c \langle \pnna{\pn}\chia u, (-\dr \qmnt{\qm}\dr - \epsilon L^2  g_{L^m} V_L') \pnna{\pn}\chia u\rangle
+ \ErrorTermsn
\end{align}
\end{proposition}
\begin{proof}
The commutator is given by the Leibniz rule as the sum of five terms. 
\begin{align}
\label{ePSECommutator}
[H_1+H_3,\Gammanm]
=&[H_1+H_3,\chia]\pnna{\pn} L^{n-\epsilon} \gamma_{L^m} \pnna{\pn}\chia +
\chia\pnna{\pn} L^{n-\epsilon}\gamma_{L^m}\pnna{\pn}[H_1+H_3,\chia] \\
& + \chia[H_1+H_3,\pnna{\pn}] L^{n-\epsilon}\gamma_{L^m}\pnna{\pn}\chia +
\chia\pnna{\pn} L^{n-\epsilon}\gamma_{L^m}[H_1+H_3,\pnna{\pn}]\chia \\
& +\chia\pnna{\pn} L^{n-\epsilon}[H_1+H_3,\gamma_{L^m}]\pnna{\pn}\chia 
\end{align}
Since the first two terms are adjoints of each other, it is sufficient to estimate only one of them. The remaining three terms are more complicated and must be estimated together.

Since $\chia$ is a function, it commutes with $V_L$, and we only need to determine the contribution from $-\drr$, which can be explicitly calculated. 
\begin{align}
&\langle u,\chia\pnna{\pn}\gamma_{L^m}L^{n-\epsilon}\pnna{\pn}[-\drr,\chia]u\rangle\\
&= \langle u,\chia\pnna{\pn}\gamma_{L^m}L^{n-\epsilon}\pnna{\pn}(-2\dr\chia'-\chia'')u\rangle
\end{align}
To analyse this term, $\gamma_{L^m}$ is expanded into two terms so that the most significant term can be dealt with first. The slightly less common expansion $\gamma_{L^m}=\dr g_{L^m}-\frac12 g_{L^m}'$ is used. There are a total of four terms to consider, two from expanding $[H_1,\chia]$ times two from expanding $\gamma_{L^m}$. 

The term that appear to be highest order is the one involving $\dr g_{L^m}$ from $\gamma_{L^m}$ and $-2\dr\chia'$ from $[-\drr,\chia]$. It will be referred to as $I_{0}$. It can be simplified by moving one of the derivative operators through $\pnna{\pn}$ (since they commute), and then moving it and a factor of $\chia$ to the other side of the inner product. 
\begin{align}
I_{0}\defin &\langle u,\chia\pnna{\pn} \dr g_{L^m}L^{n-\epsilon}\pnna{\pn}(-2\dr\chia')u\rangle\\
=&-\langle \dr\chia u,\pnna{\pn} g_{L^m}L^{n-\epsilon}\pnna{\pn}(-2\dr\chia')u\rangle
\end{align}
Applying the Cauchy-Schwartz estimate gives the sum of two terms. One of these is $\|\dr \chi u\|^2$ which is $\Oint$ (integrable in time) by corollary \ref{cII6.1}. The other can be rewritten in terms of $\pn$. 
\begin{align}
I_{0}\geq& -\|\dr\chia u\|^2-\|\pnna{\pn} g_{L^m}L^{n-\epsilon}\pnna{\pn}\dr\chia' u\|^2\\
\geq& -\|L^{1-\epsilon} \pnna{\pn} g_{L^m}\pn\pnna{\pn}\chia' u\|^2 +\Oint
\end{align}
The remaining norm is going to be manipulated with the goal of applying lemma \ref{lYnwelldef} to control $\pn\pnna{\pn}^2$ by powers of $L$. This manipulation will introduce additional error terms from commuting the phase space localisations. To commute $\pnna{\pn}$ and $g_{L^m}$ while avoiding any error terms involving factors of $\pn$ outside a commutator, it is necessary to commute $\pn$ through $g_{L^m}$ and then commute $\pn\pnna{\pn}$ back through $g_{L^m}$. 
\begin{align}
I_{0}
\geq& -\|L^{1-\epsilon} \pn\pnna{\pn} g_{L^m}\pnna{\pn}\chia'u\|^2\\
&-\|L^{1-\epsilon} \pnna{\pn} [g_{L^m},\pn]\pnna{\pn}\chia' u\|^2 +\Oint\\
\geq& -\|L^{1-\epsilon}  g_{L^m}\pn\pnna{\pn}^2\chia'u\|^2 -
\|L^{1-\epsilon} [\pn\pnna{\pn}, g_{L^m}]\pnna{\pn}\chia'u\|^2 \\
&-\|L^{1-\epsilon} \pnna{\pn} [g_{L^m},\pn] \pnna{\pn}\chia' u\|^2+\Oint
\end{align}
By lemma \ref{lSomeOtherFunctionsInDRCommutatorsAreGood}, the commutators are of the form $L^{m+n-1}\bndd$. Since $\pnna{\pn}$ is also a bounded operator, all the error terms from commuting are lower order terms. In the remaining term, the factor of $g_{L^m}$ can be dropped for a constant and then lemma \ref{lYnwelldef} can be applied to control $\pn\pnna{\pn}^2$. 
\begin{align}
I_{0} \geq& -C\|L^{1-\epsilon} g_{L^m}\pn\pnna{\pn}^2\chia'u\|^2 
-C\|L^{1-\epsilon} L^{n+m-1}\bndd\chia'u\|^2+\Oint\\
\geq& -C(\|L^\frac{1+n\epsilon-2\epsilon}{1-\epsilon}\chia'u\|^2+\|\dr\chia'u\|^2)-\ErrorTermsn
\end{align}
The angular energy away from $r=\alpha$ is time integrable by corollary \ref{cReductionToLonPhotonSphere}. Since $\chia$ is compactly supported and identically $1$ in a neighbourhood of $r=\alpha$, $|\chia'|^2$ is dominated by a multiple of $(\rho_* g(\rho_*))$ and $\|L \chia'u\|$ is integrable in time. Interpolating between the local decay result and the angular energy decay away from the photon sphere result, $\|L^\frac{1+n\epsilon-2\epsilon}{1-\epsilon}\chia'u\|^2$ is time integrable, and hence in $\Oint$. By corollary \ref{cII6.1}, $\|\dr\chia'u\|^2$ is time integrable and hence in $\Oint$. All the terms in $I_{0}$ are now controlled. 
\begin{align}
I_{0}\defin &\langle u,\chia\pnna{\pn} \dr g_{L^m}L^{n-\epsilon}\pnna{\pn}(-2\dr\chia')u\rangle \geq \ErrorTermsn .
\end{align}

The term that appears to be the next highest order is also $\ErrorTermsn$ not merely $\Oint$. This is the term involving $g_{L^m}'$ from $\gamma_{L^m}$ and $-2\dr\chia'$ from $[-\drr,\chia]$. It is estimated by moving most of the localisation to the other side of the inner product, applying the Cauchy-Schwartz inequality, replacing bounded operators by a constant, and then finding that the remaining terms are each $\ErrorTermsn$. 
\begin{align}
\langle u,\chia \pnna{\pn} (-\frac12& \frac{L^{m+n-\epsilon}}{(1+L^m|\rho_*|)^2})\pnna{\pn}\dr \chia'u\rangle\\
=& -\frac12\langle \pnna{\pn} L^{m+n-\epsilon} \qmn{\qm} \pnna{\pn}\chia u,\dr \chia'u\rangle\\
\geq&-C \|L^{m+n-\epsilon}\chia u\| \|\dr \chia' u\|\\
\geq&-C(\|L^{m+n-\epsilon}\chia u\|^2+\|\dr \chia' u\|^2)\\
\geq&-\ErrorTermsn
\end{align}

The terms involving $\chia''$ from $[-\drr,\chia]$ are the lowest order terms and are quickly shown to be lower order. The same method as was used for the term involving $g_{L^m}'$ and $-2\dr\chia'$ can be applied. 
\begin{align}
\langle u,\chia\pnna{\pn}\dr L^{n-\epsilon}g_{L^m}\pnna{\pn} \chia'' u\rangle
=&\langle g_{L^m}\dr \pnna{\pn}\chia u,L^{n-\epsilon}\pnna{\pn}\chia''u\rangle \\
 \geq&-C\|\dr\chia u\|\|L^{n-\epsilon}\chia''u\|\\
\geq&-\Oint\\
=&\ErrorTermsn
\end{align}
\begin{align}
\langle u,\chia\pnna{\pn}&\frac{L^{m+n-\epsilon}}{1+(\frac{\rho_*L^m}{2M})^2}\pnna{\pn}\chia''u\rangle\\
=& \langle\pnna{\pn}\chia u,L^{m+n-\epsilon} \qmn{\qm} \pnna{\pn}\chia''u\rangle\\
\geq&-C\|L^{\frac{m+n-\epsilon}{2}}\pnna{\pn}\chia u\|\|L^{\frac{m+n-\epsilon}{2}}\pnna{\pn}\chia u\|\\
\geq&-C\|L^{\frac{m+n-\epsilon}{2}}\chia u\|^2\\
\geq&-\ErrorTermsn
\end{align}


We now turn to the three remaining terms in \eqref{ePSECommutator}. For the commutators with $\pnna{\pn}$, since $-\drr$ commutes with $\pn$, it is sufficient to compute the commutator with $V_L$. No distinction is made between $[H_3,\Gammanm]=[V_L(L^2-1),\Gammanm]$ and $[V_L L^2,\Gammanm]$ since the difference is lower order by an identical argument to that in lemma \ref{lfPSH2}. These commutators are expanded to third order. 
\begin{align}
[V_L,\pnna{\pn}] = \pnnaWITHADERIV{\pn}[V_L,\pn] + \pnnaWITHTWODERIV{\pn}[[V_L,\pn],\pn] + R_3 .
\end{align}
In each of the third and fourth terms of \eqref{ePSECommutator}, we make this substitution and expand $\gamma_{L^m}$ as a sum of two terms. This gives twelve terms, which are combined with the fifth term in \eqref{ePSECommutator}, for a total of thirteen. These are grouped as
\begin{align}
\chia[H_3,\pnna{\pn}] L^{n-\epsilon}\gamma_{L^m}\pnna{\pn}\chia \\
+\chia\pnna{\pn} L^{n-\epsilon}\gamma_{L^m}[H_3,\pnna{\pn}]\chia \\
+\chia\pnna{\pn} L^{n-\epsilon}[H_1+H_3,\gamma_{L^m}]\pnna{\pn}\chia 
=&\chia\pnna{\pn} L^{n-\epsilon}[-\drr,\gamma_{L^m}]\pnna{\pn}\chia \\
&+ \chia( I_1 + I_2 + I_3 + I_4 + I_5 + I_6 ) \chia.
\label{eI1toI6introduced}
\end{align}
where 
\begin{align}
I_1=L^{2+n-\epsilon} (&\pnnaWITHADERIV{\pn} [V_L,\pn]\dr g_{L^m} \pnna{\pn}\\
& + \pnna{\pn} g_{L^m} [V_L,\dr] \pnna{\pn}\\
& + \pnna{\pn} g_{L^m}\dr \pnnaWITHADERIV{\pn}[V_L,\pn])\\
I_2=L^{2+n-\epsilon} \frac12 (&-\pnna{\pn}'[V_L,\pn]g_{L^m}'\pnna{\pn} +\pnna{\pn} g_{L^m}'\pnna{\pn}'[V_L,\pn])\\
I_3=L^{2+n-\epsilon}(&\pnnaWITHTWODERIV{\pn}[[V_L,\pn],\pn]\dr g_{L^m}\pnna{\pn}\\
&+\pnna{\pn} g_{L^m}\dr\pnnaWITHTWODERIV{\pn}[[V_L,\pn],\pn])\\
I_4=L^{2+n-\epsilon}\frac12(&-\pnnaWITHTWODERIV{\pn}[[V_L,\pn],\pn]g_{L^m}'\pnna{\pn}\\
& +\pnna{\pn} g_{L^m}'\pnnaWITHTWODERIV{\pn}[[V_L,\pn],\pn])\\
I_5=L^{2+n-\epsilon}(&R_3 \dr g_{L^m}\pnna{\pn}+\pnna{\pn} g_{L^m}\dr R_3)\\
I_6=L^{2+n-\epsilon}\frac12(&-R_3 g_{L^m}'\pnna{\pn} +\pnna{\pn} g_{L^m}\dr R_3)
\end{align} 

The leading order term\footnote{From here to the end of the estimate on $I_6$, the proof is cut and pasted from the proof of prop 8.15 in the previous notes. }, $I_1$, is composed of two types of terms. One type comes from the commutator $[L^2 V_L,\gamma_{L^m}]$, and the other comes from combining  the first order term in the expansion of $[H_3,\pnna{\pn}]$ with the derivative terms in $\gamma_{L^m}$. To get the desired estimate, it is necessary to rearrange the phase space localising functions, which introduces additional lower order error terms. In the calculations, leading order terms are written first here, followed by those which will be shown to be lower order. 

The term $I_1$ is rearranged so that all the leading order subterms contain the $\rho_*$ localisation factor, $ (-V_L'g_{L^m})$, which first appears in the commutator $[V_L,\gamma_{L^m}]$, and is the localisation in the statement of the lemma. The commutators from rearranging localisations are lower order by lemma \ref{lSomeOtherFunctionsInDRCommutatorsAreGood}. 
\begin{align}
\label{ExpansionOfH3CommutatorOnceWasLine303}
I_1=
L^{2+n-\epsilon}(&\pnna{\pn}'(-V_L')(-iL^{n-1})\dr g_{L^m}\pnna{\pn}\\
&+\pnna{\pn} g_{L^m}[V_L,\dr]\pnna{\pn}\\
&+\pnna{\pn} g_{L^m}\dr\pnna{\pn}'(-V_L')(-iL^{n-1}))\\
= L^{2+n-\epsilon}(&\pnna{\pn}'\pn (-V_L'g_{L^m})\pnna{\pn}\\
&+\pnna{\pn}(-V_L'g_{L^m})\pnna{\pn} \\
&+\pnna{\pn} (-V_L'g_{L^m})\pn\pnna{\pn}')\\
+L^{2+n-\epsilon}(&\pnna{\pn}'[-V_L',\pn]g_{L^m}\pnna{\pn}+\\
&\pnna{\pn} g_{L^m}[\pn\pnna{\pn}',-V_L'])\\
=L^{2+n-\epsilon}(&\pn\pnna{\pn}'(-V_L' g_{L^m})\pnna{\pn}\\
&+\pnna{\pn}(-V_L' g_{L^m})\pnna{\pn}\\
&+\pnna{\pn}(-V_L' g_{L^m})\pn\pnna{\pn}')\label{leadingH3commutatorfirstnoon23Dec}\\
+L^{1+2n-\epsilon}&\bndd
\end{align}
Following this, the $\pn$ localising functions are gathered to the right of the $\rho_*$ space localisation. Once they are gathered this way, their sum is the operator $\pnnc{\pn}^2$. One factor of $\pnnc{\pn}$ can be moved back through the $\rho_*$ localisation, so that $\pnnc{\pn}$ appears on both sides of the $\rho_*$ localisation. Rearranging in this way introduces new commutators, which are not immediately estimated since they involve $(-V_L'g_{L^m})$. This is a localising function in $\rho_*$ and $\rho_*L^m$ so that the previous lemmas which estimate commutators can not be applied. 
\begin{align}
I_1
=L^{2+n-\epsilon}(&-V_L' g_{L^m})\pnna{\pn}(\pnna{\pn}+2\pn\pnna{\pn}')\\
+L^{2+n-\epsilon}&([\pn\pnna{\pn}',-V_L' g_{L^m}]\pnna{\pn}\\
&+[\pnna{\pn},-V_L' g_{L^m}]\pnna{\pn}\\
&+[\pnna{\pn},-V_L' g_{L^m}]\pn\pnna{\pn}')\\
+L^{1+2n-\epsilon}&\bndd\\
=L^{2+n-\epsilon}(& \pnnc{\pn}(-V_L'g_{L^m})\pnnc{\pn})\\
L^{2+n-\epsilon}(&[(-V_L'g_{L^m}),\pnnc{\pn}]\pnnc{\pn}\\
&+[\pn\pnna{\pn}',-V_L' g_{L^m}]\pnna{\pn}\\
&+[\pnna{\pn},-V_L' g_{L^m}]\pnna{\pn}\\
&+[\pnna{\pn},-V_L' g_{L^m}]\pn\pnna{\pn}')\\
+L^{1+2n-\epsilon}&\bndd
\end{align}
Of the terms in the previous expression, the first is a positive term we wish to keep, and the expectation value of the last is $\ErrorTermsn$. The other terms all involve a commutator of the form $[F(\pn),(-V_L'g_{L^m})]$ multiplied by $L^{2+n-\epsilon}$ and bounded operators. 

By rewriting the commutator $[F(\pn),(-V_L'g_{L^m})]$, the four remaining commutator terms can also be shown to be $\ErrorTermsn$. It is necessary to separate $-V_L'g_{L^m}$ into a function of $\rho_*$ and of $\rho_*L^m$ to estimate the commutator. Naively applying the Leibniz rule for commutators, gives a contribution from the commutator with $g_{L^m}$ which is $L^{1+m+2n-\epsilon}\bndd$, which is not $\ErrorTermsn$ in expectation value. The decomposition 
\begin{align}
-V_L' g_{L^m} = (-V_L' \rho_*^{-1}) (\rho_*L^m g_{L^m}) L^{-m}
\end{align}
ensures that the commutators are all $\ErrorTermsn$ in expectation. This decomposition is possible since $-V_L'$ vanishes linearly at $\rho_*=0$. Taylor's theorem guarantees $-V_L'\rho_*^{-1}$ is $C^\infty$. 
\begin{align}
L^{2+n-\epsilon}[F(\pn),(-V_L'g_{L^m})]
=& L^{2-m+n-\epsilon}[F(\pn),(-V_L'\rho_*^{-1}) L^m\rho_*]
\label{eVlgmcommutator}
\end{align}
The Leibniz rule can now be applied and the resulting commutators shown to be lower order operators by lemma \ref{lSomeOtherFunctionsInDRCommutatorsAreGood}. Note that the class of $F$ considered here consists of $\pnnc{x}$, $\pnnb{x}$, and $\pnna{x}$, all of which are covered by lemma \ref{lSomeOtherFunctionsInDRCommutatorsAreGood}, and that $(V_L' \rho_*^{-1})$ and $L^m\rho_* g_{L^m}=\qmf{\qm}^2$ both have bounded first and second derivative as required by lemma \ref{lSomeOtherFunctionsInDRCommutatorsAreGood}. 
\begin{align}
L^{2+n-\epsilon}[F(\pn),(-V_L'g_{L^m})]
=& L^{2-m+n-\epsilon}(-V_L' \rho_*^{-1}) [F(\pn),\qmf{\qm}^2]\\
&+ L^{2-m+n-\epsilon}[F(\pn),(-V_L' \rho_*^{-1})]\qmf{\qm}^2\\
=& L^{2+n-\epsilon}L^{n+m-1}\bndd + L^{2-m+n-\epsilon}L^{m+n-1}\rho_*\bndd \\
=& L^{1+2n-\epsilon}\rho_* \bndd +L^{1+2n-\epsilon}\bndd
\end{align}
Since $\chia$ has compact support, and applying the self-adjointness of $\rho_*$, we are left with an error term like 
\begin{align}
|\langle\rho_*\chia u, L^\frac{1+2n-\epsilon}{2}\bndd\chia u\rangle |
\leq& C\|L^\frac{1+2n-\epsilon}{2}\chia u\|\\
=&\ErrorTermsn
\end{align}
This gives the following estimate of the leading order term, $I_1$. 
\begin{equation}
\langle u,\chia I_1 \chia u\rangle= \| L^\frac{2+n-\epsilon}{2}(-V_L'\rho_*^{-1})^\frac12 \qmf{\qm} \pnnc{\pn} \chia u\|^2 +\ErrorTermsn
\label{eVlgmcommutatora}
\end{equation}
Using the estimate in equations \eqref{eVlgmcommutator}-\eqref{eVlgmcommutatora} again, we can rearrange the localisation as 
\begin{equation}
\langle u,\chia I_1 \chia u\rangle= \| L^\frac{2+n-\epsilon}{2}\pnnc{\pn}(-V_L'\rho_*^{-1})^\frac12 \qmf{\qm} \chia u\|^2 +\ErrorTermsn
\end{equation}
The localisation $\pnnc{\pn}$ can be replaced by $\epsilon^{1/2} \pnna{\pn}$, so the estimate becomes, 
\begin{equation}
\langle u,\chia I_1 \chia u\rangle= \epsilon \| L^\frac{2+n-\epsilon}{2}\pnna{\pn}(-V_L'\rho_*^{-1})^\frac12 \qmf{\qm} \chia u\|^2 +\ErrorTermsn .
\end{equation}
To conclude the estimate on $I_1$, we commute the localisation again, 
\begin{equation}
\langle u,\chia I_1 \chia u\rangle= \epsilon \| L^\frac{2+n-\epsilon}{2}(-V_L'\rho_*^{-1})^\frac12 \qmf{\qm} \pnna{\pn}\chia u\|^2 +\ErrorTermsn .
\end{equation}

It now remains to show that $I_2, \ldots, I_6$ are $\ErrorTermsn$. Discussion of these continues from highest order to lowest. The first of these is $I_2$ which contains terms involving $g_{L^m}'$ and $\pn\pnnaWITHADERIV{\pn}$ terms. Note that because two different expansions of $\gamma_{L^m}$ were used in equation \eqref{eI1toI6introduced}, the signs on the two factors of $g_{L^m}'$ are different. To simplify $I_2$, the factors of $g_{L^m}'$ are expanded. 
\begin{align}
I_2=L^{2+n-\epsilon}\frac12(&-\pnnaWITHADERIV{\pn}[V_L,\pn]g_{L^m}'\pnna{\pn}\\
&+\pnna{\pn} g_{L^m}'\pnnaWITHADERIV{\pn}[V_L,\pn])\\
=L^{2+n-\epsilon}\frac12(&-\pnnaWITHADERIV{\pn}(-iL^{n-1})(-V_L')L^m \qmn{\qm}\pnna{\pn}\\
&+\pnna{\pn}\qmn{\qm}\pnnaWITHADERIV{\pn}(-iL^{n-1})(-V_L'))\\
=\frac{i}{2}L^{1+m+2n-\epsilon}(&-\pnnaWITHADERIV{\pn}V_L'\qmn{\qm}\pnna{\pn}\\
&+\pnna{\pn} \qmn{\qm}\pnnaWITHADERIV{\pn}V_L')
\end{align}
This leaves an expression for $I_2$ as the sum of two terms. The goal now is to rearrange these terms so that the highest order parts from each cancel. As usual, the rearrangement will introduce commutators. Naively attempting to commute $\pnna{\pn}$ or $\pnnaWITHADERIV{\pn}$ with $g_{L^m}$ will leave a commutator of the form $L^{2m+3n-\epsilon}\bndd$, which is not $\ErrorTermsn$ in expectation value. Once again, an additional factor of $\rho_*\rho_*^{-1}$ is introduced and used to absorb a factor of $L^m$ into a function of $\rho_*L^m$. 
\begin{align}
I_2=\frac{i}{2}L^{1+2n-\epsilon}(&-\pnnaWITHADERIV{\pn} \rho_*L^m \qmn{\qm} (V_L' \rho_*^{-1})\pnna{\pn}\\
&+\pnna{\pn} L^m \qmn{\qm}\pnnaWITHADERIV{\pn}\rho_*(V_L' \rho_*^{-1}))
\end{align}
Functions of $\rho_*$ are now moved to obtain the same localising functions, $(V_L' \rho_*^{-1})$ and $\rho_*L^m \qmn{\qm} = \qm\qmn{\qm}$, in both terms. This rearrangement introduces additional commutators which can be estimated by lemma \ref{lSomeOtherFunctionsInDRCommutatorsAreGood}. 
\begin{align}
I_2=\frac{i}{2}L^{1+2n-\epsilon}(&-\pnnaWITHADERIV{\pn}(\qm\qmn{\qm}) \pnna{\pn}(V_L' \rho_*^{-1})\\
&+\pnna{\pn} (\qm\qmn{\qm}) \pnnaWITHADERIV{\pn}(V_L' \rho_*^{-1}))\\
+\frac{i}{2}L^{1+2n-\epsilon}(&-\pnnaWITHADERIV{\pn}(\qm\qmn{\qm}) [(V_L' \rho_*^{-1}),\pnna{\pn}]\\
&+\pnna{\pn}\qmn{\qm}L^m [\pnna{\pn},\rho_*](V_L' \rho_*^{-1}))\\
=\frac{i}{2}L^{1+2n-\epsilon}(&-\pnnaWITHADERIV{\pn}(\qm\qmn{\qm}) \pnna{\pn}(V_L' \rho_*^{-1})\\
&+\pnna{\pn} (\qm\qmn{\qm}) \pnnaWITHADERIV{\pn}(V_L' \rho_*^{-1}))\\
+L^{3n-\epsilon}&\bndd_1+ L^{m+3n-\epsilon}\bndd_2
\end{align}
The $\pn$ localisation is now grouped to the left of all $\rho_*$ and $\rho_*L^m$ localisation. The leading order terms are both $\pnna{\pn}\pnnaWITHADERIV{\pn}(\qm\qmn{\qm}) (V_L' \rho_*^{-1})$, but with opposite signs, so they cancel. The remaining terms are commutators which are found to be $L^{m+3n-\epsilon}\bndd$ in the standard way. Although $\qmn{x}$ is not $C^2$, we will only need to use $x\qmn{x}$ which has its second derivative in $L^\infty$. 
\begin{align}
I_2
=\frac{i}{2}L^{1+2n-\epsilon}(&-\pnnaWITHADERIV{\pn}[(\qm\qmn{\qm}) ,\pnna{\pn}](V_L' \rho_*^{-1})\\
&+\pnna{\pn} [(\qm\qmn{\qm}) ,\pnnaWITHADERIV{\pn}](V_L' \rho_*^{-1}))\\
+L^{3n-\epsilon}&\bndd_1 +L^{m+3n-\epsilon}\bndd_2\\
=L^{m+3n-\epsilon}&\bndd_3 +L^{3n-\epsilon}\bndd_1
\end{align}
Because the powers of $L$ involved are all less than $1+2n$, the 
desired estimate on $I_2$ holds.
\begin{equation}
\langle \chia u, I_2 \chia u\rangle =\ErrorTermsn
\end{equation}

Next $I_3$ is considered. This involves a second commutator $[[V_L,\pn],\pn]$ which can be explicitly computed as $(V_L'')(-L^{2n-2})$. 
\begin{align}
I_3=L^{2+n-\epsilon}(&\pnnaWITHTWODERIV{\pn}[[V_L,\pn],\pn]\dr g_{L^m}\pnna{\pn}\\
&+\pnna{\pn} g_{L^m}\dr\pnnaWITHTWODERIV{\pn}[[V_L,\pn],\pn])\\
= L^{2+n-\epsilon}(&\pnnaWITHTWODERIV{\pn}(V_L'')(-L^{2n-2})\dr g_{L^m} \pnna{\pn}\\
&+\pnna{\pn} g_{L^m} \dr\pnnaWITHTWODERIV{\pn}(V_L'')(-L^{2n-2}))
\end{align}
To further simplify this, it is rewritten in terms of $\pn$ instead of $\dr$. 
\begin{align}
I_3 = iL^{1+2n-\epsilon}(&\pnnaWITHTWODERIV{\pn}(V_L'')\pn g_{L^m} \pnna{\pn}\\
&+\pnna{\pn} g_{L^m} \pn\pnnaWITHTWODERIV{\pn}(V_L''))
\end{align}
This is rearranged to group $\pn$ and $\pnnaWITHTWODERIV{\pn}$ together. Since $x\pnnaWITHTWODERIV{x}$ is smooth and decays like $x^{\frac{-5+\epsilon}{4}}$ as $x\rightarrow\infty$, $x\pnnaWITHTWODERIV{x}$ is a bounded function and $\pn\pnnaWITHTWODERIV{\pn}$ is a bounded operator. The rearrangement introduces commutators which are estimated in the standard way by lemma \ref{lSomeOtherFunctionsInDRCommutatorsAreGood}. 
\begin{align}
I_3 = L^{1+2n-\epsilon}(&(-i\pn\pnnaWITHTWODERIV{\pn})V_L'' g_{L^m}\pnna{\pn}\\
&+\pnna{\pn} g_{L^m}(-i\pn\pnnaWITHTWODERIV{\pn})V_L'')\\
+L^{1+2n-\epsilon}(&\pnnaWITHTWODERIV{\pn}[V_L'',-i\pn]g_{L^m}\pnna{\pn})\\
= L^{1+2n-\epsilon}&\bndd_1+L^{3n-\epsilon}\bndd_2
\end{align}
Since the powers of $L$ involved are bounded by $1+2n$, the expectation value of $I_3$ with respect to $\chia u$ is $\ErrorTermsn$. 

To estimate $I_4$, the commutator $[[V_L,\pn],\pn]$ and the derivative $g_{L^m}'$ are explicitly expanded to show $I_4$ is a product of bounded operators and powers of $L$. 
\begin{align}
L^{2+n-\epsilon}\frac12(&-\pnnaWITHTWODERIV{\pn}[[V_L,\pn],\pn]g_{L^m}'\pnna{\pn}\\
&+  \pnna{\pn} g_{L^m}'\pnnaWITHTWODERIV{\pn}[[V_L,\pn],\pn])\\
= \frac12 L^{2+n-\epsilon}(&-\pnnaWITHTWODERIV{\pn}(V_L')(-L^{2n-2})L^m\qmn{\qm}\pnna{\pn}\\
& +\pnna{\pn}L^m\qmn{\qm}\pnnaWITHTWODERIV{\pn}(V_L')(-L^{2n-2}))\\
=L^{m+3n-\epsilon}&\bndd
\end{align}
Since $m+3n-\epsilon<1+2n$, the expectation value of $I_4$ with respect to $\chia u$ is $\ErrorTermsn$. 

Finally the terms $I_5$ and $I_6$ involving $R_3=L^{3n-3}\bndd$ are considered. These are each the sum of two inner products, but by the symmetry and anti symmetry properties of the operators involved, it is sufficient to consider just one of the inner products for each. The expectation value of $I_5$ is considered, and the operators $\chia$, $\pnna{\pn}$ and $\dr$ are moved to the left side of the inner product. 
\begin{align}
\langle  u,\chia \pnna{\pn} L^{2+n-\epsilon}\dr g_{L^m} R_3 \chia u\rangle
=& \langle -\dr\pnna{\pn}\chia u,L^{2+n-\epsilon}g_{L^m}R_3\chia u\rangle
\end{align}
The Cauchy-Schwartz inequality is applied and the resulting norms are found to be of a bounded operator acting on $\dr\chia u$ and on $L^{4n-1-\epsilon}\chia u$. By corollary \ref{cII6.1}, $\|\dr\chia u\|^2$ is time integrable and hence $\Oint$ and $\ErrorTermsn$. Since $4n-1-\epsilon<\frac{1+2n}{2}$, $\|L^{4n-1-\epsilon}\chia u\|^2$ is $\ErrorTermsn$. 
\begin{align}
\langle  u,\chia \pnna{\pn} L^{2+n-\epsilon}\dr g_{L^m} R_3 \chia u\rangle
\geq&-(\|\pnna{\pn} \dr\chia u\| \|L^{2+n-\epsilon}L^{3n-3}\bndd\chia u\|)\\
\geq&-(\|\pnna{\pn} \dr\chia u\|^2+ \|L^{2+n-\epsilon}L^{3n-3}\bndd\chia u\|^2)\\
\geq& -(\|\dr\chia u\|^2+ \|L^{4n-1-\epsilon}\chia u\|^2)\\
\geq& -\ErrorTermsn
\end{align}

The term $I_6$ is estimated in expectation value by expanding $g_{L^m}'$ and finding the commutator is a bounded operator times $L^{4n+m-1-\epsilon}$. Since $4n+m-1-\epsilon<1+2n$, this expectation value in $\ErrorTermsn$. 
\begin{align}
\langle u, \chia L^{2+m+n-\epsilon}g_{L^m}' R_3 \chia u\rangle
=& \langle u, \chia L^{2+n-\epsilon}\qmn{\qm} R_3 \chia u\rangle\\
=& \langle u,\chia L^{2+m+n-\epsilon} L^{3n-3}\bndd \chia u\rangle\\
=&  \langle u,\chia L^{4n+m-1-\epsilon} \bndd\chia u\rangle \\
=& \Oint
\end{align}
This term is $\Oint$ because $4n+m-1-\epsilon\leq\frac32$. This shows the contributions from $I_2, \ldots, I_6$ are lower order. 

Combining the results so far, we have that 
\begin{align}
\langle u,[H_1+H_3,\Gammanm]u\rangle
\geq& 
\langle u,\chia\pnna{\pn}L^{n-\epsilon}[-\drr,\gamma_{L^m}]\pnna{\pn}\chia u\rangle \\
&+\langle u,\chia\pnna{\pn} L^{n-\epsilon} \epsilon (-V_L'  g_{L^m}) \pnna{\pn} \chia u \rangle + \ErrorTermsn  \\
\geq&
\langle \pnia\chia u, L^{n-\epsilon}[H_1+\epsilon H_3,\gamma_{L^m}]\pnna\chia u\rangle + \ErrorTermsn .
\label{eNearTheEndOfPSH1EpsilongVL}
\end{align} 

From the angular modulation commutator estimate, lemma \ref{lAngModResultLikePreviouslII4.2}, we have the desired result, 
\begin{align}
\langle u,[H_1+H_3,\Gammanm]u\rangle
\geq& c \| L^{\frac{n-\epsilon}{2}} \qmnt{\qm}^\frac12 \dr \pnna{\pn}\chia u\|^2 \\
&+ c\epsilon \langle \pnna{\pn}\chia u, L^{2+n-\epsilon}  g_{L^m} V_L') \pnna{\pn}\chia + \ErrorTermsn. 
\end{align}
Note that lemma \ref{lAngModResultLikePreviouslII4.2} does not require the full $-g_{L^m} V_L'$ term, since there's a $-(1-\epsilon)g_{L^m} V_L'$ term appearing on the right. Thus, even though we only have $- \epsilon g_{L^m} V_L'$ in equation \eqref{eNearTheEndOfPSH1EpsilongVL}, this is sufficient to get positivity and control over $-(\epsilon/2) g_{L^m} V_L'$. 
\end{proof}

\subsection{Derivative Bounds}
\label{ssDerivativeLocalisedEstimates}

%
%

In this subsection, we remove the localisation on the $\qm$ phase space variable. The resulting estimates are localised in $\pn$ only. Our estimates will be in the regions corresponding to $\pnnNOARG$ and $\pniNOARG$. The first lemma shows that $\Gammanhalf$ majorates $L^\frac{\frac32+m-\epsilon}{2}$, and the second, that $\Gammann$ majorates $L^{1-\delta-\frac\epsilon2}$. 

The first lemma is a quick application of lemma \ref{lfPSH1}, with lemma \ref{lfPSH2} used to show the $H_2$ contributions are lower order. 

\begin{lemma}
\label{lEstimatingThePhaseSpaceCommutatorBeneathStrip}
For \uasoluinS, $\epsilon>0$, $n\in[0,\frac12]$, and $\chi$ compactly supported, there is a constant $C$ such that, 
\begin{align}
\langle u,[\linH,\Gammanhalf]u\rangle\geq& C\|L^{\frac{\frac32+n-\epsilon}{2}}\pnn{\pn}\chi u\|^2-\ErrorTermsn
\end{align}
\end{lemma}
\begin{proof}
By the results of the angular modulation argument, we have that for any smooth, compactly supported function, $\chi$, 
\begin{align}
-\dr \qmn{\qm}\dr - L^2 g_{L^m}V_L' \geq& L^{3m} \chi
\end{align}
Taking $\tilde\chi$ to be smooth, positive, compactly supported, and identically one on the support of $\chia$ and $m=\frac12$, we have, from the result of lemmas \ref{lfPSH1} and \ref{lfPSH2}, and lemma \ref{lAngModResultLikePreviouslII4.2}, 
\begin{align}
\langle u,[\linH,\Gammanhalf]u\rangle
\geq& \langle \pnna{\pn}\chia u, L^{n-\epsilon} (-\dr \qmnt{\qhalf}\dr - \epsilon L^2  g_{L^m} V_L')\pnna{\pn} \chia u\rangle + \ErrorTermsn \\
\geq& \epsilon\langle \pnna{\pn}\chia u, L^{\frac32 +n-\epsilon} \tilde\chi \pnna{\pn} \chia u\rangle + \ErrorTermsn \\
\geq& C\| L^\frac{3/2+n-\epsilon}{2} \pnna{\pn} \chia u \|^2 \\
&+ C\| L^\frac{3/2+n-\epsilon}{2} L^{n-1} \bndd \chia u \|^2 + \ErrorTermsn .
\end{align}
Since the second term, from $[\tilde\chi,\pnna{\pn}]$, is also a lower order term, the desired result holds. 
\end{proof}

%
%

Our next goal is a result with localisation $\pni{\pn}$. We begin by introducing new $\pnifn$ type localisations, which are used as smooth approximations to $\pniNOARG$. 

\begin{definition}
The function $\pnibNOARG:[0,\infty)\rightarrow[0,1]$ is defined to be a smooth $C^\infty$  function which has support on $[\frac12,\infty)$ and is identically one on $[1,\infty)$. The function $\pnicNOARG:[0,\infty)\rightarrow[0,1]$ is defined to be a $C^\infty$ function which has support on $[0,2]$ and is identically $1$ on $[0,1]$. These are extended as even functions. The extensions are also Schwartz class, since the original functions are constant in a neighbourhood of zero. The operator $\pnia{\pn}$ is defined by
\begin{align}
\pnia{\pn} = \pnib{\pn} \pnic{\pn} 
\end{align}
\end{definition}

The following lemma allows us to replace one $\pnnfn$ type localisation with $\pnifn$ type localisation and to move this replacement through $\qm$ localisation. 

\begin{lemma}
\label{Aconvenientrugunderwhichtosweepalltheproblemswithfunctionsoftwooperators}
Given $\delta$ and $\epsilon$ positive, there is a constant $C_1$ such that for all $v\in\solset$ and for $F_2$ either $\pnnaNOARG$ or $\pnncNOARG$, 
\begin{align}
\| \pni{\pn} v\|\leq& \|\pnia{\pn} v\|\\
C_1 \| L^{-\delta}\pniaOverF{\pn} v\| \leq& \|\pn v\|
\end{align}

In addition, for $n\in[0,\frac12]$ and $m\in[0,\frac12]$, there is a constant $C_2$, such that for any differentiable function $F_1$ and for $F_2$ either $\pnnaNOARG$ or $\pnncNOARG$, 
\begin{align}
\| L^{1-m-n-\delta} [F_1(\qm),\pniaOverF{\pn}] \|_\opnorm \leq C_2\|F_1'\|_\infty 
\end{align}
\end{lemma}
\begin{proof}
Since $L$, $\dr$, and $\pn$ are all commuting operators, any functions of these operators, defined by the spectral theorem, also commute. The spherical harmonic decomposition is into orthogonal subspaces preserved by these operators. Therefore, it is sufficient to prove the first two results for functions with a single spherical harmonic component and to prove the third result by considering the operator on the right as an operator on a single spherical harmonic. 

For a fixed value of $\tl$, 
\begin{align}
\pniaWITHl{x}=&\pnibWITHl{|x|} \pnic{|x|}\\
\geq& \chi([1,\infty),\tl^{\delta}|x|)\chi([0,1],x)\\
\geq& \chi([\tl^{-\delta},\infty),|x|)\chi([0,1],x)\\
\geq& \chi([\tl^{-\delta},1],|x|)\\
\geq& \pniWITHl{x}
\end{align}
Therefore, by the spectral theorem, for any function $v\in\solset$, on each spherical harmonic
\begin{align}
\| \pni{\pn} v\| \leq \|\pnia{\pn} v\|
\end{align}

Now the case when $F_2$ is either $\pnnaNOARG$ or $\pnncNOARG$ is considered. Both $\pnna{x}$ and $\pnnc{x}$ are smooth, strictly positive functions, so in either case $F_2(x)$ has a bounded inverse for $|x|\in\supp(\pnicNOARG)\subset[0,2]$, and $\pnicOverF{\pn}$ is a bounded operator. Since all the operators involved commute and $\pnibWITHl{|x|}$ is bounded, $\pniaOverF{\pn}$ is a well defined, bounded operator. Again, for a fixed $l$,
\begin{align}
|x| \geq& \frac12 \tl^{-\delta} \chi([\frac12 \tl^{-\delta},\infty),|x|) \\
\geq& \frac12 \tl^{-\delta} \chi([\frac12,\infty),\tl^{\delta}|x|) \\
\geq& \frac12 \tl^{-\delta} \pnibWITHl{|x|}\\
\geq& C\tl^{-\delta} \pnibWITHl{|x|} \pnicOverF{x} \\
\geq& C\tl^{-\delta} \pniaOverF{x}
\end{align}
On each spherical harmonic, by the spectral theorem, 
\begin{align}
\| \pn v\|= \| (|\pn|) v\| \geq C\| L^{-\delta} \pniaOverF{\pn} v \|
\end{align}

Finally the commutator is calculated. Again this is proven on a single spherical harmonic shell. 
\begin{align}
\pniaOverF{\pn}
=& \pnib{\pn} \pnicOverF{\pn} \\
L^{1-m-n-\delta} [F_1(\qm),\pniaOverF{\pn}]
=& L^{1-m-n-\delta} [F_1(\qm), \pnib{\pn} \pnicOverF{\pn}]\\
=&L^{1-m-n-\delta}[F_1(\qm), \pnib{\pn}] \pnicOverF{\pn} \\
&+ L^{1-m-n-\delta} \pnib{\qm} [F_1(\pn), \pnicOverF{\pn}]
\end{align}
Each of these terms is now estimated. Since only one spherical harmonic is being considered, the operator $L$ can be replaced with the constant $l$. This simplifies the discussion of the commutator involving $\pnib{\pn}$. 

Since $\pnibNOARG$ is $C^\infty$ and constant outside a compact interval, it follows that $\pnibNOARG'$ is Schwartz  class and that $\|\Fourier{\pnibNOARG'}\|_1\leq\infty$. 
\begin{align}
\| L^{1-m-n-\delta}[F_1(\qm),\pnibWITHl{\pn}]\pnicOverF{\pn}\|_\opnorm
\leq& Cl^{-\delta} \| l^{1-2n}[F_1(\qn),\pnibWITHl{\pn}]\|_\opnorm\\
\leq& C l^{-\delta} \|F_1'\|_\infty \|\Fourier{ (\pnibWITHl{|\bullet|})'}\|_1 
\end{align}
At this point, the derivative of the function in a scaled variable, $(\pnibWITHl{\bullet})'$, is evaluated to be a scaled version of the derivative evaluated at the scaled variable, $l^\delta \pnibNOARG'(l^\delta\bullet)$. To evaluate the Fourier transform of this, it is noted that $\|\Fourier{f(\lambda\bullet)}\|_1=\|\Fourier{f(\bullet)}\|_1$. 
\begin{align}
\| L^{1-m-n-\delta}[F_1(\qm),\pnib{\pn}] \pnicOverF{\pn}\|_\opnorm
\leq& C l^{-\delta} \|F_1'\|_\infty \|l^\delta \Fourier{ \pnibDERIVWITHl{|\bullet|}}\|_1 \\
\leq& C \|F_1'\|_\infty \| \Fourier{\pnibDERIVNOARGNOSCALE} \|_1 \\
\leq& C \|F_1'\|_\infty
\end{align}
This completes the estimate on the first of the commutator terms. 
 
Since $\pnicNOARG$ is Schwartz class, and $F_2$ is smooth and has bounded inverse on $\supp(\pnicNOARG)$, $\pnicOverFNOARG$ is Schwartz class, and $\|\Fourier{\pnicOverFNOARG'}\|_1\leq\infty$.
\begin{align} 
\|L^{1-m-n-\delta}\pnib{\pn}[F_1(\qm),\pnicOverF{\pn})\|_\opnorm 
\leq& C\|L^{1-m-n-\delta}[F_1(\qn),\pnicOverF{\pn}]\|_\opnorm\\
\leq& C\|L^{1-m-n}[F_1(\qm),\pnicOverF{\pn}\|_\opnorm\\
\leq& C\|F_1'\|_\infty \|\Fourier{\pnicOverFNOARG'}\|_1 \\
\leq& C\|F_1'\|_\infty
\end{align}
This completes the estimate on the second of the commutator terms. 
\end{proof}

We now prove a $\pniNOARG$ localised estimate. This is analogous to the $\pnnNOARG$ localised estimate in lemma \ref{lEstimatingThePhaseSpaceCommutatorBeneathStrip}; although, the proof is more complicated because the $\pnifn$ type localisations are more complicated to work with. Because $\pn$ is bounded below on the support of $\pni{\pn}$, we can dominate more powers of $L$ with this localisation. Here, we take $m=n$. 

\begin{lemma}
\label{lEstimatingThePhaseSpaceCommutatorInStrip}
For \uasoluinS, $\epsilon>0$, $\delta>0$, and $n\in[0,\frac12]$,
\begin{align}
\langle u,[\linH,\Gammann]u\rangle\geq& \|L^{1-\frac\epsilon2-\delta}\pni{\pn}\chia u\|^2 -\ErrorTermsn
\end{align}
\end{lemma}
\begin{proof}
We take the result of lemmas \ref{lfPSH1} and \ref{lfPSH2} with $m=n$, 
\begin{align}
\langle u,[\linH,\Gammann]u\rangle
\geq& \langle \pnna{\pn}\chia u, L^{n-\epsilon}(-\dr \qmnt{\qn}\dr + \epsilon L^2(- g_{L^n}V_L')) \pnna{\pn}\chia u\rangle + \ErrorTermsn
\end{align}

The term involving the derivatives is estimated first. The first step is weakening the estimate to eliminate the square root. 
\begin{align}
\| L^\frac{2n-\epsilon}{2} \qmnt{\qn}^\frac12 \dr\pnna{\pn}\chia u\| 
\geq&\| L^{n-\frac\epsilon2} \qmnt{\qn} \dr \pnna{\pn}\chia u\|
\end{align}
Now $\dr$ is treated as if it were a localisation in $\dr$ and is commuted to the right of the chain of operators.
\begin{align}
\| L^\frac{2n-\epsilon}{2} \qmnt{\qn} \dr\pnna{\pn}\chia u\| 
\geq&\| L^{n-\frac\epsilon2} \dr \qmnt{\qn} \pnna{\pn}\chia u\| \\
&- \|L^{n-\frac\epsilon2}[\qmnt{\qn},\dr] \pnna{\pn}\chia u\| 
\end{align}
The commutator from moving $\dr$ through $\qmnt{\qn}$ can be computed explicitly. It involves the function $\qmntWITHDERIV{x}=-2 x\qmnt{x}^2$, which is bounded. 
\begin{align}
[\qmnt{\qn},\dr]
=&\dr \left(1+(L^n\rho_*)^2\right)^{-1}\\
=&-2 L^n \qn\qmnt{\qn}^2\\
=& L^n \bndd .
\end{align}
The ``localisation'' $L^n\dr$ can be replaced by $L^{1-\delta} \pniaOverpnna{\pn}$ using lemma \ref{Aconvenientrugunderwhichtosweepalltheproblemswithfunctionsoftwooperators}. 
\begin{align}
\| L^\frac{2n-\epsilon}{2} \qmnt{\qn} \dr\pnna{\pn}\chia u\| 
\geq& C\|L^{1-\frac\epsilon2-\delta} \pniaOverpnna{\pn} \qmnt{\qn} \pnna{\pn}\chia u\| \\
&-  \|L^{2n-\frac\epsilon2} \bndd \chia u\|
\end{align}

The $\pniaOverpnna{\pn}$ is commuted back through the $\qn$ localisation. The commutator involving $\pniaOverpnna{\pn}$ is shown to be lower order by lemma \ref{Aconvenientrugunderwhichtosweepalltheproblemswithfunctionsoftwooperators}. Lemma \ref{lSomeOtherFunctionsInDRCommutatorsAreGood} can not be directly applied to commutators involving $\pniaOverpnna{\pn}$, because $\pnia{\pn}=\pnib{\pn}\pnic{\pn}$ is a product of two operators with different localisations. 
\begin{align}
\|L^\frac{2n-\epsilon}{2} \qmnt{\qn} \dr\pnna{\pn}\chia u\|
\geq& C\|L^{1-\frac\epsilon2-\delta} \qmnt{\qn} \pnia{\pn}\chia u\|\\
&- \|L^{1-\frac\epsilon2-\delta}[\pniaOverpnna{\pn},\qmnt{\qn}] \pnna{\pn}\chia u\|\\
&- \|L^{2n-\frac\epsilon2} \bndd \chia u\| \\
\geq& C\|L^{1-\frac\epsilon2-\delta} \qmnt{\qn} \pnia{\pn}\chia u\| \\
&- C\|L^{1-\frac\epsilon2-\delta}L^{2n-1+\delta} \bndd\chia u\|
- \|L^{2n-\frac\epsilon2}\bndd\chia u\|
\end{align}
Since $2n-\frac\epsilon2$ is the power of $L$ appearing in both error terms, and this exponent is less than $\frac{1+2n}{2}$, the error terms are $\ErrorTermsn$. To simplify the summation of the current estimate with the next step in this proof, an additional localisation by $(-V_L' \rho_*^{-1})$ is introduced. 
\begin{align}
&\| L^\frac{2n-\epsilon}{2} \qmnt{\qn} \dr\pnna{\pn}\chia u\|\\
\geq& C\|L^{1-\frac\epsilon2-\delta} \qmnt{\qn} (-V_L\rho_*^{-1})^\frac12 \pnia{\pn}\chia u\|-\ErrorTermsn 
\label{EstimateInStripFromH1}
\end{align}

Now the $ g_{L^n}V_L'$ terms are estimated. The goal is again to replace the localisation $\pnna{\pn}$ by $\pnia{\pn}$. This requires commuting $\pniaOverpnna{\pn}$ through the $\rho_*$ localisation. 
A factor of $L^\delta$ is dropped to control commutator terms involving $\pnia{\pn}$ at a later stage. 
\begin{align}
\|L^{\frac{2-\epsilon}{2}}(-V_L'\rho_*^{-1})^\frac12 \qmf{\qn} \pnna{\pn}\chia u\|
=&\|L^{1-\frac\epsilon2-\delta} \qmf{\qn} (-V_L'\rho_*^{-1})^\frac12 \pnna{\pn}\chia u\|
\end{align}
The bounded localisation $\pniaOverpnna{\pn}$ is introduced, and then commuted through the localisation in $\qn$. 
\begin{align}
\|L^{\frac{2-\epsilon}{2}}&(-V_L' \rho_*^{-1})^\frac12 \qmf{\qn} \pnna{\pn}\chia u\|\\
\geq&C\| \pniaOverpnna{\pn} L^{1-\frac\epsilon2-\delta} \qmf{\qn}(-V_L'\rho_*^{-1})^\frac12 \pnna{\pn}\chia u\|\\
\geq&C_1 \|L^{1-\frac\epsilon2-\delta} \qmf{\qn} \pniaOverpnna{\pn} (-V_L'\rho_*^{-1})^\frac12 \pnna{\pn}\chia u\| \\
& -C_2\|L^{1-\frac\epsilon2-\delta}[\pniaOverpnna{\pn},\qmf{\qn}] (-V_L'\rho_*^{-1})^\frac12 \pnna{\pn}\chia u\|
\end{align}
Now the operator $\pniaOverpnna{\pn}$ is commuted through the localisation $(-V_L'\rho_*^{-1})$. This eliminates the $\pnna{\pn}$ localisation. 
\begin{align}
\|L^{\frac{2-\epsilon}{2}}&(-V_L'\rho_*^{-1})^\frac12 \qmf{\qn} \pnna{\pn}\chia u\|\\
\geq&C \|L^{1-\frac\epsilon2-\delta} \qmf{\qn} (-V_L'\rho_*^{-1})^\frac12 \pnia{\pn}\chia u\| \\
& -C \|L^{1-\frac\epsilon2-\delta} \qmf{\qn} [\pniaOverpnna{\pn}, (-V_L'\rho_*^{-1})^\frac12] \pnna{\pn}\chia u\| \\
& -\|L^{1-\frac\epsilon2-\delta}[\pniaOverpnna{\pn}, \qmf{\qn}] (-V_L'\rho_*^{-1})^\frac12 \pnna{\pn}\chia u\|
\end{align}
The commutator terms can now be estimated using lemma \ref{Aconvenientrugunderwhichtosweepalltheproblemswithfunctionsoftwooperators} and are found to be lower order error terms. Note that $\qmf{x}$ has a $L^\infty$ derivative. 
\begin{align}
\|L^{\frac{2-\epsilon}{2}}& (-V_L'\rho_*^{-1})^\frac12 \qmf{\qn} \pnna{\pn}\chia u\|\\
\geq&C \|L^{1-\frac\epsilon2-\delta} \qmf{\qn} (-V_L'\rho_*^{-1})^\frac12\pnia{\pn} \chia u\| \\
& -C \|L^{1-\frac\epsilon2-\delta}\qmf{\qn} L^{n+\delta-1}\bndd_1 \pnna{\pn}\chia u\| \\
& -\|L^{1-\frac\epsilon2-\delta}L^{2n+\delta-1}\bndd_2(-V_L'\rho_*^{-1})^\frac12 \pnna{\pn}\chia u\|\\
\geq&C
\|L^{1-\frac\epsilon2-\delta}\qmf{\qn} (-V_L'\rho_*^{-1})^\frac12 \pnia{\pn}\chia u\|\\
 &+\ErrorTermsn
\label{EstimateInStripFromH3}
\end{align}

Equations \eqref{EstimateInStripFromH1} and \eqref{EstimateInStripFromH3}
can be combined with the initial estimate on the expectation value of $[\linH,\Gammann]$ to produce an intermediate result. 
\begin{align}
\langle u,[\linH,\Gammann]u\rangle
\geq& C
(\|L^{1-\frac\epsilon2-\delta}\qmf{\qn} (-V_L'\rho_*^{-1})^\frac12 \pnia{\pn}\chia
u\|^2\\
&+\|L^{1-\frac\epsilon2-\delta}\qmnt{\qn} (-V_L'\rho_*^{-1})^\frac12 \pnia{\pn}\chia u\|^2)\\
& -\ErrorTermsn
\end{align}

Both these terms involve localisation in $\qn$ acting on $(-V_L'\rho_*^{-1})^\frac12 \pnia{\pn}\chia u$. The sum of the localisations is bounded below by a constant. 
\begin{align}
\qmnt{x}^2 + \qmf{x}^2 
= (\frac{1}{1+x^2})^2+(x\sqrt{\frac{g(x)}{x}})^2
\geq C
\end{align}
Therefore, the two terms can be combined to provide a better estimate.  
\begin{align}
\langle u,[\linH,\Gammann]u\rangle
\geq & C \|L^{1-\frac\epsilon2-\delta} (-V_L'\rho_*^{-1})^\frac12 \pnia{\pn}\chia u\|- \ErrorTermsn
\end{align}

To eliminate the factor of $(-V_L'\rho_*^{-1})^\frac12 $, a new function is introduced. This function, $f$, is a smooth, compactly supported function and equal to the inverse of $(-V_L'\rho_*^{-1})^\frac12 $ on $\supp{\chia}$. Since it is bounded, it can be freely introduced into the norm. The $\rho_*$ localisation can then be commuted through $\pnia{\pn}$. 
\begin{align}
\langle u,[\linH,\Gammann]u\rangle
\geq & C \|L^{1-\frac\epsilon2-\delta} f (-V_L'\rho_*^{-1})^\frac12 \pnia{\pn}\chia u\|- \ErrorTermsn\\
\geq& C \|L^{1-\frac\epsilon2-\delta} \pnia{\pn} f (-V_L'\rho_*^{-1})^\frac12 \chia u\|\\
& -\|L^{1-\frac\epsilon2-\delta} [f (-V_L'\rho_*^{-1})^\frac12 ,\pnia{\pn}]\chia u\|\\
& - \ErrorTermsn
\end{align}
From the definition of $f$, the product of all the $\rho_*$ localisation reduces to $\chia$. The new commutator terms can be estimated by lemma \ref{Aconvenientrugunderwhichtosweepalltheproblemswithfunctionsoftwooperators} and are found to be lower order error terms.  
\begin{align}
\langle u,[\linH,\Gammann]u\rangle
\geq& C \|L^{1-\frac\epsilon2-\delta} \pnia{\pn} \chia u\|-\|L^{1-\frac\epsilon2-\delta} L^{n+\delta-1}\bndd \chia u\|\\
&- \ErrorTermsn\\
\geq& C \|L^{1-\frac\epsilon2-\delta} \pnia{\pn} \chia u\|- \ErrorTermsn
\end{align}
Lemma \ref{Aconvenientrugunderwhichtosweepalltheproblemswithfunctionsoftwooperators} can now be used to replace the localisation in $\pnia{\pn}$ by localisation in $\pni{\pn}$. 
\begin{align}
\langle u,[\linH,\Gammann]u\rangle
\geq& C \|L^{1-\frac\epsilon2-\delta} \pni{\pn} \chia u\|- \ErrorTermsn
\label{PartIOflEstimatingThePhaseSpaceCommutator}
\end{align}
\end{proof}

%
%
\subsection{Phase Space Induction}
\label{ssPhaseSpaceInduction}

The previous section shows that $\Gammanhalf$ and $\Gammann$ majorate $L^\frac{\frac32+n-\epsilon}{2}$ and $L^{1-\delta-\frac{\epsilon}{2}}$ respectively. We would like to integrate the Heisenberg identity to conclude that the time integral of the expectation value of these powers of $L$ are bounded. However, our definition of majoration allows the domination to occur only in a region of phase space and the lower order terms to be unlocalised. 

The lower order terms are $\ErrorTermsn$. To control these terms, we use a finite induction on $n$, to eventually control $L^{1-\varepsilon}$ without phase space localisation. 

Following this, we use the control of $L^{1-\varepsilon}$ in the conformal estimate to prove point wise in time, weighted $L^6$ decay, with bounds involving an additional $L^\varepsilon$ factor. 

\begin{theorem}[Phase Space Induction]
\label{tPhaseSpaceInduction}
If $\varepsilon>0$, then \foruasoluinS, there is a constant, $K_{\text{PS}}$, such that
\begin{align}
\int_1^\infty \|L^{1-\frac\varepsilon2} \chia u(t)\|^2 dt \leq K_{\text{PS}}( E[u,\dot{u}] + \|u(1)\|^2)
\end{align}
\end{theorem}
\begin{proof}
This is proven for $\varepsilon= \delta+\epsilon$. 

We induct on $n$ to prove, simultaneously, the two statements
\begin{align}
\|L^{\frac{\frac32+n-\delta-\epsilon}{2}} \chia u\|^2=&\Oint ,
\label{PhaseSpaceInductionR2} \\
\|L^{1-\delta}\chi([L^{1-n},\infty),\dr)\chia u\|^2=&\Oint .
\label{PhaseSpaceInductionR1}
\end{align}
The induction will run from $n=0$, in steps of size $\delta$, as long as
\begin{align}
n\leq \frac12 -2\delta-\epsilon .
\label{PSInductionCondition}
\end{align}
Each step in the induction will be proven by Morawetz type arguments using $\Gammann$ and $\Gammanhalf$.

The base case of \eqref{PhaseSpaceInductionR2} follows from the angular modulation result, theorem \ref{tII5.1}, which says, in the notation of this section $\| L^\frac32 \chia u\|^2 = \Oint$. The base case of \eqref{PhaseSpaceInductionR1} follows from the angular modulation theorem, theorem \ref{tII5.1}. For $n=0$, since 
\begin{align}
l^{1-\delta} \chi([l,\infty),x) 
\leq & l^{1-\delta} \chi([l^{1-\delta},\infty),x) 
\leq x ,
\end{align}
by the spectral theorem and corollary \ref{cII6.1}, 
\begin{align}
\| L^{1-\delta} \chi([L,\infty), \dr) \chia u\|^2
\leq& \| \dr \chia u \|^2 
=\Oint .
\end{align}

The inductive step is now considered. 


Morawetz type estimates with $\Gammann$ and $\Gammanhalf$ will be used. Before these estimates are proven, $\|\Gammanm u\|^2$ is shown to be bounded by the energy and a local decay term, under the condition $0\leq n \leq m\leq \frac12$. 

To begin estimating the norm of $\Gammanm u$, $\Gammanm$ is expanded as a product of operators, and the factor of $\gamma_{L^m}$ is expanded as a sum of two terms. These terms can be further simplified by eliminating bounded functions. 
\begin{align}
\|\Gammanm u\|=&\| \chia\pnna{\pn}\gamma_{L^m} L^{n-\epsilon}\pnna{\pn}\chia u\|\\
\leq& \|\chia\pnna{\pn} \dr g_{L^m}L^{n-\epsilon} \pnna{\pn}\chia u\|\\
&+\|\chia\pnna{\pn}\frac12 \frac{L^m}{1+(\frac{\rho_*L^m}{2M})^2}L^{n-\epsilon}\pnna{\pn}\chia u\|\\
\leq& \|\pnna{\pn} \dr g_{L^m}L^{n-\epsilon}\pnna{\pn}\chia u\|+\|L^{m+n-\epsilon}\chia u\|
\end{align}
Interpolation can be used to simplify $\|L^{m+n-\epsilon}\chia u\|$. 
\begin{align}
\|\Gammanm u\| \leq& \|\pnna{\pn} \dr g_{L^m}L^{n-\epsilon}\pnna{\pn}\chia u\|+C(E[u]+\|\weakloc^{-1}u\|^2)
\end{align}
To control $\|\pnna{\pn} \dr g_{L^m}L^{n-\epsilon}\pnna{\pn}\chia u\|$, it is rewritten in terms of $\pn$ and $L$, and its sub factors are rearranged. 
\begin{align}
\|\pnna{\pn}\dr L^{n-\epsilon}g_{L^m}\pnna{\pn}\chia u\|
=&\|L^{1-\epsilon}\pnna{\pn}\pn g_{L^m}\pnna{\pn}\chia u\|\\
\leq& \|L^{1-\epsilon} g_{L^m}\pn\pnna{\pn}^2\chia u\|\\
&+\|L^{1-\epsilon}[\pn\pnna{\pn},g_{L^m}]\pnna{\pn}\chia\|
\end{align}
The first of these two norms can be estimated by lemma \ref{lIIqLcontrolsGamma}. The second can be estimated using lemma \ref{lSomeOtherFunctionsInDRCommutatorsAreGood} which states that the commutator is $L^{n+m-1}\bndd$. 
\begin{align}
\|\pnna{\pn}\dr L^{n-\epsilon}g_{L^m}\pnna{\pn}\chia u\|
\leq& \|L^{ \frac{1+n\epsilon-2\epsilon}{1-\epsilon}} \chia u\| +
\|\dr \chia u\|\\
&+\|L^{n+m-\epsilon}\bndd\pnna{\pn}\chia u\|
\end{align}
Since $\frac{1+n\epsilon-2\epsilon}{1-\epsilon}\leq 1$ and $n+m-\epsilon\leq1$, each of the terms involving powers of $L$ in this expression can be estimated by interpolation. 
\begin{align}
\|\pnna{\pn}\dr L^{n-\epsilon}g_{L^m}\pnna{\pn}\chia u\|\leq&C(E[u]+\|\weakloc^{-1}u\|^2)
\end{align}
This completes the estimate on $\|\Gammanm u\|$. 

We can now integrate the Heisenberg like identity to find
\begin{align}
\int_1^T \langle u, [\linH, \Gammanm] u\rangle dt
=& \int \Dt (\langle u,\Gammanm \dot{u}\rangle - \langle \dot{u}, \Gammanm u\rangle) dt\\
=& -2 \langle \dot{u}, \Gammanm u\rangle |_1^T \\
\leq& C(E[u]+\|\weakloc^{-1}u(T)\|^2 +\|\weakloc^{-1}u(1)\|^2) 
\end{align}
The second term on the right can be dropped since, by the local decay result, the norm $\|\weakloc^{-1}u(T)\|^2 \rightarrow 0$ on a sequence of times. 
\begin{align}
\int_1^T \langle u, [\linH, \Gammanm] u\rangle dt
\leq& C(E[u] +\|\weakloc^{-1}u(1)\|^2) 
\end{align}

The derivative localisation results, lemmas \ref{lEstimatingThePhaseSpaceCommutatorInStrip} and \ref{lEstimatingThePhaseSpaceCommutatorBeneathStrip} can now be applied. From the inductive hypothesis \eqref{PhaseSpaceInductionR2} and condition \eqref{PSInductionCondition}, the lower order terms are integrable:
\begin{align}
\| L^{\frac{1+2n}{2}} \chia u \|^2 =& \Oint ,\\
\ErrorTermsn =& \Oint .
\end{align}

From inductive hypothesis \eqref{PhaseSpaceInductionR1} and the second derivative localising result, lemma \ref{lEstimatingThePhaseSpaceCommutatorInStrip}, 
\begin{align}
\|L^{1-\delta}\chi([L^{1-n-\delta},\infty),\dr)\chia u\|^2
=& \|L^{1-\delta}\chi([L^{1-n-\delta},L^{1-n}),\dr)\chia u\|^2\\
&  +\|L^{1-\delta}\chi([L^{1-n},\infty),\dr)\chia u\|^2 \\
=& \|L^{1-\delta}\pni{\pn}\chia u\|^2\\
&  +\|L^{1-\delta}\chi([L^{1-n},\infty),\dr)\chia u\|^2 \\
&\leq \langle u,[\linH,\Gammann]u\rangle +\ErrorTermsn\\
&+\Oint . 
\end{align}
Integration in time extends inductive hypothesis \eqref{PhaseSpaceInductionR1} to $n+\delta$
\begin{align}
\int_1^\infty \|L^{1-\delta}\chi([L^{1-n-\delta},\infty),\dr)\chia
  u\|^2 dt
\leq& C(E[u] + \|u(1)\|^2) .
\end{align}

From condition \eqref{PSInductionCondition}, it follows that $\frac{\frac32+n-\epsilon}{2} < 1-\delta$. From the first derivative localising result, lemma \ref{lEstimatingThePhaseSpaceCommutatorBeneathStrip}, 
\begin{align}
\|L^{\frac{\frac32+n-\epsilon}{2}} \chia u\|^2
=& \|L^{\frac{\frac32+n-\epsilon}{2}} \chi([0,L^{1-n}],\dr)\chia u\|^2\\
&+\|L^{\frac{\frac32+n-\epsilon}{2}} \chi([L^{1-n},\infty),\dr) \chia u\|^2\\
\leq& \|L^{\frac{\frac32+n-\epsilon}{2}} \pnn{\pn} \chia u\|^2
  +\|L^{1-\delta} \chi([L^{1-n},\infty),\dr) \chia u\|^2\\
\leq& \langle u,[\linH,\Gammanhalf]u\rangle +\ErrorTermsn+\Oint .
\end{align}
Integration in time extends inductive hypothesis \eqref{PhaseSpaceInductionR2} up to $n+\delta$. 
\begin{align}
\label{PSInductionToExtendR2LpowerIntegration}
\int_1^\infty \|L^{\frac{\frac32+n-\epsilon}{2}} \chia u\|^2 dt
\leq& C(E[u] + \|u(1)\|^2)
\end{align}

Since the induction continues as long as condition \eqref{PSInductionCondition} holds, after the last application of equation \eqref{PSInductionToExtendR2LpowerIntegration}, 
\begin{align}
\int_1^\infty \|L^\frac{2-2\delta-2\epsilon}{2} \chia u\|^2 dt 
\leq& C(E[u] + \|u(1)\|^2)
\end{align}
This proves the desired result with $\varepsilon=\delta-\epsilon$. Since $\delta$ and $\epsilon$ can be taken to be arbitrarily small, so can $\varepsilon$. 
\end{proof}

Except for the loss of $L^\epsilon$, this provides the control on the angular energy near the photon sphere required by corollary \ref{cReductionToLonPhotonSphere} to control the weighted $L^6$ norm. Since the wave equation is linear and $L$ commutes with $\linH$, $L^\epsilon u$ is also a solution, and the energy of this function can be used to recover the additional factor of $L^\epsilon u$. 

\begin{theorem}
\label{cL6Control}
If \tildeutruesolu, and $u=r\tilde{u}$, then 
\begin{align}
\| F^\frac12 \tilde{u} \|_{L^6(\tildestarman)} 
\leq& t^{-\frac{1}{3}} C(\|u_0\|_2^2 + E[u_0,u_1] + \confchrg[u_0,u_1] + \|L^\epsilon u\|^2 + E[L^\epsilon u_0, L^\epsilon u_1])^\frac12 \\
\| (\rho_*^2+1)^\frac{-1}{2} \tilde{u} \|_{L^2(\tildestarman)}
\leq& t^{-\frac{1}{2}} C(\|u_0\|_2^2 + E[u_0,u_1] + \confchrg[u_0,u_1] + \|L^\epsilon u\|^2 + E[L^\epsilon u_0, L^\epsilon u_1])^\frac12
\end{align}
\end{theorem}
\begin{proof}
If $\tilde{u}$ is a solution to equation \eqref{tildeLW}, then $u$ is a solution to equation \eqref{LW}. The weighted $L^6$ norms are related. 
\begin{align}
\| F^\frac12 \tilde{u} \|_{L^6(\tildestarman)}^6
=& \int F^3 |\tilde{u}|^6 \dmu \\
=& \int F^3 r^{-4} |u| \dmu .
\end{align}

Since the wave equation equation \eqref{LW} is linear and $\linH$ commutes with $L$, if $u$ is a solution, then so is $L^\epsilon u$. The phase space induction theorem, theorem \ref{tPhaseSpaceInduction}, can be applied to $L^\epsilon u$. 
\begin{align}
\|L^{1-\epsilon} \chia (L^\epsilon u)\|^2 =& \Oint, \\
\int \|L \chia u\|^2 dt\leq& C(\|L^\epsilon u\|^2 + E[L^\epsilon u,L^\epsilon \dot{u}] ).
\end{align}

This is the term which is needed to control the conformal charge in the reduction to angular energy on the photon sphere, corollary \ref{cReductionToLonPhotonSphere}. 
\begin{align}
\confchrg[u(t),\dot{u}(t)] 
\leq&\confchrg[u_0,u_1]  + Ct E[u_0,u_1]^\frac12(E[u_0,u_1]^\frac12+\|u_0\|) \\
&+ Ct \int^t \int \chia |\Lu|^2 \dmu d\tau\\
\leq& \confchrg[u_0,u_1]  + Ct E[u_0,u_1]^\frac12(E[u_0,u_1]^\frac12+\|u_0\|) \\
&+ Ct (\|L^\epsilon u_0\|^2 + E[L^\epsilon u_0, L^\epsilon u_1]) 
\end{align}
Since the conformal charge controls the weighted $L^6$ norm, we have
\begin{align}
\| F^\frac12 \tilde{u} \|_{L^6(\tildestarman)}
\leq& (\int F^3 r^{-4} |u| \dmu )^\frac16 \\
\leq& Ct^\frac{-1}{3} (\|u_0\|^2 + E[u_0,u_1] + \confchrg[u_0,u_1] + \|L^\epsilon u_0\|^2 + E[L^\epsilon u_0, L^\epsilon u_1]) . 
\end{align}

Similarly, since lemma \ref{lWeightedL2ByConformal} controls the weighted $L^2$ norm by the conformal charge, 
\begin{align}
\langle \tilde{u}, \frac{1}{\rho_*^2+1} \tilde{u}\rangle_{L^2(\tildestarman)}
=&\langle u, \frac{1}{\rho_*^2+1} u\rangle\\
\leq& C t^{-2} \confchrg[u(t),\dot{u}(t)] \\
\leq& C t^{-2} \confchrg[u_0,u_1] + Ct^{-1}(E[u_0,u_1]+\|u_0\|^2+\|L^\epsilon u_0\|^2 + E[L^\epsilon u_0, L^\epsilon u_1])
\end{align}

\end{proof}

\vspace{.5in}

{\Large \bf Acknowledgements}

We wish to thank the anonymous referee and J. Sterbenz. 

These results originally appeared as part of the Ph.D. Dissertation of P. Blue, ``Decay estimates and phase space analysis for wave equations on some black hole metrics''\cite{PBDissertation}, submitted in October 2004. Both authors would like to thank the NSF for partial support under grants NSF DMS 0100490 (P. Blue) and 0501043 (A. Soffer). \providecommand{\bysame}{\leavevmode\hbox to3em{\hrulefill}\thinspace}
\providecommand{\MR}{\relax\ifhmode\unskip\space\fi MR }
\providecommand{\MRhref}[2]{%
  \href{http://www.ams.org/mathscinet-getitem?mr=#1}{#2}
}
\providecommand{\href}[2]{#2}

 \end{document}